\documentclass[leqno,12pt]{amsart}
\usepackage{amsmath,mathrsfs,mathtools}
\usepackage{amssymb, amsmath,amsfonts,amsthm}
\usepackage{txfonts}
\usepackage{subfigure}
\usepackage{amscd}
\usepackage{dsfont}
\usepackage{amssymb,supertabular,a4}

\def\prof{{\sc Proof.\ \ }}

\def\inte#1{
\displaystyle\mathop{#1\kern0pt}^\circ }

\let\t=\tau

\def\t{\tau}

\def\virgp{\raise 2pt\hbox{,}}
\def\cdotpv{\raise 2pt\hbox{;}}

\def\C{\mathop{\mathbb C\kern 0pt}\nolimits}
\def\DD{\mathop{\mathbb D\kern 0pt}\nolimits}
\def\EE{\mathop{{\mathbb E \kern 0pt}}\nolimits}
\def\K{\mathop{\mathbb K\kern 0pt}\nolimits}
\def\N{\mathop{\mathbb N\kern 0pt}\nolimits}
\def\Q{\mathop{\mathbb Q\kern 0pt}\nolimits}
\def\R{\mathop{\mathbb R\kern 0pt}\nolimits}
\def\SS{\mathop{\mathbb S\kern 0pt}\nolimits}
\def\ZZ{\mathop{\mathbb Z\kern 0pt}\nolimits}
\def\TT{\mathop{\mathbb T\kern 0pt}\nolimits}
\def\P{\mathop{\mathbb P\kern 0pt}\nolimits}

\newcommand{\beq}{\begin{equation}}
\newcommand{\eeq}{\end{equation}}
\newcommand{\ben}{\begin{eqnarray}}
\newcommand{\een}{\end{eqnarray}}
\newcommand{\beno}{\begin{eqnarray*}}
\newcommand{\eeno}{\end{eqnarray*}}
\newtheorem{defi}{Definition}[section]
\newtheorem{thm}{Theorem}[section]
\newtheorem{lem}{Lemma}[section]

\newtheorem{pro}{Proposition}[section]

\begin{document}

\title[Neutron transport equation]{\bf{Existence and BV-regularity for Neutron transport equation in non-convex domain}}
\author[Y. Guo]{Yan Guo}
\address[Y. Guo]{Division of Applied Mathematics, Brown University,
\newline\indent Providence, RI 02812,USA }
\email{Yan\_Guo@sjtu.edu.cn}

\author[X.F. Yang]{Xiongfeng Yang}
\address[X.F. Yang]{
   \newline\indent Department of Mathematics, MOE-LSC and SHL-MAC, Shanghai Jiao Tong University,
\newline\indent Shanghai, 200240, P.R. China}
\email{xf-yang@sjtu.edu.cn}

\thanks{Y. Guo's research is supported in part by NSFC grant 10828103, NSF grant DMS  1209437, Simon Research Fellowship and BICMR.
X.F. Yang is supported by NNSF of China under the Grant 11171212 and
the SJTU's SMC Projection A}

\date{}
\begin{abstract}
This paper considers the neutron transport equation in bounded domain with a combination of the diffusive boundary condition and the in-flow boundary condition. We firstly study the existence of solution in any fixed time by $L^2-L^{\infty}$ method, which was established to study Boltzmann equation in \cite{[Guo2]}. Based on the uniform estimates of the solution, we also consider the BV-regularity of the solution in non-convex domain. A cut-off function, which aims to exclude all the characteristics emanating from the grazing set $\mathfrak{S}_B$, has been constructed precisely.
\end{abstract}



 \maketitle
%\tableofcontents

\section{Introduction}
\setcounter{section}{1} \setcounter{equation}{0}

The neutron transport equation is a type of radiative transport equation, which is a balance statement that the neutrons conserved. This equation is commonly used to determine the behavior of nuclear reactor cores and experimental or industrial neutron beams. For more details, see \cite{[DL]}. In this paper, we consider the following neutron transport equation
\begin{eqnarray}\label{1origin equation}
\frac{\partial u}{\partial t}+v \cdot \nabla u + \Sigma(x,v)u = Ku + q(t,x,v).
%\\
%\label{Initial data}
%u(t,x,v)|_{t=0} = u_0(x,v), ~~~~(x,v)\in X \times V.
\end{eqnarray}
The function $u(t,x,v)$ represents a density of the number of particles. $\Sigma(x,v)\geq0 $ describes the effective total cross section, which is a given positive function of $x$ and $v$, the given operator $K$ is defined as
\begin{eqnarray}\label{1operator K}Ku(x,v)=\int_{V} k(x,v,v')u(t,x,v')dv',\end{eqnarray}
The nonnegative kernel $k(x, v, v')$ models a transfer of
a density of numbers of neutrons from one speed to another. It depends on
the state of the material at the point $x,~v,~v'$, and it is isotropic if the kernel only depends on
the variables $v$ and $v'$. $q(t,x,v)$ is a source of a finite total number of neutrons at each moment $t$.

In bounded domain, the equation describes the evolution of a population of neutrons in a domain $\Omega$  occupied by a medium which interacts with the neutrons. Here $\Omega$ be a bounded open and connected subset of $\mathbb{R}^3$, $\partial \Omega$ is denoted as its boundary. The domain $V$ is the velocity space, which generally is the form $ V=\{v\in \mathbb{R}^3| ~a \leq |v|\leq b\}$ or a finite union of spheres.

There are extensive developments in the study of the neutron transport equation. The existence of the solution for both steady neutron transport equation and time-dependent neutron transport equation have been studied in \cite{[Gibbs]}, \cite{[Pao1]}, \cite{[Pao2]} by constructing a maximal and a minimal solutions. An abstract theorem is also established in \cite{[BP]}. The existence of the solution for the neutron transport equation was constructed in the Banach spaces $L^{p}, ~1 \leq  p <\infty$, it means that $L^{\infty}$ solution has not been treated. For the asymptotic expansions in transport theory, we refer to \cite{[BSS]}, \cite{[Lar1]}, \cite{[Lar2]}, \cite{[LK]} and \cite{[PT]}. Most of the above works considered the neutron transport while the neutron flux entering $\Omega$ at each point of $\partial \Omega$ is zero, that is, the zero in-flow boundary condition. Later, the existence of the neutron transport equation with different types of boundary condition has also appeared  in \cite{[Lod1]}, \cite{[Lod2]}, \cite{[Lod3]}, \cite{[LS]} and \cite{[ZL]} and referees therein. Moreover, the authors \cite{[BLP]} had studied the existence and the asymptotic expansion of the solution for neutron transport equation with in-flow boundary condition together with a weak diffusive boundary condition by the probabilistic theory. In \cite{[WG]}, the very recent result show the asymptotic expansion of the solution for the neutron transport equation in 2-D unit disc, which gave a more precise approximation of the solution around the boundary by modifying the Milne problem.

  In the following, we list some assumptions on the phase space $\Omega \times V$. We assume that the boundary $\partial \Omega$ is locally a graph of a given $C^2$ function: for each point $x_0 \in \partial \Omega$, there exist $r>0$ and a $C^2$ function $\eta: R^2 \rightarrow R$ such that, up to a rotation and relabeling, we have
 \begin{eqnarray}
 &&\partial X \cap B(x_0; r)= \{ x\in B(x_0;r): x_3=\eta(x_1,x_2)\},\label{1Non-convex Domain boundary}\\
 && \partial X \cap B(x_0; r)= \{ x\in B(x_0;r): x_3>\eta(x_1,x_2)\}.\label{1Non-convex Domain}
 \end{eqnarray}
In this case, the outward normal direction $n$ at $x \in \partial \Omega$ can be expressed as
\begin{eqnarray*}\label{1Outer normal D}
n(x_1,x_2) = \frac{1}{\sqrt{1+|\nabla_x \eta(x_1,x_2)|^2}}\bigg(\partial_{x_1}\eta(x_1,x_2),\partial_{x_2}\eta(x_1,x_2),1\bigg).
 \end{eqnarray*}
The domain $\Omega$ is called a {\bf strictly~ non-convex ~domain} if there exists at least one point $x_0 \in \partial \Omega$ and nonzero $u\in \mathbb{R}^2$ such that (\ref{1Non-convex Domain boundary})-(\ref{1Non-convex Domain}) hold
 and
 \begin{eqnarray}\label{1Strictly non-convex}
 \sum_{i,j=1,2} u_i u_j \partial_i \partial_j \eta(x_0) <0.
\end{eqnarray}
The phase boundary in the phase space $\Omega \times V$ is denoted as $\gamma=  \partial \Omega\times V$, and we split it into the outgoing boundary $\gamma_+$, the incoming boundary $\gamma_-$, and the grazing boundary $\gamma_0$
\begin{eqnarray*}
\gamma_{\pm} &=& \{(x,v) \in \partial \Omega\times V : n(x) \cdot v \gtrless 0\},\\
\gamma_0 &=& \{(x,v) \in \partial \Omega\times V : n(x) \cdot v = 0\}.
 \end{eqnarray*}
It is known that $\gamma_+$ and $\gamma_-$ (resp $\gamma_0$) are open subsets (resp. closed) of $\gamma = \partial \Omega \times V $ such that
\begin{eqnarray}
\gamma = \partial \Omega \times V = \gamma_+ \cup \gamma_0 \cup \gamma_- .
\end{eqnarray}
In this paper, we assume that $\Omega$ is strictly non-convex domain, $V$ is a bounded domain and it can be locally expressed as (\ref{3Decomposition of velocity space}).

 Let us explain the difficulty to study the regularity of the kinetic equation in bounded domain with boundary condition. It partly dues to the characteristic nature of boundary conditions. To make it clear, we consider the transport equation with the given boundary condition
\begin{eqnarray}\label{1Test equation}
v\cdot \nabla_x f(x,v)=0,~~~f|_{\Gamma_-}=g.
\end{eqnarray}
Given $(x,v)\in X$, let $[X(s),V(s)]=[X(s;t,x,v),V(s;t,x,v)]=[x-(t-s)v,v]$ be a trajectory for the transport equation:
\begin{eqnarray}\label{1Trajectory}
\frac{d X(s)}{ds}=V(s),~~~~~~\frac{d V(s)}{ds}=0,
\end{eqnarray}
with the initial condition $[X(t;t,x,v),V(t;t,x,v)]=[x,v]$. Then, we solve (\ref{1Test equation}) as $f(x,v)=g(x_b(x,v),v)=g(x-t_b(x,v)v,v)$, where $t_{\mathbf{b}}(x,v)\geq 0$ is the backward exit time, or the last moment at which the back-time straight line
$[X(s;0,x,v),V(s;t,x,v)]$ remains in the interior of $X$. It is defined as
\begin{eqnarray}\label{1Backward exit time}
 t_{\mathbf{b}}(x,v):=\sup\big(\{0\} \cap \{\tau: x-sv \in X ~\text{for ~all}~0<s<\tau\}\big).
 \end{eqnarray}
The backward exit position on the boundary $\partial X$ is
\begin{eqnarray}\label{1Backward exit P}
x_{\mathbf{b}}(x,v)=x-t_{\mathbf{b}}(x,v)v,
\end{eqnarray}
and we always have $v\cdot n(x_{\mathbf{b}}(x,v))\leq 0$.
Similarly the forward exit time $t_{\mathbf{f}}$ and the forward exit position are defined as
 \begin{eqnarray}\label{1Forward exit time}
\begin{array}{lll}
 t_{\mathbf{f}}(x,v):=\sup\big(\{0\} \cap \{\tau: x+sv \in X ~\text{for ~all}~0<s<\tau\}\big),\\
  x_{\mathbf{f}}(x,v)=x+t_{\mathbf{f}}(x,v)v.
 \end{array}
 \end{eqnarray}
 Generally, it is difficult to determine $t_{\mathbf{b}},~x_{\mathbf{b}}$ as well as the solution to (\ref{1Test equation}) with the diffusive boundary condition. This was solved by introducing the probability measure on the boundary in \cite{[Guo2]}.

There are a few results about the regularity of the solutions to the kinetic equation in bounded domain. The first one has been appeared in  \cite{[Guo1]}, Guo constructed the singular solutions of the Vlasov-Maxwell equation on a half line. Recently, Guo \cite{[Guo2]} developed the $L^2-L^{\infty}$ estimate for the solution of Boltzmann equation in convex domain with different boundary conditions, and it was show that the solution are continuous away from the grazing set $\gamma_0$. It should be pointed out that the domain needs not to be convex for the diffusive reflection condition case. Later, the regularity of the solution for Boltzmann equation was studied in \cite{[GKTT1]}. The authors established the $C^1$ solution in convex domain and show that the solution should not be $C^2$. In the above two papers, it could be proved that $x_{\mathbf{b}}(x,v)$ has singular behavior if $n(x_{\mathbf{b}}(x,v))\cdot v =0$, and the solution might be singular on the set:
 \begin{eqnarray}\label{1singular set}
 \mathfrak{S}_B:=\{ (x,v)\in \overline{\Omega} \times V: n(x_{\mathbf{b}}(x,v))\cdot v=n(x-t_{\mathbf{b}}(x,v)v)\cdot v=0\},
 \end{eqnarray}
 which is the collection of all the characteristics emanating from the grazing set $\gamma_0$. In a non-convex domain,
 Kim \cite{[Kim]} discovered that the singularity (discontinuity) of the solution of Boltzmann equation always occurs, and such singularity propagates along the singular set $\mathfrak{S}_B$. More precisely, let the concave (singular) grazing boundary in the grazing boundary to be defined as
\begin{eqnarray}\label{1Concave grazing B}
\gamma_0^{\mathbf{S}}=\{(x,v)\in \gamma_0: ~t_{\mathbf{b}}(x,v) \neq 0~\text{and}~t_{\mathbf{b}}(x,-v)\neq 0\} \subset \gamma_0.
\end{eqnarray}
It was proved that $\gamma_0^{\mathbf{S}}$ is the only part at which discontinuity can
be created or propagates into the interior of the phase space $\Omega\times V$. So the discontinuity set of the solution in $\overline{\Omega} \times V$ is
\begin{eqnarray}\label{1Discontinuity set}
\mathfrak{D} =\gamma_0 \cup \{(x,v) \in \overline{\Omega}\times V: (x_{\mathbf{b}}(x,v),v) \in \gamma_0^{\mathbf{S}}\}.
 \end{eqnarray}
It implies that we can not get the classical solution of Boltzmann equation. A nature problem is the regularity of the solutions  in non-convex domain. Very recently, the BV-regularity of the solution to Boltzmann equation in non-convex domain has been studied in \cite{[GKTT2]}. Moreover, it was proved that the singular set to the characteristics emanating from the strictly non-convex points
  \begin{eqnarray*}\label{1singular set1}
  \big\{(x,v)\in \mathfrak{S}_B: (x_b(x,v),v) ~\text{is a strictly non-convex point} \big\}
  \end{eqnarray*}
  is a co-dimension 1 submanifold of $\Omega\times V$. This means that the BV regularity is the best regularity for Boltzmann equation in the non-convex domain.

Similarly, we expect to establish the existence and BV-regularity of the
solution for the neutron transport equation in non-convex domain. The large time behavior and the regularity of the solution to the Neutron transport equation would be considered in forthcoming papers. In this paper, we assume $\Omega$ is anon-convex domain and we consider
\begin{eqnarray}\label{1Main problem}
\frac{\partial u}{\partial t} + v\cdot \nabla u +\Sigma u =  K u + q,
\end{eqnarray}
with the initial-boundary condition
\begin{eqnarray}\label{1Boundary condition}
u(0,x,v)=u_0(x,v),\qquad u(t)|_{\gamma_-}= \mathcal{P}_{\gamma} u+ r.
\end{eqnarray}
Here $r, ~u_0$ are given functions and $\mathcal{P}_{\gamma}$ is the diffusive reflection: for $(x,v)\in \gamma_-$,
\begin{eqnarray*}\label{1Diffusive BC}
\mathcal{P}_{\gamma} u(t,x,v) = c \int_{v'\in V: n(x)\cdot v'>0} u(t,x,v')\{n(x) \cdot v'\}dv'.
\end{eqnarray*}
%with the Maxwellian distribution of the
%wall $M_w$:
%\begin{eqnarray*}
% M_w (v) =\frac{1}{(2\pi \theta)^{N/2}} \exp\bigg( -\frac{v^2}{2\theta}\bigg),~~~~v\in V,
%\end{eqnarray*}
%where $\theta$ being the temperature of the surface $\partial\Omega$ (which is assumed to be constant).
Here the constant $c$ is normalized as
\begin{eqnarray} c\int_{v'\in V: n(x)\cdot v'>0} \{ n(x) \cdot v'\}dv' = 1. \label{1normalized condition}\end{eqnarray}

%  The main purpose of this paper is to establish the existence and BV regularity estimate of the solution for the neutron transport problem in %non-convex domain. For the study of the neutron transport equation with diffusive boundary condition,  by setting $ u= \sqrt{M_w} \tilde{u}$ and %drop off the $\tilde{ }$ , then $u$ satisfies
%\begin{eqnarray}\label{1Main problem}
%\frac{\partial u}{\partial t} + v\cdot \nabla u +\Sigma u =  K u+q,,~~~u(0,x,v)=u_0(x,v)
%\end{eqnarray}
%with the boundary condition
%\begin{eqnarray}\label{1Boundary condition}
%u(t)|_{\gamma_-}=\mathcal{P}_{\gamma} u+ r.
%\end{eqnarray}
%Here $r,~ u_0, $ are given functions, the operator $K$ is defined with kernel $\tilde{k} = M_w^{-1/2}(v)k(x,v,v')M_w^{1/2}(v')$ with $k$ is the %kernel defined in (\ref{1operator K}). By dropping the $\tilde{ }$ again, the operator is rewritten as
%\begin{eqnarray}\label{1Operator kernel2}
%Ku = \int_{V} k(x,v,v') u(t,x,v')dv',
%\end{eqnarray}
%nd the diffusive boundary operator is rewritten as
%\begin{eqnarray*}
%\mathcal{P}_{\gamma} u = c M_w^{1/2}(v) \int_{v'\in V: n(x)\cdot v'>0} u(t,x,v')M_w^{1/2}(v')\{n(x) \cdot v'\}dv'.
%\end{eqnarray*}
%
 The operator $\mathcal{P}_{\gamma}$ could be viewed as function on $\{v\in V: v\cdot n(x)>0\}$ for any fixed $x\in \partial \Omega$, it is a $L^2_{v}-$ projection with respect to the measure $|n(x)\cdot v|dv$ for any boundary function $u$ defined on $\gamma_+$.

 Before state the main results, we give some notations. We denote $||\cdot||_{\infty}$ the norm of $L^{\infty}(\overline{\Omega}\times V)$, while $||\cdot||_p$ is the norm of the $L^p(\Omega\times V)$. In particular, $(\cdot,\cdot)$ is the inner product of the space $L^2(\Omega\times V)$. We also denote $|\cdot|_p$ the norm of $L^p(\partial \Omega \times V, dS_x dv)$ and $|\cdot|_{\gamma,p}$ the norm of $L^p(\partial \Omega \times V)= L^p(\partial \Omega \times V, d\gamma)$ with $d\gamma=|n(x)\cdot v|dS_x dv$ with the surface measure $dS_x$ on $\partial X$. We write $|\cdot|_{\gamma_{\pm},p}=|\cdot I_{\gamma_{\pm}}|_{\gamma,p}$. For a
function on $\Omega\times V$, we denote $f_{\gamma}$ to be its trace on $\gamma$ whenever it exists. $ f \lesssim g$ means $f=O(g)$.

 We show that the $L^2$ estimates of the solution for the neutron transport equation with the mixing boundary condition can be obtained by the tracing theorem. These estimates can be applied to achieve the estimates of the solution in $L^{\infty}$ norm by using the general characteristics curves of the equation with the same boundary condition. Thus, we get the existence of the solution for the neutron transport equation, which is stated as follows.

\begin{thm}\label{1Existence of the solution} Let $\Omega$ be a bounded open subset of $\mathbb{R}^3$ with $C^2$ boundary $\partial \Omega$ as in (\ref{1Non-convex Domain boundary})-(\ref{1Non-convex Domain}). Assume that $ ||u_0||_{\infty}, \sup_{0\leq t\leq T}|r(t)|_{\infty}, \sup_{0\leq t\leq T}||g(t)||_{\infty}$ are bound for any fixed $T>0$. Suppose further that there exist some constant $M_a,~M_b$, for all $(x,v) \in \Omega\times V$, it holds
\begin{eqnarray}\label{1Operator kerne1}
~0\leq \Sigma(x,v)\leq M_a , \quad 0\leq  \int_{V} k(x,v',v)d v',  \int_{V} k(x,v,v')d v' \leq M_b.
\end{eqnarray}
Then, there exists a unique solution $u\in L^{\infty}([0,T]; L^{\infty}(\Omega\times V))$ of the problem (\ref{1Main problem}) with (\ref{1Boundary condition}) such that
\begin{eqnarray}\label{1Uniform bound of solution}
||u(t)||_{\infty} \lesssim ||u_0||_{\infty} + \sup_{0\leq t\leq T}|r(t)|_{\infty}+ \sup_{0\leq t\leq T}||q(t)||_{\infty}, \quad \text{for ~all}\,\,~ 0\leq t\leq T.
\end{eqnarray}
\end{thm}
%\begin{rmk} In general, let $q=0$. If  neutron transport

%\end{rmk}
Based on the uniform estimate of solution (\ref{1Uniform bound of solution}), we study the BV-regularity of the solutions in non-convex domain. A function $f \in L^1(\Omega \times V)$ has {\bf bounded variation} in $ \Omega \times V$ if
 \begin{eqnarray}\label{1Bounded variation}
 ||f||_{\widetilde{BV}}=:\sup\bigg\{\iint_{\Omega \times V} f div \psi dxdv~:~ \psi \in C_c^1(\Omega \times V; \mathbb{R}^3),|\psi|\leq 1\bigg\}<\infty
 \end{eqnarray}
 The $BV$ space is defined as follow
\begin{eqnarray}
\bigg\{ f\in  L^1(\Omega \times V) \,~|\,~ ||f||_{BV} = ||f||_{L^1}+||f||_{\widetilde{BV}}<\infty \bigg\}.
\end{eqnarray}

For the estimate of BV norm of the solution, we should imposed some additional regularity on the data. Let $\partial =(\partial_x,\partial_v)$.
For any fixed $T$, we assume
\begin{eqnarray} \label{1Regularity on data}
||u_0||_{BV} + \sup_{0\leq t\leq T}\Big[|r|_{\infty} + |\partial_t r(t)|_1+  |\partial r(t)|_1 + ||q(t)||_{BV}\Big]< \infty,
\end{eqnarray}
 and there exist some constant $M_a',~M_b'$, it holds, for all $(x,v) \in \Omega\times V $
\begin{eqnarray}\label{1Operator kerne2}
~ \partial \Sigma(x,v)\leq M_a' , \quad   \int_{V} \partial k(x,v',v)d v',  \int_{V} \partial k(x,v,v')d v' \leq M_b'.
\end{eqnarray}
The second main result in this paper is the following.

\begin{thm} \label{1BV regularity fro neutron transport}  Let $\Omega$ is a non-convex domain with $C^2$  boundary $\partial \Omega$ as in (\ref{1Non-convex Domain boundary})-(\ref{1Non-convex Domain}). Suppose that all the conditions in Theorem \ref{1Existence of the solution} hold. Moreover, we assume that (\ref{1Regularity on data})-(\ref{1Operator kerne2}) hold.
Then, there exists a unique solution $u\in L^{\infty}([0,T]; BV(\Omega\times V))$ of the problem (\ref{1Main problem}) with (\ref{1Boundary condition}) satisfies, for all $ 0\leq t\leq T$,
\begin{eqnarray}\label{1BV bound of solution}
||u(t)||_{BV} \lesssim ||u_0||_{BV} + \sup_{0\leq t\leq T}\Big[|r|_{\infty} + |\partial_t r(t)|_1 +  |\partial r(t)|_1 + ||q(t)||_{BV}\Big],
\end{eqnarray}
and $\nabla_{x,v}u d\gamma$ is a Radon measure $\sigma_t$ on $\Omega \times V$ such that
$\int_0^T |\sigma_t(\partial \Omega \times V)|dt \lesssim ||u_0||_{BV} + \sup_{0\leq t\leq T}\Big[|r|_{\infty} + |\partial_t r(t)|_1+  |\partial r(t)|_1 + ||q(t)||_{BV}\Big].$
\end{thm}

Since the boundary operator associated with the diffusive condition is of norm exactly one, the standard theory
of transport equation in bounded domains \cite{[BP]} fails. The ideas of the proof of Theorem \ref{1Existence of the solution} is similar to that in \cite{[Guo2]} and \cite{[EGKM]}. Here, we replace the original unknown function $u$ with et $U=e^{\lambda t} u$ for any fixed time $T$. It is more convenient because there is a diffusion term $\lambda U$ in the equation. To prove the existence of the solution, we first study the equation about $U$ with a reduced diffusive reflection boundary condition, which is set up to establish a contracting map argument. Then, we take the limit and get the solution based on the uniformly estimates of the sequence. In this paper, we give a more precise estimate of the sequence $U_m$ (see Lemma \ref{2Main result for sequaences}), it can be easily applied to both the bound of the sequence and its convergence.

We now illustrate the main ideas of the proof of Theorem \ref{1BV regularity fro neutron transport}, which is similar to that in \cite{[GKTT2]}. For simplicity, we assume that $u$ satisfies the following simpler problem
\begin{eqnarray*}
\partial_t u +v\cdot \nabla_x u +\Sigma u=H,\quad u|_{t=0}=u_0, \quad u|_{\gamma_-}=\mathcal{P}_{\gamma}u +r,
\end{eqnarray*}
 and where $\Sigma\geq 0$, $H$ and $r$ are smooth enough. In general, solutions $u$ are discontinuous on
 $\mathfrak{S}_B$ and (distributional) derivatives do not exist.
In order to take derivatives, we construct some smooth cut-off function $\chi_{\varepsilon}(x, v)$ vanishing on an open neighborhood of $\mathfrak{S}_B$ and
  consider the following problem
\begin{eqnarray*}
\partial_t u^{\varepsilon} +v\cdot \nabla_x u^{\varepsilon} +\Sigma u^{\varepsilon}= \chi_{\varepsilon}H, \\
 u|_{t=0}= \chi_{\varepsilon}u_0,\quad
 u^{\varepsilon}|_{\gamma_-}= \chi_{\varepsilon} \mathcal{P}_{\gamma}u^{\varepsilon} +\chi_{\varepsilon} r.
\end{eqnarray*}
Due to the cut-off $\chi_{\varepsilon}$, the solution of $u^{\varepsilon}$ vanishes on some open subset of
$\overline{\Omega} \times V$  containing the singular set $\mathfrak{S}_B$. Therefore $u^{\varepsilon}$ is smooth and we can apply (distributional) derivatives $\partial  $ to the approximation equation.
Once we can show that $u^{\varepsilon}$ is uniformly bounded in $L^{\infty}$ and $\partial u^{\varepsilon}$ is uniformly bounded in $L^1(\Omega \times V)$, then we conclude that $u^{\varepsilon}$ converges to the solution $u$ weak-$*$ in $L^{\infty}$ and $BV$.

So, we should firstly construct the smooth cut-off function $\chi_{\varepsilon}$ such that it vanished on an open neighborhood $\mathcal{O}_{\varepsilon, \varepsilon_1}$ of $\mathfrak{S}_B$.  we can show that $\mathcal{O}_{\varepsilon,\varepsilon_1}$
 contains all points whose distance from $\mathfrak{S}_B$ is less than $\varepsilon$.
Such $\varepsilon$ thickness is important for constructing cut-off functions.
Secondly, we should control the outgoing term for the estimate of the derivatives. For this purpose, we split the outgoing boundary $\gamma_+$ into
  the (outgoing) almost grazing set
  \begin{eqnarray}\label{1Almost grazing set}
   \gamma_+^{\delta}:=\{ (x,v) \in \gamma_+,~v\cdot n(x) \leq \delta\},
  \end{eqnarray}
and the (outgoing) non-grazing set
\begin{eqnarray}\label{1Non-grazing set}
\gamma_+ \setminus \gamma_+^{\delta}:=\{ (x,v) \in \gamma_+,~v\cdot n(x) \geq \delta\}.
\end{eqnarray}
 The set $\gamma_+ \setminus \gamma^{\delta}_+$ contribution can be controlled by the bulk integration and the initial data
 by the trace theorem. While the $\gamma_+^{\delta}$ contribution cannot be bounded by the bulk integration and $\int_0^t |\partial u^{\epsilon}|_{\gamma_+,1}$ of the energy-type estimate directly. Fortunately, we extract an extra small constant in front of the term $\int_0^t |\partial u^{\epsilon}|_{\gamma_+,1}$ to bound $\partial u^{\varepsilon}$ on $\gamma_-$ by using the Duhamel formula along the trajectory (Double iteration schedule).

 The plan of this paper is the following: In section 2, we obtain the solution $U$ with a reduced diffusive reflection boundary condition. Based on the uniformly estimates of the approximation solution in $L^2$,  we take the limit and get the solution of the neutron transport equation. For the uniform bound of the approximation solution, we follows the abstract scheme appeared in \cite{[EGKM]}. Here we give a new estimate in Lemma \ref{2Main result for sequaences} which can be applied to obtain both the bound of the sequence and its convergence. In section 3, we firstly construct the desired $\varepsilon$- neighborhood of the singular set and its smooth cut-off functions $\chi_{\varepsilon}$. Then, we analyse the estimates of $\chi_{\varepsilon}$ and their derivatives in the bulk and on the boundary. Moreover, the new trace theorem is achieved by using double iteration, and we give the estimates of the approximation sequence $u^{\varepsilon,m}$ in $L^{\infty}$ and its derivatives in $L^1(\Omega\times V)$. At last, some useful geometric results will be listed in Section 4.

\section{Existence of the solution}
\setcounter{section}{2} \setcounter{equation}{0}

In this section, we consider the existence of the solution to (\ref{1Main problem}) with the boundary condition (\ref{1Boundary condition}) for all $0\leq t\leq T$. The result is as follows.
\begin{pro}\label{2Solution for the original problem}  Let $T>0$. Suppose that $||u_0||_{\infty}, \sup_{0\leq t\leq T}|r(t)|_{\infty},~ \sup_{0\leq t\leq T} ||q||_{\infty} $ are all bound. Then, there is a solution $u$ of (\ref{1Main problem})-(\ref{1Boundary condition}) such that
 for any $0\leq t\leq T$, it satisfies
\begin{eqnarray}
||u(t)||_{\infty} \lesssim ||u_0||_{\infty}+\sup_{0\leq t\leq T} \Big(|r(t)|_{\infty} + || q(t)||_{\infty}\Big). \label{2Estimate of solution}
\end{eqnarray}
\end{pro}

Set $U=e^{-\lambda t}u$ for some constant $\lambda \gg 1 $ which will be determined later, it satisfies the modified problem
\begin{eqnarray}\label{2Modified problem}
\partial_t U + v\cdot \nabla U +\lambda U = (-\Sigma + K)U+q^{\lambda},~~~
U=u_0,
\end{eqnarray}
 with $U(t)|_{\gamma_-}=\mathcal{P}_{\gamma}U + r^{\lambda}$, and
$q^{\lambda}= e^{-\lambda t}q ,~ r^{\lambda}=e^{-\lambda t}r$. The following result is obvious, we omit the proof here.

\begin{pro}\label{2equilibrium} The problem (\ref{1Main problem})-(\ref{1Boundary condition}) has a unique solution if and only if
the problem (\ref{2Modified problem}) has a unique solution. Moreover, solutions of
(\ref{1Main problem}) correspond to solutions of (\ref{2Modified problem}).
\end{pro}

The existence of the solution to (\ref{2Modified problem}) can be obtained from the following proposition.

\begin{pro}\label{2Solution for the modified problem} Let $T>0$.
$||u_0||_{\infty}, \sup_{0\leq t\leq T}|r^{\lambda}(t)|_{\infty},~ \sup_{0\leq t\leq T} ||q^{\lambda}||_{\infty} $ are all bound. Then, there is a solution $U$ of (\ref{2Modified problem}) such that
 for any $0\leq t\leq T$, it satisfies
\begin{eqnarray}
||U(t)||_{\infty} \lesssim ||u_0||_{\infty}+\sup_{0\leq t\leq T} \Big(|r^{\lambda}(t)|_{\infty} + || q^{\lambda}(t)||_{\infty}\Big). \label{2Estimate of Modified solution}
\end{eqnarray}
\end{pro}

In the following, we mainly establish the solution to (\ref{2Modified problem}).

\subsection{$L^2$ estimates of the solution}

In order to prove Proposition \ref{2Solution for the modified problem}, we need to study the $L^2$ of the solution to (\ref{2Modified problem}). The estimate are obtained from the following two lemmas. Firstly, we start with the existence of the solution to the simple transport equation. Secondly, we construct the solution to (\ref{2Modified problem}) by assuming $\lambda$ is sufficiently large.

\begin{lem} \label{2Existence with penalization inflow} Let $\lambda>0$ and $T>0$. Suppose that $Q \in L^2\cap L^{\infty}([0,T)\times \Omega\times V)$, $R\in L^2\cap L^{\infty}([0,T)\times \gamma_-), ~U_0\in L^2\cap L^{\infty}(\Omega\times V)$.
% and the compatibility condition on $\gamma_- $
%$U_0(x,v)=R(0,x,v) ~~\text{for}~~(x,v) \in \gamma_- .$
Then, for all $0\leq t\leq T$, there exists a unique solution $U(t,x,v) $ to
 \begin{eqnarray}\label{2MSimply Linear problem}
[\partial_t + v\cdot \nabla  +\lambda] U = Q,  ~~U(t)|_{\gamma_-}=R,~~U|_{t=0}=U_0.
\end{eqnarray}
such that
\begin{eqnarray}
||U(t)||_2^2 + \int_0^t \big( |U(s)|_{2}^2 + ||U(s)||_2^2 \big)ds \lesssim ||U_0||_2^2+\int_0^t\big(|R(s)|_{2,-}^2 + ||Q(s)||_2^2\big)ds  \label{2L2 Estimate of penelization0}
\end{eqnarray}
and
\begin{eqnarray}\label{2LI Estimate of penalization}
||U(t)||_{\infty} + |U(s)|_{\infty}   \lesssim ||U_0||_{\infty} + \sup_{0\leq s \leq T} |R(s)|_{\infty} + \int_0^T ||Q(s)||_{\infty}ds. \label{L2 Estimate of penelization0}
\end{eqnarray}
%Moreover if $Q$, $R$ and $U_0$ are continuous away from the grazing set $\gamma^0$,
%then $U$ is continuous away from $\mathfrak{D}$. In particular, $\mathfrak{D} =\gamma_0$ if
%$\Omega $ is convex.
\end{lem}
\prof From integration along the characteristic lines of
\begin{eqnarray*}\frac{d x}{ds} = v  \in V\quad ~~\text{and} \quad \frac{d v}{ds} = 0, \label{Characteristic lines} \end{eqnarray*}
the solution of $U(t,x,v)$ can be rewritten in the integration form as
\begin{eqnarray}\label{2Solution formula Modified linear Prob}
U(t,x,v) &=& \mathbf{1}_{\{t < t_{\mathbf{b}}\}}e^{-\lambda t} U_0(x-tv,v) \notag + \mathbf{1}_{\{t> t_{\mathbf{b}}\}}e^{-\lambda t_{\mathbf{b}}} R(t-t_{\mathbf{b}},x_{\mathbf{b}},v)\vspace{3pt}\notag\\
&&+ \int_0^{\min\{t,t_{\mathbf{b}}(x,v)\}} e^{-\lambda s}Q(t-s,x-sv,v)ds.
\end{eqnarray}
Here $t_{\mathbf{b}}(x,v)$ and $x_{\mathbf{b}}$ are defined in (\ref{1Backward exit time})-(\ref{1Backward exit P}). We then show that $U(t,x,v)$ is a weak solution of (\ref{2MSimply Linear problem}) in the sense of distributions.

 Now, we will establish the $L^2$ estimate of the solution of (\ref{2MSimply Linear problem}). Multiplying (\ref{2MSimply Linear problem}) with $U$ and integrating over $ ]0,T[\times\Omega\times V$, then Green¡¯s formula gives
\begin{eqnarray*}
 ||U(t)||_2^2 +\int_0^t\Big(|U(s)|_{2,+}^2 + 2 \lambda ||U(s)||_2^2\Big) ds = ||U_0||_2^2+\int_0^t\Big(|U(s)|_{2,-}^2 + 2 (Q,U)\Big)ds.
\end{eqnarray*}
Since $ \lambda > 0$ and $U|_{\gamma_-} =R$, by the Cauchy inequality, we get
\begin{eqnarray}\label{2L2Estimate Inequality}
 ||U(t)||_2^2 +\int_0^t\Big(|U(s)|_{2}^2 +||U(s)||_2^2\Big) ds \lesssim_{\lambda} ||U_0||_2^2+\int_0^t\Big(|R(s)|_{2,-}^2  + ||Q(s)||_2^2\Big) ds.
\end{eqnarray}
This gives the inequality (\ref{2L2 Estimate of penelization0}). The uniqueness of the solution follows from (\ref{2L2Estimate Inequality}) when $U_0=0,~R=0$ and $Q=0$. The inequality (\ref{2LI Estimate of penalization}) is easily derived from (\ref{2Solution formula Modified linear Prob}) since $\lambda >0$.$\hfill \square$\\

%In order to prove the continuity, we adapt the discussion in \cite{[EGKM]}. Let $(x, v) \in \Omega \times V \setminus \mathfrak{D}$. By the %definition of $\mathfrak{D}$, then $n(x_{\mathbf{b}}(x, v)) \cdot v < 0$, and hence $t_{\mathbf{b}}(x, v)$ is smooth from Lemma \ref{Derivatives}. %Therefore, if $Q$ $R$ and $U_0$ are continuous, $U(t,x, v)$ is continuous at $(t,x, v) \in [0,T]\times (\overline{\Omega}\times V\setminus %\mathfrak{D})$. Similar to (\cite{[EGKM]}), we can prove $\mathfrak{D}=\gamma_0$ if $\Omega$ is convex.  $\hfill \square$\\

In the next lemma, we firstly study the solution of (\ref{2Modified problem}) with a reduced diffusive reflection boundary condition,  which is necessary to establish a contracting map argument. Then, we take the limit and get the solution based on the uniformly estimates of the sequence.
\begin{lem} \label{2Existence with Collision operator} Let $T>0$ and
\begin{eqnarray}\label{2range of parameter} \lambda > \lambda_0 = 1+ M_a + M_b. \end{eqnarray}
 Suppose that $q^{\lambda} \in L^2([0,T]\times \Omega\times V)$, $r^{\lambda}\in L^2([0,T]\times \gamma_-)$ and $u_0\in L^2(\Omega\times V)$.
 %and also the compatibility condition
%$\tilde{u}_0(x,v)=\mathcal{P}_{\gamma}\tilde{u}_0(x,v)+\tilde{r}(0,x,v) ~~\text{for}~~(x,v) \in \gamma_- .$
  Then, there exists a unique solution $U$ to (\ref{2Modified problem}).
Moreover, for any $0\leq t\leq T$, the solution satisfies
\begin{eqnarray}
||U(t)||_2^2 + \int_0^t\big(|U(s)|_{2}^2 +||U(s)||_2^2\big) ds \lesssim ||u_0||_2^2+\int_0^t\Big(|r^{\lambda}(s)|_{2,-}^2+ ||q^{\lambda}(s)||_2^2\Big) ds . \label{2L2 Estimate of penelization0}
\end{eqnarray}
%and
%\begin{eqnarray}
%||U(t)||_{\infty} \lesssim ||\tilde{u}_0||_{\infty}+ \sup_{0\leq s\leq %T}(|\tilde{r}(s)|_{\infty,-}+||\tilde{q}(s)||_{\infty}).\label{LI Estimate of penalization}
%\end{eqnarray}
%Finally the limit $U $ as $j \rightarrow \infty$ of the sequence $\{ U^j\}$ exists and solves (\ref{MSimply problem %diffusive Ref Bound}) uniquely. Also, if $Q$, $R$ and $U_0$ are continuous away from the grazing set $\gamma^0$,
%then $U(t,x,v)$ is continuous away from $\mathfrak{D}$.
\end{lem}
\prof Firstly,  for any $j>0$, we consider the existence of (\ref{2Modified problem}) with the reduced diffusive reflection boundary condition
\begin{eqnarray}\label{2Reduced diffusive ref condition}
U(t,x,v)|_{\gamma_-} = (1-\frac{1}{j}) \mathcal{P}_{\gamma}U+ r^{\lambda}(t,x,v), ~~\text{for}~j>0.
\end{eqnarray}
By applying Lemma \ref{2Existence with penalization inflow} to the following iteration in both $j$ and $l$: $U^0= u_0$, and for $l\geq 0$,
\begin{eqnarray} \label{2MSimply problem jl diffusive Ref Bound}
(\partial_t + v\cdot \nabla +\lambda) U^{l+1}= (-\Sigma + K)U^l+q^{\lambda},~~~~~~U^{l+1}=u_0,
\end{eqnarray}
with $U^{l+1}|_{\gamma_-} = (1-\frac{1}{j}) \mathcal{P}_{\gamma}U^l+r^{\lambda}$.

{\it Step 1}. We fix $j>0$ and take $l \rightarrow \infty$ of the solution of (\ref{2MSimply problem jl diffusive Ref Bound}) with (\ref{2Reduced diffusive ref condition}).  Multiply $U^{l+1} $ on both sides (\ref{2MSimply problem jl diffusive Ref Bound}) and integrate over $[0,T]\times \Omega \times V$, from Green¡¯s identity, it holds that
\begin{eqnarray}\label{2Basic L2 Estimate}
&&||U^{l+1}(t)||_2^2 + \int_0^t|U^{l+1}(s)|_{2,+}^2 ds + 2\lambda \int_0^t ||U^{l+1}(s)||_2^2ds \notag\\
&&\hspace{2mm} = ||u_0||_2^2+\int_0^t|(1-\frac{1}{j}) \mathcal{P}_{\gamma}U^l+r^{\lambda}|_{2,-}^2 ds+ 2\int_0^t\Big[\big((-\Sigma+ K)U^{l}+q^{\lambda},U^{l+1}\big)\Big] ds.
\end{eqnarray}
By the choice of $\lambda$ and $\lambda_0$ in (\ref{2range of parameter}), we derive that
 \begin{eqnarray*}
&&2|\big((-\Sigma+K)U^{l},U^{l+1}\big)| \leq 2||U^{l+1}||_2 ||(-\Sigma+K)U^{l}||_2 \leq \lambda_0 ||U^{l+1}||_2^2+ \lambda_0||U^{l}||_2^2~\vspace{3pt}\\
&&2|(q^{\lambda},U^{l+1})| \leq (\lambda-\lambda_0)||U^{l+1}||_2^2 +\frac{4}{\lambda-\lambda_0}||q^{\lambda}||_2^2.
\end{eqnarray*}
 Moreover, there is $C_j > 0$ such that
\begin{eqnarray}\label{2Estimate of Boundary term}
|(1-\frac{1}{j})\mathcal{P}_{\gamma}U^l+r^{\lambda}|_{2,-}^2 \leq |(1-\frac{1}{j})\mathcal{P}_{\gamma}U^l|_{2,-}^2+
\frac{1}{2j^2} |\mathcal{P}_{\gamma}U^l|^2_{2,-}+C_j|r^{\lambda}|^2_{2,-}.
\end{eqnarray}
For simplicity, we denote that
\begin{eqnarray*}
\mathcal{E}:=  ||u_0||_2^2+\int_0^t\Big(|r^{\lambda}(s)|_{2,-}^2 +||q^{\lambda}(s)||_2^2\Big)ds.
\end{eqnarray*}
It is easy to know that $|\mathcal{P}_{\gamma} U^l |_{2,-}^2 \leq | U^l|^2_{2,+}$. (\ref{2Basic L2 Estimate}) derives to
\begin{eqnarray*}
&&||U^{l+1}(t)||_2^2 +\int_0^t\Big[|U^{l+1}(s)|_{2,+}^2 + \lambda ||U^{l+1}(s)||_2^2\Big] ds \notag\\
&& \hspace{5mm}\leq \int_0^t \Big[(1-\frac{2}{j}+\frac{3}{2j^2})|U^l(s)|_{2,+}^2 +\lambda_0 ||U^{l}(s)||^2_2\Big] ds + C_{\lambda,j}\mathcal{E}.
\end{eqnarray*}

Since  $\lambda >\lambda_0$ and $1-\frac{2}{j}+\frac{3}{2j^2}<1$, there is some $\eta_{\lambda, j} < 1$ such that
\begin{eqnarray*}\label{2Main estimate of the approximate sequence}
 &&||U^{l+1}(t)||_2^2 +\Big(\int_0^t|U^{l+1}(s)|_{2,+}^2 ds + \lambda \int_0^t||U^{l+1}||_2^2 ds\Big) \notag\\
&& \hspace{5mm}\leq\eta_{\lambda,j} \Big(\int_0^t( |U^l|_{2,+}^2 +\lambda||U^{l}||_2^2)ds\Big) + C_{\lambda,j} \mathcal{E}.
\end{eqnarray*}
Since $U^0=u_0$, we iterate again to obtain
\begin{eqnarray*}\label{Basic L2 Estimate}
 &&||U^{l+1}(t)||_2^2 +\Big(\int_0^t|U^{l+1}(s)|_{2,+}^2 + \lambda ||U^{l+1}||_2^2 ds\Big) \notag\\
&& \leq \eta_{\lambda,j}^2  \Big(\int_0^t( |U^{l-1}|_{2,+}^2 +\lambda||U^{l-1}||_2^2)ds\Big) + (1+\eta_{\lambda,j})C_{\lambda,j}\mathcal{E}\notag\\
& & \cdots \notag\\
& & \leq \eta_{\lambda,j}^{l+1}  \Big(\int_0^t( |U^0|_{2,+}^2 + \lambda ||U^0||_2^2)ds\Big)+ \frac{1+\eta_{\lambda,j}^{l+1}}{1-\eta_{\lambda,j}}C_{\lambda,j}\mathcal{E} \notag\\
&& \leq \eta_{\lambda,j}^{l+1} t( |U^0|_{2,+}^2 + \lambda||U^0||_2^2)ds + \frac{1+\eta_{\lambda,j}^{l+1}}{1-\eta_{\lambda,j}}C_{\lambda,j}\mathcal{E}.
\end{eqnarray*}
So, we get the following  uniform estimates of $U^{l}$ with respect to $l$
\begin{eqnarray}\label{2Basic Uniform L2 Estimate}
 ||U^{l+1}(t)||_2^2 +\int_0^t|U^{l+1}(s)|_{2,+}^2 ds + \lambda \int_0^t||U^{l+1}||_2^2 ds ~\lesssim_{\lambda,j,T}~ \mathcal{E}.
\end{eqnarray}
Now, taking the difference of $U^{l+1} - U^l$, it satisfies
\begin{eqnarray}\label{2Difference of sequaence}
 \Big[\partial_t + v\cdot \nabla  +\lambda\Big] (U^{l+1}-U^l)= (-\Sigma + K)(U^l-U^{l-1}),
\end{eqnarray}
with $(U^{l+1}-U^l)|_{\gamma_-}=(1-\frac{1}{j})\mathcal{P}_{\gamma}(U^l-U^{l-1}), (U^{l+1}-U^l)(0)=0.$ Similar to the estimate of (\ref{2Basic Uniform L2 Estimate}), we yield
\begin{eqnarray}\label{2Basic L2 Estimate of Difference}
&&||U^{l+1}(t)-U^l(t)||_2^2 +\int_0^t\Big(|(U^{l+1}(s)-U^l(s)|_{2,+}^2 + \lambda ||U^{l+1}(s)-U^l(s)||_2^2 \Big)ds \notag\\
&& \leq \eta_{\lambda,j}^{l}\int_0^t \Big(|(U^1-U^0)(s)|_{2,+}^2 + \lambda ||(U^1-U^0)(s)||_2^2\Big)ds
\end{eqnarray}
From (\ref{2Basic Uniform L2 Estimate}), we know that $\int_0^t( |(U^1-U^0)(s)|_{2,+}^2 +||(U^1-U^0)(s)||_2^2)ds <\infty$ for fixed $t$. It concludes that $\{U^{l}\}_{l \geq 0}$ is a Cauchy sequence with respect to $l$ for $\eta_{\lambda,j}<1$.

Let $l \rightarrow \infty $ to obtain $U_j$ as a solution of
 \begin{eqnarray} \label{2Reduced boundary problem}
 (\partial_t +v\cdot \nabla +\lambda)U_j=(\Sigma-K)U_j+q^{\lambda},~~U_j(0)=u_0
 \end{eqnarray}
 with $U_j|_{\gamma_-} = (1-\frac{1}{j})\mathcal{P}_{\gamma}U_j+r^{\lambda}$ for all $0\leq t\leq T$. Moreover, the estimate of (\ref{2Main estimate of the approximate sequence}) derives to
\begin{eqnarray}\label{2Estimate of the approximateS1}
&&||U_j(t)||_2^2 +\int_0^t\Big[\frac{1}{j}|\mathcal{P}_{\gamma} U_j(s)|_{2,+}^2 + |(I-\mathcal{P}_{\gamma})U_j(s)|_{2,+}^2 \Big]ds + \lambda \int_0^t||U_j(s)||_2^2 ds  \lesssim_{\lambda,j,T} \mathcal{E}.
\end{eqnarray}

 {\it Step 2}. Let $j\rightarrow \infty$ for $U_j$. It needs a uniform estimate of $U_j$ w.r.t. $j$. Multiply (\ref{2Reduced boundary problem}) with $U_j$ and integrate over $[0,T]\times \Omega \times V$, then Green's identity derives to
\begin{eqnarray*}
&&||U_j(t)||_2^2 +\int_0^t|U_j(s)|_{2,+}^2 ds + 2\int_0^t\lambda||U_j(s)||_2^2ds \notag\\
&& = ||u_0||_2^2+\int_0^t|\big((1-\frac{1}{j}) \mathcal{P}_{\gamma}U_j+r^{\lambda}\big)(s)|_{2,-}^2 ds + 2\int_0^t\Big[((\Sigma-K)U_j,U_j)+(q^{\lambda},U_j)\Big]ds.
\end{eqnarray*}
 For any $\eta>0$ and $j> 0$, $|(1-\frac{1}{j})\mathcal{P}_{\gamma}u_j+r^{\lambda}|_{2,-}^2 $ can be rewritten as
 \begin{eqnarray}
&&  |(1-\frac{1}{j}) \mathcal{P}_{\gamma}U_j+r^{\lambda}|_{2,-}^2 \notag\\ \
&& = (1-\frac{1}{j})^2 |\mathcal{P}_{\gamma} U_j|_{2,-}^2 +2(1-\frac{1}{j}) \int_{\gamma-}\mathcal{P}_{\gamma}U_j r^{\lambda} d\gamma +|r^{\lambda}|_{2,-}^2 \\
 && \leq  (1+\eta)|\mathcal{P}_{\gamma} U_j|_{2,-}^2  +C_{\eta} |\tilde{r}|_{2,-}^2. \notag \label{2Estimate of Diffusion boundary}
 %&\leq & |\mathcal{P}_{\gamma} u_j|_{2,-}^2 +C_{\eta} |r|_{2,-}^2+\eta |u_j|_{2,+}^2
\end{eqnarray}
Because $|\mathcal{P}_{\gamma} U_j|_{2,-}^2=|\mathcal{P}_{\gamma} U_j|_{2,+}^2 $ and $| U_j|_{2,+}^2 = |\mathcal{P}_{\gamma} U_j|_{2,+}^2+ |(I-\mathcal{P}_{\gamma}) U_j|_{2,+}^2 $, together with the fact that $2((\Sigma-K)U_j,U_j) \leq 2\lambda_0 ||U_j(s)||_2^2$ and $2\int U_j q^{\lambda}\leq (\lambda-\lambda_0)||U_j||_2^2+\frac{1}{\lambda-\lambda_0}||q^{\lambda}||_2^2$,  we derive that
\begin{eqnarray}
&&||U_j(t)||_2^2 +\int_0^t |(\mathbf{1}-\mathcal{P}_{\gamma})U_j(s)|_{2,+}^2 ds +  \int_0^t ||U_j(s)||_2^2 ds \notag\\
&& \leq  \eta \int_0^t |\mathcal{P}_{\gamma} U_j(s)|_{2,+}^2ds + C_{\eta,\lambda,T} \mathcal{E}. \label{2Basic L2 Estimate3}
\end{eqnarray}
%Since $||P u_j||_2^2 \leq ||u_j||_2^2$, integrating (\ref{Basic L2 Estimate3}) from $0$ to $t$, we have
%\begin{eqnarray}
%\int_0^t||Pu_j(s)||_2^2ds \lesssim_t ||u_I||_2^2+ C_{\varepsilon,\eta}\int_0^t\Big(|r(s)|_{2,-}^2+||h(s)||_2^2\Big) %ds+\eta \int_0^t |\mathcal{P}_{\gamma} u_j(s)|_{2,-}^2ds
%\end{eqnarray}
%From (\ref{Diffusive Boundary Condition}) we can write $\mathcal{P}_{\gamma} u_j = z_{\gamma}(x) $ for a suitable %function $z_{\gamma}(x)$.
Now, we study the estimate of $\mathcal{P}_{\gamma}U_j$ by the trace theorem which has appeared in \cite{[EGKM]} and \cite{[GKTT1]}. For the purpose of it, we consider the boundary contribution
\begin{eqnarray*}
\int_0^t |\mathcal{P}_{\gamma}U_j(s)|_{2,\pm}^2ds = c^2 \int_0^t\int_{\gamma_{\pm}}\Big[ \int_{\{v': n(x)\cdot v'>0 \}}U_j(s,x,v') \{n(x)\cdot v'\}dv'\Big]^2 d\gamma ds.
\end{eqnarray*}
We split the domain of inner integration as
\begin{eqnarray*}
& &\{v' \in V : n(x) \cdot v' > 0\} = \{v' \in V : 0 < n(x) \cdot v' < \epsilon ~\text{or}~|v'|\leq \epsilon\}\vspace{2mm}\\
&& \hspace{4cm}\cup \, \{v' \in V : n(x) \cdot v' \geq \epsilon ~\text{and}~|v'|\geq \epsilon\}.
\end{eqnarray*}
The first set¡¯s contribution (grazing part) of $\int_0^t |\mathcal{P}_{\gamma}U_j(s)|_{2,\pm}^2ds $ is bounded by the H\"{o}lder inequality,
\begin{eqnarray}\label{2Estimate of grazing}
 && c \int_0^t\int_{\partial \Omega} c \int_{ V } \{n\cdot v\} dv\int_{\{ 0 < n \cdot v' < \epsilon
 ~\text{or}~|v'|\leq \epsilon \}} \{n\cdot v'\} dv' \notag\\
 && \hspace{3cm} \times  \int_{\{v':n\cdot v'>0\}} |U_j(s,x,v')|^2\{n\cdot v'\}dv' dS_x ds \notag \\
 &&\lesssim_{\Omega} \epsilon \times \int_0^t\int_{\gamma_+} |U_j(s)|^2 d\gamma ds.
\end{eqnarray}
Here we have used the fact that
\begin{eqnarray*}
\int_{|n\cdot v'|\leq \epsilon}  |\{ n\cdot v'\}|dv' \lesssim \epsilon, \quad c\int_{V} |\{ n\cdot v\}|dv =1.
\end{eqnarray*}
For the bound of the second set¡¯s contribution
(non-grazing part) of $\int_0^t |\mathcal{P}_{\gamma}U_j(s)|_{2,\pm}^2ds $, we use Lemma \ref{Trace theorem} and (\ref{2Basic L2 Estimate3}). From the equation, $(\partial_t+ v \cdot \nabla_x)( U_j)^2 = -2 \lambda[U_j]^2 -2 U_j (\Sigma-K)U_j + 2 U_j q^{\lambda}$. Taking the absolute value and integrating on $\Omega\times V$, for $\lambda \gg 1$,  we have
\begin{eqnarray*}
  \int_0^t ||(\partial_t+ v \cdot \nabla_x)( U_j)^2 (s)||_1 ds  \leq    4 \int_0^t \Big[\lambda ||U_j(s)||_2^2+ ||q^{\lambda}(s)||_2^2 \Big]ds.
\end{eqnarray*}
The trace theorem \ref{Trace theorem} gives
\begin{eqnarray}\label{2Estimate of Non-grazing}
 \int_0^t |U_j\mathbf{1}_{\gamma_+ \setminus \gamma_+^{\varepsilon} }(s)|_2^2ds  &\lesssim _{\varepsilon,\Omega,\lambda,T}& ||u_0||_2^2 + \int_0^t \Big[||U_j(s)||_2^2+||(\partial_t+ v\cdot \nabla_x ) (U_j)^2(s)||_1\Big] ds \notag\\
& \lesssim _{\varepsilon,\Omega,\lambda,T}&  ||u_0||_2^2+ \int_0^t \Big[ ||q^{\lambda}(s)||_2^2+ ||U_j(s)||_2^2 \Big]ds.
\end{eqnarray}
From (\ref{2Estimate of grazing}) and (\ref{2Estimate of Non-grazing}), we have
\begin{eqnarray}\label{2Estimate of Projection on boundary}
 \int_0^t |\mathcal{P}_{\gamma}U_j(s)|_{2,\pm}^2ds \lesssim_{\varepsilon,\Omega,\lambda,T}   \int_0^t \Big(\varepsilon|U_j(s)|^2_{2,+} + ||U_j(s)||_2^2\Big) ds + \mathcal{E}.
\end{eqnarray}
Combining (\ref{2Basic L2 Estimate3}) and (\ref{2Estimate of Projection on boundary}) with small $\eta > 0$ and $\varepsilon > 0$, we have the following uniform estimate
\begin{eqnarray} \label{2Uniform Estimate ACT j}
||U_j(t)||_2^2 +\int_0^t |U_j(s)|_2^2 ds + \lambda \int_0^t ||U_j(s)||_2^2 ds &\lesssim _{\varepsilon,\Omega,\lambda,T} & \mathcal{E}.
\end{eqnarray}

By taking a weak limit, we obtain a weak solution $U$ to (\ref{2Modified problem}) with the same bound
(\ref{2Uniform Estimate ACT j}). Taking the difference, we have
\begin{eqnarray}
\Big(\partial_t +v\cdot \nabla_x +\lambda\Big)[U_j-U]=(-\Sigma-K)(U_j-U), ~~~[U_j-U](0)=0
\end{eqnarray}
with $[U_j-U]|_{\gamma_-} = \mathcal{P}_{\gamma}[U_j-U]+ \frac{1}{j}\mathcal{P}_{\gamma}U_j $. Applying (\ref{2Uniform Estimate ACT j}) with $r^{\lambda} = \frac{1}{j} \mathcal{P}_{\gamma}U_j$, we obtain, as $j\rightarrow \infty$
\begin{eqnarray}
||U_j(t)||_2^2 +\int_0^t \big( |U_j(s)|_{2,+}^2 + ||U_j(s)||_2^2\big) ds \lesssim \frac{1}{j^2} \int_0^t |\mathcal{P}_{\gamma}U_j(s)|^2 ds \rightarrow 0
\end{eqnarray}
 because of (\ref{2Estimate of the approximateS1}). We yield that the limit $U$ is $L^2$ solution to (\ref{2Modified problem}) for all $0 \leq t\leq T $. Moreover, the estimate of (\ref{2L2 Estimate of penelization0}) is easily obtained from (\ref{2Uniform Estimate ACT j}). The proof of Lemma \ref{2Existence with Collision operator} is completed. $\hfill\square$

\subsection{$L^{\infty}$  estimate of the solution}
We would study the uniform $L^{\infty}$ estimates of the solution for the problem (\ref{2Modified problem}). To bootstrap $L^2$ estimate into $L^{\infty}$ estimate, we need to define the stochastic cycles for the generalized characteristic lines interacted with the boundary. This method was firstly introduced by Guo in \cite{[Guo2]}, which is a canonical way to treat $L^{\infty}$  estimate of the solution to Boltzmann equation with diffusive boundary condition. In the following, we construct the stochastic cycles for neutron transport equation with diffusive boundary condition, which is similar to that for Boltzmann equation in \cite{[Guo2]} and \cite{[EGKM]}. Then, we show $L^{\infty}$ estimate of the solution to the neutron transport equation (\ref{2Modified problem}) by this stochastic cycles.

 Let $\mathcal{V}(x)=\{v'\in D: v'\cdot n(x) >0\}$, the probability measure $d\sigma=d\sigma(x)$ is given by
\begin{eqnarray}\label{Probability measure on boundary}
d\sigma(x)=c \{n(x)\cdot v'\} dv', ~~~\text{with}\quad c\int_{\mathcal{V}(x)} d\sigma(x)=1. \label{2Measure and Weighted Function}
\end{eqnarray}

 \begin{defi}\label{2Stochastic Cycles} (Stochastic Cycles).
 Fix any point $t>0$ and $ (x,v) \notin \gamma_0 \cap \gamma_-$. Let $(t_0,x_0,v_0)= (t,x,v)$. For $v_{k+1} \in \mathcal{V}_{k+1} =\{v_{k+1} \cdot n(x_{k+1}) >0 \}$, define the $(k+1)-$ component of the back-time cycle as
\begin{eqnarray}
(t_{k+1},x_{k+1},v_{k+1}) =(t_k-t_{\mathbf{b}}(x_k,v_k),x_{\mathbf{b}}(x_k,v_k),v_{k+1}). \label{2Subsequence of back time cycle}
\end{eqnarray}
And the stochastic cycle is defined as
\begin{eqnarray*}
X_{\mathbf{cl}}(s;t,x,v)=\sum_{l}\mathbf{1}_{[t_{l+1},t_l)}(s)\{x_l+(s-t_l)v_l \},\quad V_{\mathbf{cl}}(s;t,x,v)=\sum_l\mathbf{1}_{[t_{l+1},t_l)}(s)v_l.
\end{eqnarray*}
We define the iterated integral for $k\geq 2$,
\begin{eqnarray}
\int_{\Pi_{l=1}^{k-1}\mathcal{V}_l } \Pi_{l=1}^{k-1} d\sigma_l \equiv \int_{\mathcal{V}_1} \cdots \Big\{\int_{\mathcal{V}_{k-1}}d\sigma_{k-1}\Big\} \cdots d\sigma_1.
\end{eqnarray}
\end{defi}
Here $v_l~(l=1,2,\cdots, k)$ are all independent variables and $t_k,~x_k$ depend on $t_l,~x_l,~v_l$ for
$l\leq k-1$. The phase space $\mathcal{V}_{l}$ implicitly depends on $(t,x,v,v_1,v_2,\cdots,v_{l-1})$. We show that the set in the phase space $\Pi_{l=1}^{k-1}\mathcal{V}_l$ not reaching $t=0$ after $k$ bounces is small when $k$ is large.

 \begin{lem}\label{2Small for large-time bounce}
Fixed $T>0$. For any $\varepsilon >0$, there exists $k_0(\varepsilon, T)$ such that for  $k\geq k_0$, for all  $(t,x,v) \in [0,T] \times \overline{\Omega}\times V$,
\begin{eqnarray}
\int_{\Pi_{l=m}^{k-1}\mathcal{V}_l} \mathbf{1}_{\{t_k(t,x,v,v_1,v_2\cdots,v_{k-1})>0\}}\Pi_{m=1}^{k-1}d\sigma_m \leq \varepsilon. \label{2Measure Bounce For Many times}
\end{eqnarray}
%Moreover, we also have,
%\begin{eqnarray}\label{Measure on boundary}
%\int_{\Pi_{m=1}^{k-1}\mathcal{V}_j}\sum_{m=1}^{k-1} \mathbf{1}_{\{t_{m+1}\leq 0<t_m\}}\Pi_m^{k-1} d\sigma_m \leq
%\end{eqnarray}
%Furthermore, for $T_0$ sufficiently large, there exist constant $C_1, ~C_2>0$, independent of $T_0$, such that for $k=C_1 T_0^{5/4}$,
%\begin{eqnarray}
%\int_{\Pi_{l=1}^{k-1}\mathcal{V}_l} \mathbf{1}_{\{t_k(t,x,v,v_1,v_2\cdots,v_{k-1})>0\}}\Pi_{l=1}^{k-1}d\sigma_l \leq \{\frac{1}{2}\}^{C_2T_0^{5/4}}. %\label{Measure Bounce For large time}
%\end{eqnarray}
\end{lem}
\prof Choosing $0<\delta <1$ sufficiently small, the non-grazing sets for $1\leq m\leq k-1$ is defined as
\begin{eqnarray*}
\mathcal{V}_m^{\delta}=\{v_m\in \mathcal{V}_m: v_m\cdot n(x_m)\geq \delta\} .%\cap \{v_l\in \mathcal{V}_l:|v_l| \leq \frac{1}{\delta}\}.
\end{eqnarray*}
For any $m$, by a change of variable $v_{\parallel}= \{n(x_m)\cdot v_m\} n(x) $ and $v_{\perp}=v_m-v_{\parallel}$ for $|v_{\parallel}|\leq |n(x_m)\cdot v_m| \leq \delta$, the measure of the grazing set is estimated as
\begin{eqnarray*}
\int_{\mathcal{V}_m \setminus \mathcal{V}_m^{\delta}}d\sigma_m \leq \int_{v_m\cdot n(x_m) \leq \delta} M_w(v_m)dv_m \leq C \int_{-\delta}^{\delta}dv_{\parallel} \int_{\mathbb{R}^2} e^{-\frac{|v_{\perp}|^2}{2\theta}} dv_{\perp}  \leq C\delta,
\end{eqnarray*}
where $C$ is independent of $m$. On the other hand, if $v_m\in \mathcal{V}_m^{\delta}$, then from diffusive back-time cycle, we have $x_m-x_{m+1}=(t_m-t_{m+1})v_m$. From Lemma \ref{Derivatives}, since  $v_m\cdot n(x_m) \geq \delta$ and $v_m$ is bounded, then
$(t_m-t_{m+1}) \geq \frac{\delta}{C_{\xi}}.$
Therefore, if $0< t_k(t,x,v,v_1,v_2,\cdots,v_{k-1})  \leq T$, then there can be at most $[\frac{C_{\xi}T}{\delta}]+1$ number of $v_m\in \mathcal{V}_m^{\delta}$ for $1\leq m\leq k-1$. We therefore have
\begin{eqnarray*}
&&\int_{\mathcal{V}_l}\cdots\Big\{\int_{\mathcal{V}_{k-1}}\mathbf{1}_{t_k>0}d\sigma_{k-1} \Big\}d\sigma_{k-2}\cdots d\sigma_1 \\
&&\leq \sum_{m=1}^{[\frac{C_{\xi}T}{\delta}]+1}\int_{\{\text{There\,are \,exactly\,m \,of}~ v_{l_i} \in \mathcal{V}_{l_i}^{\delta} \,\text{and}\,k-1-m \,\text{of}\,v_{l_i} \notin \mathcal{V}_{l_i}^{\delta}\}} \Pi_{m=1}^{k-1}d\sigma_m \\
&&\leq \sum_{m=1}^{[\frac{C_{\xi}T}{\delta}]+1}\Big(\begin{array}{ccc}k-1\\m \end{array}\Big)\Big| \sup_j\int_{\mathcal{V}_j^{\delta}} d\sigma_j\Big|^m \Big\{\sup_{j}\int_{\mathcal{V}_j \setminus \mathcal{V}_j^{\delta}}d\sigma_m\Big\}^{k-m-1}.
\end{eqnarray*}
Since $d\sigma_m$ is a probability measure $\int_{\mathcal{V}_m^{\delta}}d\sigma_m \leq 1$ and
\begin{eqnarray*}
\Big\{\int_{\mathcal{V}_j \setminus \mathcal{V}_j^{\delta}}d\sigma_j \Big\}^{k-m-1} \leq
\Big\{\int_{\mathcal{V}_j \setminus \mathcal{V}_j^{\delta}}d\sigma_j \Big\}^{k-2-[\frac{C_{\xi}T}{\delta}]}
\leq (C \delta)^{k-2-[\frac{C_{\xi}T}{\delta}]}.
\end{eqnarray*}
But \begin{eqnarray*}
\sum_{m=1}^{[\frac{C_{\xi}T}{\delta}]+1}\Big(\begin{array}{ccc}k-1\\m \end{array}\Big) \leq \sum_{m=1}^{[\frac{C_{\xi}T}{\delta}]+1}\frac{(k-1)^m}{m!} \leq (k-1)^{[\frac{C_{\xi}T_0}{\delta}]+1} \sum_{m=1}^{[\frac{C_{\xi}T}{\delta}]+1}\frac{1}{m!}\leq (k-1)^{[\frac{C_{\xi}T}{\delta}]+1} . \end{eqnarray*} it deduces that
\begin{eqnarray}
\int_{\Pi_{l=m}^{k-1}\mathcal{V}_m} \mathbf{1}_{\{t_k >0\}} \Pi_{l=1}^{k-1}d\sigma_m \leq (k-1)^{[\frac{C_{\xi}T}{\delta}]+1} (C\delta)^{k-2-[\frac{C_{\xi}T}{\delta}]}. \label{Measure for K times}
\end{eqnarray}
For $\varepsilon >0$, (\ref{2Measure Bounce For Many times}) follows for $C\delta <1$ and $k \gg [\frac{C_{\xi}T}{\delta}]+1 $. For example, we can choose
$k= 15\{[\frac{C_{\xi}T}{\delta}]+1 \}+2$ for small $\delta>0$.
$\hfill\square$\\
%To prove (\ref{Measure Bounce For large time}), we let $k-2=N\{ [\frac{C_{\xi}T_0}{\delta^3}]+1\}$, so that if
%$\frac{C_{\xi}T_0}{\delta^3} \geq 1$, (\ref{Measure for K times}) can be further estimated by
%\begin{eqnarray*}
%\left\{N(\frac{C_{\xi}T_0}{\delta^3}+1)(C \delta)^N \right\}^{[\frac{C_{\xi}T_0}{\delta^3}]+1} \leq
%\left\{\frac{2NC_{\xi}T_0}{\delta^3} (C \delta)^N \right\}^{[\frac{C_{\xi}T_0}{\delta^3}]+1} \leq
%\left\{C_{N,\xi}T_0 \delta^{N-3} \right\}^{[\frac{C_{\xi}T_0}{\delta^3}]+1}.
%\end{eqnarray*}
%We choose $C_{N,\xi}T_0 \delta^{N-3} =\frac{1}{2}$, so that $\delta =(\frac{1}{2C_{N,\xi}T_0})^{\frac{1}{N-3}}$ is small for $T_0$ large and for %$N\geq 3$. Moreover
%\begin{eqnarray*}
% [\frac{C_{\xi}T_0}{\delta^3}]+1 \backsim C_{N,\xi}T_0^{1+\frac{3}{N-3}}
%\end{eqnarray*}
%and $\frac{C_{N,\xi}T_0}{\delta^3}\geq 2 $ if $T_0$ is large so that we can close estimate.

%Finally we choose $N=15$. For $T_0$ sufficiently large, $[\frac{C_{\xi}T_0}{\delta^3}]+1 \backsim CT_0^{5/4}$ and
%$k=16\{[\frac{C_{\xi}T_0}{\delta^3}]+1\}+2 \backsim CT_0^{5/4} $ and (\ref{Measure Bounce For large time}) follows.

Let $h$ satisfies the following neutron transport equation with in-flow boundary condition,
\begin{eqnarray*}\{\partial_t +v\cdot \nabla_x +\lambda +\Sigma \}h =q, \quad h|_{t=0}=h_0,\quad h|_{\gamma_-}=g.\end{eqnarray*}
For $(x,v) \notin \gamma_0$, we denote its backward exit point as $[t-t_{\mathbf{b}}(x,v),x_{\mathbf{b}}(x,v),v]$. Since  $\frac{d}{d\tau} \{e^{\int_0^{\tau}(\lambda+\Sigma)} h \}=q$ along the characteristic $\frac{d x}{d\tau}=v,~\frac{d v}{d\tau}=0$. Thus, if $t-t_{\mathbf{b}}<0$,
\begin{eqnarray*}
h(t,x,v) = e^{-\int_0^t (\lambda+\Sigma)}h_0(x-vt,v)+\int_0^t e^{\int_0^{t-s}(\lambda+\Sigma)} q(s,x-v(t-s),v)ds.
\end{eqnarray*}
If  $t-t_{\mathbf{b}}>0$, we have
\begin{eqnarray*}
h(t,x,v) = e^{-\int_0^{t_{\mathbf{b}}} (\lambda+\Sigma)}g(t-t_{\mathbf{b}},x_{\mathbf{b}},v)+\int_{t-t_{\mathbf{b}}}^t e^{\int_0^{t-s}(\lambda+\Sigma)} q(s,x-v(t-s),v)ds.
\end{eqnarray*}
Recall the diffusive cycles Definition \ref{2Stochastic Cycles}, we have the following iteration scheme for the neutron transport equation with the mixing boundary condition. The proof is similar to that
of Boltzmann equation in \cite{[EGKM]}, we omit it here.
\begin{lem}\label{Stochastic cycles for solution} \cite{[EGKM]} Assume that $h$ satisfies
\begin{eqnarray*}\{\partial_t +v\cdot \nabla_x +\lambda +\Sigma \}h =q, ~~h|_{t=0}=h_0,~~h|_{\gamma_-}= \mathcal{P}_{\gamma}h +r.\end{eqnarray*}
 Then, for almost every $(x,v)\notin \gamma_0$, if $t_1(t,x,v) \leq 0$,
\begin{eqnarray*}
h(t,x,v) =e^{-\int_0^t (\lambda+\Sigma)}h_0(x-vt,v)+\int_0^t e^{-\int_0^{t-\tau} (\lambda+\Sigma)} q(\tau,x-v(t-\tau),v)d\tau.
\end{eqnarray*}
If $t_1(t,x,v) >0$, then for $ k\geq 2 $,
\begin{eqnarray*}
h(t,x,v) &= &\int_{t_1}^te^{-\int_0^{t-\tau} (\lambda+\Sigma)}q(\tau,x-v(t-\tau),v)d\tau + e^{-\int_0^{t-t_1} (\lambda+\Sigma)} |r(t_1,x_1,v)| \\
         &  & + e^{-\int_0^{t-t_1} (\lambda+\Sigma)} \int_{\Pi_{m=1}^{k-1}\mathcal{V}_m}H.
\end{eqnarray*}
where $H$ is given by
\begin{eqnarray}
&&\sum_{m=1}^{k-1} \mathbf{1}_{\{t_{m+1} \leq 0 <t_m\}} h_0(x_m- v_m t_m,v_m)d\Sigma_m(0) \notag\\
 &&+ \sum_{m=1}^{k-1} \int_0^{t_m} \mathbf{1}_{\{t_{m+1}\leq 0<t_m \}}q(\tau,x_m+(\tau-t_m)v_m,v_m)d\Sigma_m(\tau)d\tau\notag \\
&& + \sum_{m=1}^{k-1} \int_{t_{m+1}}^{t_m}\mathbf{1}_{t_{m+1}>0}q(\tau,x_m+(\tau-t_m)v_m,v_m)d\Sigma_m(\tau)d\tau
\notag\\ &&
+\mathbf{1}_{\{t_k>0\}} h(t_k,x_k,v_{k-1})d\Sigma_{k-1}(t_k) + \sum_{m=1}^{k-1}\mathbf{1}_{\{t_m>0\}}d\Sigma_m^r,
\end{eqnarray}
with $d\Sigma_m(0),~d\Sigma_{k-1}(t_k)$ are evaluated at $s=0$ and $s=t_k$ of
\begin{eqnarray} \label{2Bounce of source}
d\Sigma_m(s)=\{\Pi_{j=l+1}^{k-1}d\sigma_j \}\{e^{\nu(v_l)(s-t_m)}d\sigma_m\}\Pi_{j=1}^{m-1}\{e^{-\int_0^{t_j-t_{j+1}} (\lambda+\Sigma)}d\sigma_j\}.
\end{eqnarray}
and
\begin{eqnarray}\label{2Bounce of boundary condition}
d\Sigma_m^r = \Pi_{j=l+1}^{k-1}d\sigma_j  \{e^{-\int_0^{t_m-t_{m+1}} (\lambda+\Sigma)} r(t_{m+1},x_{m+1},v_m)d\sigma_m\} \Pi_{j=1}^{m-1} \{e^{-\int_0^{t_j-t_{j+1}} (\lambda+\Sigma)}d\sigma_j\}.
\end{eqnarray}
\end{lem}

Now, we consider the solution of (\ref{2MSimply problem jl diffusive Ref Bound}). By Lemma \ref{Stochastic cycles for solution}, we have
\begin{eqnarray}\label{2Estimate stachastic 0}
 U^{l+1}(t,x,v) & \leq & \mathbf{1}_{\{t_1 \leq 0\}} e^{-\int_0^{t} (\lambda+\Sigma)} |u_0(x-tv,v)| \notag\\
 &&+\mathbf{1}_{\{t_1 \leq 0\}}\int_0^t e^{-\int_0^{t-\tau} (\lambda+\Sigma)} |[K U^l+q^{\lambda}](\tau, x-(t-\tau)v,v)|d\tau \notag \\
 && + \mathbf{1}_{\{t_1 > 0\}} \int_{t_1}^t e^{-\int_0^{t-\tau} (\lambda+\Sigma)} |[K U^l+q^{\lambda}](\tau, x-(t-\tau)v,v)|d\tau \notag \\
&&+ \mathbf{1}_{\{t_1>0\}} e^{-\int_0^{t-t_1} (\lambda+\Sigma)} |r^{\lambda}(t_1,x_1,v)|  + e^{-\int_0^{t-t_1} (\lambda+\Sigma)} \int_{\Pi_{m=1}^{k-1}\mathcal{V}_m}|H|,
\end{eqnarray}
where $H$ is bounded by
\begin{eqnarray}
&&\sum_{m=1}^{k-1} \mathbf{1}_{\{t_{m+1} \leq 0 <t_m\}} |u_0(x_m- v_m t_m,v_m)|d\Sigma_m(0) \label{2Estimate stachastic a}\\
 &&+ \sum_{m=1}^{k-1} \int_0^{t_m} \mathbf{1}_{\{t_{m+1}\leq 0<t_m \}} |[KU^{l-m}+q^{\lambda}](\tau,x_m+(\tau-t_m)v_m,v_m)|d\Sigma_m(\tau)d\tau \label{2Estimate stachastic b}\\
&& + \sum_{m=1}^{k-1} \int_{t_{m+1}}^{t_m}\mathbf{1}_{\{t_{m+1}>0\}}|[KU^{l-m}+q^{\lambda}](\tau,x_m+(\tau-t_m)v_m,v_m)|d\Sigma_m(\tau)d\tau \label{2Estimate stachastic c}
\\ &&
+\mathbf{1}_{\{t_k>0\}} |U^{l+1-k}(t_k,x_k,v_{k-1})| d\Sigma_{k-1}(t_k) \label{2Estimate stachastic d}\\
&&+ \sum_{m=1}^{k-1}\mathbf{1}_{\{t_m>0\}}d\Sigma_m^r.\label{2Estimate stachastic e}
\end{eqnarray}

The estimate of $U^l$ in $L^{\infty}$ is as follows.
\begin{lem}\label{2LI Estimate by Stochacstic cycle}
Suppose that $||u_0||_{\infty},~\sup_{0\leq s\leq T}|r^{\lambda}(s)|_{\infty},~\sup_{0\leq s\leq T}||q^{\lambda}(s)||_{\infty}$ are bounded for fixed $T>0$. Then, there exists $C(k)>0$  such that the solution $U^{l+1}$ of (\ref{2MSimply problem jl diffusive Ref Bound}) satisfies that
\begin{eqnarray}\label{2Main Estimate in LI}
\sup_{0\leq s\leq T} ||U^{l+1}(s)||_{\infty} &\leq& \frac{1}{8} \max_{0\leq m\leq 2k} \sup_{0\leq s\leq T} ||U^{l-m}(s)||_{\infty} + C(k) \max_{1\leq m\leq 2k} \int_0^T ||U^{l-m}(s) ||_2 ds \notag\\
&& +C(k) \left[ ||u_0||_{\infty} + \sup_{0\leq s\leq T}|r^{\lambda}(s)|_{\infty} +  \sup_{0\leq s\leq T}||q^{\lambda}(s)||_{\infty}\right].
\end{eqnarray}
\end{lem}
\prof We start with $r$-contribution in (\ref{2Estimate stachastic 0}) and (\ref{2Estimate stachastic e}). Since $ d\sigma_m$ is a probability measure, $\int_{\mathcal{V}_m}d\sigma_m \leq 1$, from the definition (\ref{2Bounce of boundary condition}), the contribution of $r$ is bounded by
\begin{eqnarray}\label{2Contribution of r}
 |r^{\lambda}(t_1,x_1,v)| + \sum_{m=1}^{k-1} |r^{\lambda}(t_{m+1},x_{m+1},v_m)| \leq k \sup_{0\leq s \leq t}|r^{\lambda}(s)|_{\infty}.
\end{eqnarray}
We turn to the $q^{\lambda}$- contribution in (\ref{2Estimate stachastic 0}), (\ref{2Estimate stachastic b}) and (\ref{2Estimate stachastic c}). Since $ d\sigma_m$ is a probability measure, all the terms to $q^{\lambda}$ is bounded by
\begin{eqnarray}\label{2Contribution of g}
&&\int_0^t ||q^{\lambda}(\tau)||_{\infty} \left\{2 + \sum_{m=1}^{k-1} \int \Big[\mathbf{1}_{\{t_{m+1}\leq 0<t_m\}} + \mathbf{1}_{\{t_{m+1}> 0\}} \Big]\Pi_{m=1}^{k-1} d\sigma_m \right\}d\tau \notag\\
&& \leq
2k T \sup_{0\leq s\leq t}||q^{\lambda}(s)||_{\infty}.
\end{eqnarray}
Now, we consider the contribution of the initial data $u_0$ in (\ref{2Estimate stachastic 0}) and (\ref{2Estimate stachastic a}). It could be bounded by
\begin{eqnarray}\label{2Contribution of initial data}
 \mathbf{1}_{\{t_1 \leq 0\}} e^{-\int_0^{t} (\lambda+\Sigma)} ||u_0||_{\infty} + \int_{\Pi_{m=1}^{k-1}}\sum_{m=1}^{k-1} \mathbf{1}_{\{t_{m+1} \leq 0 <t_m\}} ||u_0||_{\infty}d\Sigma_m(0) \leq k ||u_0||_{\infty}.
\end{eqnarray}
From Lemma \ref{2Small for large-time bounce}, (\ref{2Estimate stachastic d}) can be bounded by
\begin{eqnarray}\label{2small for large bounces}
&&\int_{\Pi_{m=1}^{k-1}} \mathbf{1}_{\{t_k>0\}} |U^{l+1-k}(t_k,x_k,v_{k-1})| d\Sigma_{k-1}(t_k) \notag\\
 && \leq \sup_{0\leq s \leq t} ||U^{l+1-k}(s)||_{\infty} \int_{\Pi_{m=1}^{k-1}}  \mathbf{1}_{\{t_k>0\}} d\Sigma_{k-1}(t_k)
  \leq  \varepsilon \sup_{0\leq s \leq t} ||U^{l+1-k}(s)||_{\infty}.
\end{eqnarray}
From (\ref{2Contribution of r}), (\ref{2Contribution of g}), (\ref{2Contribution of initial data}) and (\ref{2small for large bounces}), we obtain an upper bound  that
\begin{eqnarray}\label{2Estimate L2-LI}
|U^{l+1}(t,x,v)| &\leq&   \mathbf{1}_{\{t_1 \leq 0\}}\int_0^t e^{-\int_0^{t-s} (\lambda+\Sigma)} | K U^l(\tau, x-(t-\tau)v,v)|d\tau \notag \\
 && + \mathbf{1}_{\{t_1 > 0\}} \int_{t_1}^t e^{-\int_0^{t-s} (\lambda+\Sigma)} |K U^l(\tau, x-(t-\tau)v,v)|d\tau \notag \\
 && + \sum_{m=1}^{k-1} \Big\{\int_0^{t_m} \mathbf{1}_{\{t_{m+1}\leq 0<t_m \}} |KU^{l-m}(\tau,X_{cl}(\tau),v_m)| \notag\\
&& + \int_{t_{m+1}}^{t_m}\mathbf{1}_{\{t_{m+1}>0\}}|KU^{l-m}(\tau,X_{cl}(\tau),v_m)|\Big\}d\Sigma_m(\tau)d\tau + A_l(t,x,v)
\end{eqnarray}
with $A_l(t,x,v)$ denotes
\begin{eqnarray}\label{2Error term for L2-LI}
A_l(t,x,v) =  \varepsilon \sup_{0\leq s \leq t} ||U^{l+1-k}(s)||_{\infty}+C_{k,T}\Big(||u_0||_{\infty}+ \sup_{0\leq s \leq t}\{|r^{\lambda}(s)|_{\infty}+ ||q^{\lambda}(s)||_{\infty}\}\Big).
\end{eqnarray}
Recall that the back-time cycle $(s,X_{cl}(s; t, x, v), v¡ä)$ denotes
$(t¡ä_1, x¡ä_1, v¡ä_1), (t¡ä_2, x¡ä_2, v¡ä_2 ), \cdot , (t¡ä_{m¡ä} , x¡ä_{m¡ä} , v¡ä_{m¡ä}), \cdots $, we now iterate (\ref{2Estimate L2-LI}) for $l-m$ times to get the representation for $U^{l-m}$ and then plug in $K U^{l-m}(s,X_{cl}(s), v_m)$ to obtain
\begin{eqnarray}\label{2Iterate scheme for L^2-LI}
&& |KU^{l-m}(s,X_{cl}(s), v_m)| \notag\\
&&\leq \int_{V} f(X_{cl}(s),v_m,v')|U^{l-m}(s,X_{cl}(s),v')|dv' \notag\\
 && \leq  \iint_{V\times V} \Big\{ \mathbf{1}_{\{t_1' \leq 0\}} \int_0^s +  \mathbf{1}_{\{t_1'> 0\}} \int_{t_1'}^s \Big\} e^{-\int_0^{s-s_1} (\lambda+\Sigma)} f(X_{cl}(s)-(s-s_1)v',v_m,v')\notag
  \\
 &&\hspace{5mm} \times f(X_{cl}(s)-(s-s_1)v',v',v'') |U^{l-1-m}(s_1, X_{cl}(s)-(s-s_1)v', v'')| ds_1dv'dv'' \notag\\
 && + \iint_{V\times V} \int_{\Pi_{m=1}^{k-1}\mathcal{V}_m'} \sum_{l'=1}^{k-1}e^{-\int_0^{s-t'_{l'}} (\lambda+\Sigma)}
 \Big\{ \int_0^{t_{l'}'}ds_1 \mathbf{1}_{\{t_{l'+1}' \leq 0<t'_{l'}\}} +  \int_{t_{l'}'}^{t_{l'-1}'} ds_1\mathbf{1}_{\{t_{l'}'> 0\}} \Big\} \notag\\
 &&\hspace{5mm} \times f(y_{l'},v_m,v')
 f( y_{l'},v'_{l'},v'') |U^{l-1-m-l'}(s_1, y_{l'} , v'')| d\Sigma_{l'}(s_1) dv'dv'' \notag\\
 && + \int_{V} f(X_{cl}, v_m, v')A_{l-1-m}(s,X_{cl}(s),v')dv',
  \end{eqnarray}
where $y_{l'}= x_{l'}'+(s-t_{l'}')v'_{l'}$.

The total contributions of $A_{l-m-1}$ in (\ref{2Estimate L2-LI}) are obtained via plugging (\ref{2Iterate scheme for L^2-LI}) with different
$l$  into (\ref{2Estimate L2-LI}). Since
$\int f(x,v,v')dv'<M_b$, the summation of all contributions of $A_{l-m-1}$ leads
to the bound
\begin{eqnarray} \label{2Estimate of A_l}
&&2 A_{l-1}(t)\int_0^t e^{-\int_0^{t-s}(\lambda+\Sigma)} ds + \max_{1\leq m\leq k-1} A_{l-m-1}(t)  \notag\\
&& \hspace{5mm} \times \int_{\Pi_{m=1}^{k-1}\mathcal{V}_m} \sum_{m=1}^{k-1} \left\{ \int_0^{t_m} \mathbf{1}_{\{t_{m+1}\leq 0 <t_m\}}+\int_{t_{m+1}}^{t_m} \mathbf{1}_{t_m>0} \right\} d\Sigma_m(s)ds +A_l(t)\notag\\
&& \leq C(k) \max_{0\leq m\leq k} A_{l-m}(t) \notag\\
&&\leq  \varepsilon \sup_{0\leq s \leq t} ||U^{l+1-m-k}(s)||_{\infty}+C_{k,T}\Big(||u_0||_{\infty}+ \sup_{0\leq s \leq t}|r^{\lambda}(s)|_{\infty}+ \sup_{0\leq s\leq t}||q^{\lambda}(s)||_{\infty}\Big)
\end{eqnarray}
To estimate the $U^{l-m-1}$ contribution, we separate $s - s_1 \leq \varepsilon$ and $s - s_1 \geq \varepsilon$. In the
first case, we use the fact $\int f(x,v,v')dv'<M_b $ and by (\ref{2Iterate scheme for L^2-LI}) to obtain the small contribution
\begin{eqnarray}\label{2Small time estimate}
\varepsilon \max_{1\leq l'\leq k} \sup_{0\leq s_1\leq s} ||U^{l-l'-m}(s_1)||_{\infty}.
\end{eqnarray}
Now let us treat the case of $s-s_1 \geq  \varepsilon$. It can be easily check that $\int_{ V} f(x,v,v')dv' \leq M_b <\infty$. For $N\gg1$, there is some constant $p(N)$ such that
\begin{eqnarray*}
f_p(x,v'_{l'},v'') = f \mathbf{1}_{|f|\leq p(N)}.
 \end{eqnarray*}
It satisfies that $\sup_{x,v'_{l'}}\int_{V} |f(x,v'_{l'},v'')- f_p(x,v'_{l'},v'') |  dv''  \leq \frac{1}{N} $. For any $l'$. We split $f(x,v'_{l'},v'' )
=\{f(x,v'_{l'},v'')-f_p(x,v'_{l'},v'')\}+f_p (x,v'_{l'},v'')$. Notice that $\int f(x,v,v')dv'<M_b $, the first difference of $f(x,v'_{l'},v'' )$ leads to a small contribution in (\ref{2Iterate scheme for L^2-LI})
\begin{eqnarray}\label{2Contribution on singular part}
\frac{1}{N} \max_{1\leq l'\leq k} \sup_{0\leq s_1\leq s}||U^{l-m-l'} (s_1)||_{\infty}.
\end{eqnarray}
For the remainder main contribution of $k_p(y_{l'},v'_{l'},v'')$
 the change of variable $y_{l'}= x_{l'}'+(s-t_{l'}')v'_{l'} $ ($x_{l'}'$ do not depend on $v_{l'}'$) satisfies
 $|\frac{dy}{dv'}|\geq \varepsilon^3$ for $s-s_1\geq \varepsilon$. Notice that $|f_p(y_{l'},v,v')|\leq p(N)$, the remained part can be estimated by
\begin{eqnarray} \label{2Main part L2-LI}
&& p(N)^2 \int_{v''} \int_{\mathcal{V}_{l'}'} e^{-\int_0^{s-t'_{l'}} (\lambda+\Sigma)}  |U^{l-m-l'}(s_1, y_{l'} , v'')| M_w(v'_{l'}) |n(x'_{l'})\cdot v'_{l'}| dv'_{l'}dv''\notag\\
&& \leq \frac{p(N)^2}{\varepsilon^3} \iint_{\Omega \times V} |U^{l-m-l'}(s_1, y_{l'} , v'')| dydv'' \lesssim_{N,\varepsilon} ||U^{l-m-l'}(s_1)||_2.
\end{eqnarray}
The estimates (\ref{2Small time estimate}), (\ref{2Contribution on singular part}), (\ref{2Main part L2-LI}) give a bound for (\ref{2Iterate scheme for L^2-LI}) as
\begin{eqnarray}\label{2Estimate of Iterate Scheme}
 |KU^{l+1-m}(s,X_{cl}(s), v_m)|
 &\leq&   [2 \varepsilon + \frac{C(k)}{N}] \max_{1\leq m\leq 2k} \sup_{0\leq s_1\leq s} ||U^{l-m}(s_1)||_{\infty} \notag\\
 &&+ C_{k,T}\Big(||u_0||_{\infty}+ \sup_{0\leq s \leq t}|r^{\lambda}(s)|_{\infty}+ \sup_{0\leq s\leq t}||q^{\lambda}(s)||_{\infty}\Big) \notag \\
 && + C_{\varepsilon,N} \max_{1\leq m\leq 2k} \int_0^s  ||U^{l-m}(s_1)||_2 ds_1\notag\\
 &\equiv&  B(s).
\end{eqnarray}
By plugging back (\ref{2Estimate of A_l}) and (\ref{2Estimate of Iterate Scheme}) into (\ref{2Estimate L2-LI}), we have the bounded $|U^{l+1}(t,x,v)|$ by
\begin{eqnarray}\label{2Bounded of U^l}
&&B(t) \left\{ \mathbf{1}_{t_1 \leq 0} \int_0^t e^{-\int_0^{t-s}(\lambda + \Sigma)}ds +\mathbf{1}_{t_1 > 0} \int_{t_1}^t e^{-\int_0^{t-s}(\lambda + \Sigma)}ds \right\} \notag\\
&& \hspace{3mm} + B(t) e^{-\int_0^{t-t_1}(\lambda + \Sigma)} \sum_{m=1}^{k-1}\left\{ \int_0^{t_m}\mathbf{1}_{t_{m+1} \leq 0<t_m} + \int_{t_{m+1}}^{t_m}\mathbf{1}_{t_m > 0}
 \right\}d\Sigma_m(s)ds + A_l(t)\notag\\
 &&\leq   [2 \varepsilon + \frac{C(k)}{N}] \max_{1\leq m\leq 2k} \sup_{0\leq s_1\leq s} ||U^{l-m}(s_1)||_{\infty} \notag\\
 &&\hspace{5mm}+ C_{k,T}\Big(||u_0||_{\infty}+ \sup_{0\leq s \leq t}|r^{\lambda}(s)|_{\infty}+ \sup_{0\leq s\leq t}||q^{\lambda}(s)||_{\infty}\Big)\notag\\
  && \hspace{5mm}+ C_{\varepsilon,N} \max_{1\leq m\leq 2k} \int_0^s  ||U^{l-m}(s_1)||_2 ds_1.
\end{eqnarray}
We then conclude the proof of (\ref{2Main Estimate in LI}) by choosing $\varepsilon $ sufficiently small and $N$ large sufficiently large. The proof of Lemma \ref{2LI Estimate by Stochacstic cycle} is completed. $\hfill\square$\\

Before going to prove the $L^{\infty}$ estimate of $U^{l}$, we prove a standard result similar to \cite{[GKTT1]} with more precise estimate. This lemma can be crucial to get the bound of $U^l$ as well as its convergence.
\begin{lem}\label{2Main result for sequaences}
Suppose $b_i \geq 0,~D\geq 0,~~0\leq \eta\leq 1$ and $B_i=\max\{b_i,\cdots, b_{i-(k-1)}\}$ for fixed $k \in \mathds{N}$.
 \begin{eqnarray}\label{2Assumptions of the sequence}
~b_{l+1} \leq \frac{1}{8} B_l + D \eta^l\quad \text{for~all}\quad l\geq k.\end{eqnarray}
Then, for all $p\geq 1$ and $1\leq m\leq k$, it holds that
\begin{eqnarray}\label{2Estimates of the sequence1}
\begin{array}{lll}
B_{ik+m} &\leq & \displaystyle{\Big[ \frac{1}{8^{i-1}} +  \Big(\frac{7\eta^{m}}{8^{i-1}}+\frac{7\eta^{k+m}}{8^{i-2}} +\cdots+\frac{7\eta^{(i-2)k+m}}{8}\Big)\Big]\max\{B_{k}, \frac{8}{7}D\eta^{k+1}\}},\\
b_{i k+m+1} &\leq& \displaystyle{ \Big[ \frac{1}{8^i} +  \Big(\frac{7\eta^{m}}{8^i}+\frac{7\eta^{k+m}}{8^{i-1}} +\cdots+\frac{7\eta^{(i-1)k+m}}{8}\Big)\Big]\max\{B_{k}, \frac{8}{7}D\eta^{k+1}\}.}
\end{array}
\end{eqnarray}
In particular, when $\eta=1$, it holds that, for all $l\geq 1$,
\begin{eqnarray}\label{2Estimates of the sequence2}
b_{l} \leq \max\{B_{k}, \frac{8}{7}D\}.
\end{eqnarray}
\end{lem}

%and for all $p\geq 1$
%\begin{eqnarray}\label{Estimates of the sequence2}
%B_{(p+1)k} \leq \frac{1}{8}B_{pk} +\frac{8}{7}D.\end{eqnarray}
 \prof%-(\ref{Estimates of the sequence2})
 We will prove it by induction with responding to $i$. From the definition of $B_{k}$, we know that $b_1,\cdots,b_{k} \leq B_{k} $. Then
\begin{eqnarray*}
B_k \leq \max\{B_{k}, \frac{8}{7}D\eta^{k+1}\},
\quad b_{k+1} \leq  \frac{1}{8} B_{k} + D\eta^{k+1}
 \leq \max\{B_{k}, \frac{8}{7}D\eta^{k+1}\},
\end{eqnarray*}
and
\begin{eqnarray*}
  B_{k+1} &=&\max\{b_{k+1},\cdots,b_2\} \leq \max\{B_{k}, \frac{8}{7}D\eta^{k+1}\},\\
 b_{k+2} &\leq & \frac{1}{8} B_{k+1} + D\eta^{k+2}
 \leq (\frac{1}{8}+\frac{7\eta}{8})\max\{B_{k}, \frac{8}{7}D\eta^{k+1}\}\leq \max\{B_{k}, \frac{8}{7}D\eta^{k+1}\}.
\end{eqnarray*}
Here we have used the fact that $ D\eta^{k+1} \leq \frac{7}{8}\max\{B_{k}, \frac{8}{7}D\eta^{k+1}\}  $.

Similarly, for all $1\leq m\leq k$, it derives to
\begin{eqnarray*}
B_{k+m} &=& \max\{b_{k+m},\cdots,b_{m+1}\} \leq \max\{B_{k}, \frac{8}{7}D\eta^{k+1}\}\\
b_{k+m+1} &\leq &  (\frac{1}{8}+\frac{7\eta^{m}}{8})\max\{B_{k}, \frac{8}{7}D\eta^{k+1}\}\leq \max\{B_{k}, \frac{8}{7}D\eta^{k+1}\}.
\end{eqnarray*}
Since $\eta^m \leq \eta^{m-1}$ for $\eta\leq 1$, we can also derive, for all $1\leq m\leq k$,
\begin{eqnarray*}
B_{2k+1} & = & \max\{b_{2k+1},\cdots, b_{k+2}\}\leq (\frac{1}{8} +\frac{7\eta}{8})\max\{B_{k}, \frac{8}{7}D\eta^{k+1}\},\\
b_{2k+2}&\leq &\frac{1}{8}B_{2k+1}+D\eta^{2k+2} \leq (\frac{1}{8^2} +\frac{7\eta}{8^2}+\frac{7\eta^{k+1}}{8})\max\{B_{k}, \frac{8}{7}D\eta^{k+1}\},\\
%B_{2k+2} & = & \max\{b_{2k+2},\cdots, b_{k+3}\}\leq (\frac{1}{8} +\frac{7\eta^2}{8})\max\{B_{k}, \frac{8}{7}D\eta^{k+1}\},  \\
\cdots\\
B_{2k+m} & = & \max\{b_{2k+m},\cdots, b_{k+m+1}\}\leq (\frac{1}{8} +\frac{7\eta^m}{8})\max\{B_{k}, \frac{8}{7}D\eta^{k+1}\},\\
b_{2k+m+1}&\leq &\frac{1}{8}B_{2k+m}+D\eta^{2k+m+1} \leq (\frac{1}{8^2} +\frac{7\eta^{m}}{8} +\frac{7\eta^{k+m}}{8})\max\{B_{k}, \frac{8}{7}D\eta^{k+1}\}.
\end{eqnarray*}
It means that (\ref{2Estimates of the sequence1}) is valid for $i=2$.

Suppose that (\ref{2Estimates of the sequence1}) holds for $i=p$, that is,
for all $1\leq m\leq k$
\begin{eqnarray*}
 B_{pk+m} &\leq & \Big[ \frac{1}{8^{p-1}} +  \Big(\frac{7\eta^{m}}{8^{p-1}}+\frac{7\eta^{k+m}}{8^{p-2}} +\cdots+\frac{7\eta^{(p-2)k+m}}{8}\Big)\Big]\max\{B_{k}, \frac{8}{7}D\eta^{k+1}\},\\
 b_{p k+m+1} &\leq& \Big[ \frac{1}{8^p} +  \Big(\frac{7\eta^{m}}{8^p}+\frac{7\eta^{k+m}}{8^{p-1}} +\cdots+\frac{7\eta^{(p-1)k+m}}{8}\Big)\Big]\max\{B_{k}, \frac{8}{7}D\eta^{k+1}\}.\\
\end{eqnarray*}
Then, from (\ref{2Assumptions of the sequence}), we have
\begin{eqnarray*}
B_{(p+1)k+1} &=&\max\{b_{(p+1)k+1},\cdots, b_{pk+2} \} \\
&\leq& \Big[ \frac{1}{8^p} +  \Big(\frac{7\eta}{8^p}+\frac{7\eta^{k+1}}{8^{p-1}} +\cdots+\frac{7\eta^{(p-1)k+1}}{8}\Big)\Big]\max\{B_{k}, \frac{8}{7}D\eta^{k+1}\},\\
b_{(p+1)k+2} &\leq &  \Big[ \frac{1}{8^{p+1}} +  \Big(\frac{7\eta}{8^{p+1}}+\frac{7\eta^{k+1}}{8^{p}} +\cdots+\frac{7\eta^{pk+1}}{8}\Big)\Big]\max\{B_{k}, \frac{8}{7}D\eta^{k+1}\},\\
\cdots\\
B_{(p+1)k+m} &=&\max\{b_{(p+1)k+m},\cdots, b_{pk+m+1} \} \\
&\leq& \Big[ \frac{1}{8^{p}} +  \Big(\frac{7\eta^m}{8^p}+\frac{7\eta^{k+m}}{8^{p-1}} +\cdots+\frac{7\eta^{(p-1)k+m}}{8}\Big)\Big]\max\{B_{k}, \frac{8}{7}D\eta^{k+1}\},\\
b_{(p+1)k+m+1} &\leq &  \Big[ \frac{1}{8^{p+1}} +  \Big(\frac{7\eta^m}{8^{p+1}}+\frac{7\eta^{k+m}}{8^{p}} +\cdots+\frac{7\eta^{pk+m}}{8}\Big)\Big]\max\{B_{k}, \frac{8}{7}D\eta^{k+1}\}.
\end{eqnarray*}
This implies that (\ref{2Estimates of the sequence1}) holds for $i=p+1$ and (\ref{2Estimates of the sequence1}) is true for all $p\geq 1$. (\ref{2Estimates of the sequence2}) is a consequence of (\ref{2Estimates of the sequence1}). The proof of Lemma \ref{2Main result for sequaences}$\hfill\square\\$

Now, we consider the $L^{\infty}$ estimate of $U^{l}$. For simplicity, we denote
\begin{eqnarray*}
E = ||u_0||_{\infty}+ \sup_{0\leq s \leq T}|r^{\lambda}(s)|_{\infty}+ \sup_{0\leq s\leq T}||q^{\lambda}(s)||_{\infty}.
\end{eqnarray*}
\begin{lem}\label{2LI Estimates of the solution}
Suppose that $||u_0||_{\infty},~\sup_{0\leq s\leq T}|r^{\lambda}(s)|_{\infty},~\sup_{0\leq s\leq T}||q^{\lambda}(s)||_{\infty}$ are bounded for fixed $T>0$. Then, there exists $C(k)>0$  such that the solution $U^{l+1}$ of (\ref{2MSimply problem jl diffusive Ref Bound}) satisfies that
\begin{eqnarray}\label{2Main Estimates in LI}
\sup_{0\leq s\leq T} ||U^{l+1}(s)||_{\infty} \lesssim  \max_{1\leq m\leq 2k} \int_0^T ||U^{l-m}(s) ||_2 ds + E.
\end{eqnarray}
\end{lem}
\prof {\it Step 1} Let $l \rightarrow \infty$. From Lemma \ref{2LI Estimate by Stochacstic cycle} and the estimate (\ref{2Basic Uniform L2 Estimate}), we can obtain
\begin{eqnarray*}
\sup_{0\leq s\leq T}||U^{l+1}(s)||_{\infty} \leq \frac{1}{8} \max_{0\leq m\leq 2k} \sup_{0\leq s\leq T}||U^{l-m}(s)||_{\infty} + C(k) \Big(\mathcal{E}+E\Big).
\end{eqnarray*}
Set $b_l=\sup_{0\leq s\leq T}||U^{l}(s)||_{\infty}$, $ B_{l} = \max_{1\leq m\leq 2k} \sup_{0\leq s\leq T}||U^{l-m-1}||_{\infty} $, $\eta=1$ and $D=C(k) \Big(\mathcal{E}+E\Big)$, then,  (\ref{2Estimates of the sequence1}) gives that
\begin{eqnarray}\label{2Estimate of all sequence}
&&\sup_{0\leq s\leq T}||U^{l}(s)||_{\infty} \leq \max\Big\{B_{2k},  \frac{8}{7}  C(k) \Big(\mathcal{E}+E\Big)\Big\}.
\end{eqnarray}
Now, we give the estimate of $ B_{2k}$.  Since  $||Ku||_{L^{\infty}} \leq M_b ||u||_{\infty}$ and $|\mathcal{P}_{\gamma}u|_{\infty} \leq ||u||_{\infty}$, from the iterate scheme (\ref{2MSimply problem jl diffusive Ref Bound}), we have
\begin{eqnarray*}
 |U^{l+1}(t,x,v)|  & \leq & \mathbf{1}_{\{t_1 \leq 0\}} e^{-\int_0^{t} (\lambda+\Sigma)} |u_0(x-tv,v)| \notag\\
 &&+\mathbf{1}_{\{t_1 \leq 0\}}\int_0^t e^{-\int_0^{t-\tau} (\lambda+\Sigma)} |[K U^l+q^{\lambda}](\tau, x-(t-\tau)v,v)|d\tau \notag \\
 && + \mathbf{1}_{\{t_1 > 0\}} \int_{t_1}^t e^{-\int_0^{t-\tau} (\lambda+\Sigma)} |[K U^l+q^{\lambda}](\tau, x-(t-\tau)v,v)|d\tau \notag \\
&&+ \mathbf{1}_{\{t_1>0\}} e^{-\int_0^{t-t_1} (\lambda+\Sigma)} \left|[(1-\frac{1}{j})\mathcal{P}_{\gamma}U^l + r^{\lambda}](t_1,x_1,v)\right|\notag\\
&\leq & C_1 \sup_{0\leq s \leq T}||U^l(s)||_{\infty} + 2 E.
\end{eqnarray*}
So, for fixed $k$, we get iterate a bound for $i\leq 2k$ to obtain
\begin{eqnarray*}
 \sup_{0\leq s\leq T}||U^{m+1}(s)||_{\infty} &\leq & C_1\sup_{0\leq s \leq T}||U^m(s)||_{\infty} + 2 E \\
 &\leq & C_1^2\sup_{0\leq s \leq T}||U^{m-1}(s)||_{\infty} + 2(1+C_1) E\\
 &\leq &\cdots \\
 &\leq & C_1^{m+1}\sup_{0\leq s \leq T}||U^0(s)||_{\infty} + 2[C_1^{m}+\cdots+C_1+1] E
\end{eqnarray*}
This inequality leads to
\begin{eqnarray}\label{Estimate of former 2k}
B_{2k} = \max_{0\leq m\leq 2k} \sup_{0\leq s\leq T}||U^{2k-m}||_{\infty} \leq 2 [C_1^{2k} +\cdots+ C_1+1] E <\infty.
\end{eqnarray}
Because $T<\infty $, $\Omega$ and $V$ are bounded domains, we know that $\mathcal{E} \lesssim_{T,\Omega,V} E$. From (\ref{2Estimate of all sequence}), for any $l\geq 1$, one get the following uniformly $L^{\infty}$ bound of $U^{l}$
\begin{eqnarray}
 \sup_{0\leq s \leq T}||U^{l}(s)||_{\infty} \leq \max\Big\{B_{2k}, \frac{8}{7}  C(k) \Big(\mathcal{E}+D \Big)\Big\} \lesssim_{T,\lambda,\Omega,V} E.
\end{eqnarray}
This gives the uniform estimate of the sequences $U^{l}$.

On the other hand, the difference $V^{l+1}=: U^{l+1}-U^l$ satisfies (\ref{2Difference of sequaence}). From Lemma \ref{2LI Estimate by Stochacstic cycle}, we get
\begin{eqnarray*}
\sup_{0\leq s\leq T} ||V^{l+1} (s)||_{\infty}  \leq  \frac{1}{8} \max_{0\leq m\leq 2k}\sup_{0\leq s \leq T}||V^{l-m}(s)||_{\infty} +C(k) \max_{0\leq m\leq 2k} \int_0^T ||V^{l-m}(s) ||_2ds.
\end{eqnarray*}
From the estimate (\ref{2Basic L2 Estimate of Difference}), we know there is $\eta_{\lambda,j}$ such that
\begin{eqnarray}
||V^l(s)||_2 \leq \eta_{\lambda,j}^l \int_0^t (|V^1(s)|_{2,+}^2 +\lambda ||V^1(s)||_2^2)ds.
\end{eqnarray}
Thus there exists constant $C_{k,T}$ such that
\begin{eqnarray*}
 C(k) \max_{0\leq m\leq 2k} \int_0^T ||V^{l-m}(s) ||_2ds \leq C_{k,T} (\eta_{\lambda,j})^{l/2}.
\end{eqnarray*}
Hence, let $l=2pk+m~(1\leq m\leq 2k)$, (\ref{2Estimates of the sequence1}) gives that
\begin{eqnarray*}
 \sup_{0\leq t\leq T}||V^l(t)||_{\infty} &\leq& \Big[ \frac{1}{8^i} +  \Big(\frac{7\eta_{\lambda,j}^{m/2}}{8^i}+\cdots+\frac{7\eta_{\lambda,j}^{[2(i-1)k+m]/2}}{8}\Big)\Big]\\
 &&\times \max\{\max_{1\leq l\leq 2k}\sup_{0\leq t\leq T}||V^l(t)||_{\infty}, \frac{8}{7}C_{k,T}\eta_{\lambda,j}^{(2k+1)/2}\}.
\end{eqnarray*}
Because $\eta_{\lambda,j}<1$, we know that
\begin{eqnarray}
 \sum_{l=1}^{\infty}\sup_{0\leq t\leq T}||V^l(t)||_{\infty} \leq C \max\{\max_{1\leq l\leq 2k}\sup_{0\leq t\leq T}||V^l(t)||_{\infty}, \frac{8}{7}C_{k,T}\eta_{\lambda,j}^{(2k+1)/2}\} < \infty.
 \end{eqnarray}
It means that $\{U^l\}_{l=1}^{\infty}$ is a Cauchy series. Hence, there is a limit solution $U^l \rightarrow U_j$. $U_j$ is the solution of (\ref{2Reduced boundary problem}). Thus, the difference $U^{l+1}-U_j$ satisfies
\begin{eqnarray}
(\partial_t +v\cdot \nabla+\lambda+\Sigma)(U^{l+1}-U_j) = K(U^l-U_j),~~~(U^{l+1}-U_j)|_{t=0} =0,
\end{eqnarray}
with $(U^{l+1}-U_j)|_{\gamma_-} =(1-\frac{1}{j})\mathcal{P}_{\gamma}(U^l-U_j)$. By the same argument as above, we can yield that $\sup_{0\leq t\leq T} ||(U^l-U_j)(s)||_{L^{\infty}} \rightarrow 0 $ as $l\rightarrow \infty$.\\

{\it Step 2}. We take $j \rightarrow \infty$. Let $U_j$ be the solution to (\ref{2Reduced boundary problem}). Lemma \ref{2LI Estimate by Stochacstic cycle} implies that
\begin{eqnarray*}
\sup_{0\leq s\leq T}||U_j(s)||_{\infty} &\leq& \frac{1}{8}\sup_{0\leq s\leq T}||U_j(s)||_{\infty}  +C(k) \int_0^T||U_j(s)||_2ds \\
&&+ C(k) \left( ||u_0||_{\infty}+ \sup_{0\leq s\leq T}(|r^{\lambda}(s)|_{\infty}+||q^{\lambda}(s)||_{\infty})\right).
\end{eqnarray*}
Therefore, by an induction over j,
\begin{eqnarray*}
\sup_{0\leq s\leq T}||U_j(s)||_{\infty} \leq C(k) \int_0^T||U_j(s)||_2ds
+ C(k) E.
\end{eqnarray*}
Since $\int_0^T||U_j(s)||_2ds $ is bounded from Step 2 of Lemma \ref{2Existence with Collision operator}, this implies that $\sup_{0\leq s\leq T}||U_j(s)||_{\infty}$ is uniformly bounded and we obtain the solution $U$. Taking the difference, we have
\begin{eqnarray}
(\partial_t +v\cdot\nabla_x +\lambda+\Sigma)[U_j-U]=K[U_j-U], \quad [U_j-U]|_{t=0}=0,
\end{eqnarray}
with the boundary condition $[U_j-U]|_{\gamma_-} = \mathcal{P}_{\gamma} [U_j-U] +\frac{1}{j}\mathcal{P}_{\gamma} U_j $. We regard $\frac{1}{j}\mathcal{P}_{\gamma} U_j $ as $r^{\lambda}$ in Lemma \ref{2LI Estimate by Stochacstic cycle} implies that
\begin{eqnarray}
\sup_{0\leq s\leq T}||[U_j-U](s)||_{\infty} \leq C(k)\int_0^T||[U_j-U](s)||_2ds +\frac{C(k)}{j}\sup_{0\leq s\leq T} ||U_j||_{\infty},
\end{eqnarray}
which goes to zero as $j \rightarrow \infty$. We obtain  $L^{\infty}$ solution $U$ to (\ref{2Modified problem}). $\hfill\square\\$

\noindent{\bf The proof of Theorem \ref{1Existence of the solution}}:  Let $\lambda$ is large enough such that (\ref{2range of parameter}) is satisfied. The existence of the solution $U$ of (\ref{2Modified problem}) and the estimate (\ref{2Estimate of Modified solution}) come from Lemma \ref{2LI Estimates of the solution} immediately. It is exactly the result of Proposition \ref{2Solution for the modified problem}. Form Proposition \ref{2equilibrium}, the solution $u$ of (\ref{1Main problem})-(\ref{1Boundary condition}) is also obtained. Furthermore, (\ref{2Estimate of solution}) holds true for all $0\leq t\leq T$. Proposition \ref{2Solution for the original problem} is proved. It is exactly the solution to (\ref{1Main problem}) with (\ref{1Boundary condition}). Moreover, the estimate (\ref{1Uniform bound of solution}) is valid. This completes the proof of Theorem \ref{1Existence of the solution}. $\hfill\square$

%%%%%%%%%%%%%%%%%%%%%%%%%%%%%%%%%%%%%%%%%%%%%%%%%%%%%%%%%%%%%%%%%%%%%%%%%%%%%%%%%%%%%%%%%%%%%%%%%%%%%%%%%%%%%%%%%%%%%%%%%%%%%%%%%%%%%%%%%%%%%%%%%%%%
% CHAPTER                                                                                                                                          %
%  3                                                                                                                                               %
%                                                                                                                                                  %
%%%%%%%%%%%%%%%%%%%%%%%%%%%%%%%%%%%%%%%%%%%%%%%%%%%%%%%%%%%%%%%%%%%%%%%%%%%%%%%%%%%%%%%%%%%%%%%%%%%%%%%%%%%%%%%%%%%%%%%%%%%%%%%%%%%%%%%%%%%%%%%%%%%%

\section{BV regularity of solution}
\setcounter{section}{3} \setcounter{equation}{0}
In this section, we construct an open covering of the singular set $\mathfrak{S}_B$, which is crucial to establish a
smooth approximation function that excludes the open covering of $\mathfrak{S}_B$. In particular, the measure of this singular set could be
sufficiently small from Lemma \ref{3Neighborhood of the singular set properties} and Proposition \ref{3Modified cover function}.

\subsection{The neighborhood of singular set}
In the following, we construct the neighborhoods of the singular set, which is similar to the Boltzmann equation constructed in \cite{[GKTT2]}. For the completeness of this paper, we present the details here.
\begin{lem}\label{3Neighborhood of the singular set properties}
 For $0\leq \varepsilon \leq \varepsilon_1 \ll 1$, we construct an open set $\mathcal{O}_{\varepsilon,\varepsilon_1} \subset \Omega \times V$ such that
 \begin{eqnarray}\label{3Neighborhood of the singular set}
 \mathfrak{S}_B \subset \mathcal{O}_{\varepsilon,\varepsilon_1}.
\end{eqnarray}
There exists $C_* =C_* (\Omega) \gg 1$ such that for any $0< \varepsilon \leq \varepsilon_1 \ll 1$
\begin{eqnarray}\label{3Neighborhood of the closed set}
\overline{\mathcal{O}_{\varepsilon,\varepsilon_1}} \subset \mathcal{O}_{\varepsilon,~C_* \varepsilon}.
\end{eqnarray}
Moreover, there exist $C_1 = C(\Omega, C_*)>0,$ $~C_2=C_2(\Omega,C_*)>0 $ such that
\begin{eqnarray}\label{3Neighborhood set measure}
\iint_{\Omega \times V}\mathbf{1}_{\mathcal{O}_{\varepsilon,C_*\varepsilon}}(x,v) dvdx \leq C_1 \varepsilon,
\end{eqnarray}
and
\begin{eqnarray}\label{3Neighborhood distance}
\text{dist}(\overline{\Omega} \times V \setminus \mathcal{O}_{\varepsilon, C_*\varepsilon},~\mathfrak{S}_B) >C_2 \varepsilon.
\end{eqnarray}
 \end{lem}

\noindent\prof {\it Construction of $\mathcal{O}_{\varepsilon,\varepsilon_1}$}.

{\it Step 1. Decomposition of $\partial \Omega$ } Let us fix $\theta>0$ which will be chosen later. Since the boundary
$\partial \Omega$ is locally a graph of smooth functions, from the finite covering theorem, there exists a finite number $M_{\Omega, \theta} $ of small open balls
$\mathcal{U}_1, \mathcal{U}_2, \cdots , \mathcal{U}_{M_{\Omega, \theta}} \subset \mathbb{R}^3 $ with $diam(U_m) < 2¦Ä$ for all $m$, such that
 \begin{eqnarray}\label{3Composition of Boundary}
\partial \Omega \subset \bigcup_{m=1}^{M_{\Omega,\theta}} [\mathcal{U}_m \cap \partial \Omega] ~~\text{with}~~M_{\Omega,\theta}=O(\theta^{-2}),
 \end{eqnarray}
and for every $m$, inside $\mathcal{U}_m$ the boundary $\mathcal{U}_m \cap  \partial \Omega$
 is exactly described by a smooth function $\eta_m$ defined on a (small) open set $\mathcal{A}_m \subset \mathbb{R}^2$. Up to rotations and translations and reducing the size of the ball $\mathcal{U}_m$) we will always assume that
 \begin{eqnarray}\label{3Piece of boundary}
 \mathcal{U}_m \cap \partial \Omega &=& \{ (x_1,x_2,\eta_m(x_1,x_2))\in \mathcal{A}_m \times \mathbb{R}\},\\
 \mathcal{U}_m \cap \Omega &=& \{ (x_1,x_2,x_3)\in \mathcal{A}_m \times \mathbb{R}: x_3> \eta_m(x_1,x_2)\}
 \end{eqnarray}
 with
 \begin{eqnarray*}
 (0,0) \in \mathcal{A}_m \subset [-\theta,\theta]\times [-\theta,\theta],~~~ \partial_1 \eta_m (0,0) =\partial_2 \eta_m(0,0)=0.
 \end{eqnarray*}
 Therefore, the unit out normal vector at $(0,0,\eta_m(0,0))$ is
 \begin{eqnarray*}
  n(0,0,\eta_m(0,0))=\frac{(\partial_1 \eta_m(0,0),\partial_2 \eta_m(0,0),-1)}{ \sqrt{1+|\partial_1 \eta_m(0,0)|^2+|\partial_2 \eta_m(0,0)|^2}}=(0,0,-1).
 \end{eqnarray*}
 Recall that $\partial \Omega$ is locally $C^2$. Then we can choose $\theta > 0$ small enough to satisfy for all
$m\in \{1,\cdots ,M_{\Omega,\theta}\}$ such that, for all $(x_1,x_2) \in \mathcal{A}_m$,
\begin{eqnarray} \label{3Error of boundary function}
\sum_{i=1}^2|\partial_i \eta_m(x_1,x_2) -\partial_i \eta_m(0,0)| =\sum_{i=1}^2|\partial_i \eta_m(x_1,x_2)| \leq \frac{1}{8},
~~ \sum_{i,j=1}^2 |\partial_{ij} \eta_m(x_1,x_2)| \leq C_{\eta}.
 \end{eqnarray}
Now we define the lattice point on $\mathcal{A}_m$ as
 \begin{eqnarray}\label{3Lattice points} c_{m,i,j,\varepsilon} =(\varepsilon i,~\varepsilon j)~~\text{for}~~-N_{\varepsilon}\leq i,j\leq N_{\varepsilon}=O( \varepsilon^{-1} \theta). \end{eqnarray}
Then we define the $(i, j)$-rectangular $\mathcal{R}_{m,i,j,\varepsilon,\varepsilon_1}$ which is centered at $c_{m,i,j,\varepsilon}$ and whose side is $2\varepsilon_1$:
\begin{eqnarray}\label{3Rectangular sets}
\mathcal{R}_{m,i,j,\varepsilon,\varepsilon_1} =\Big\{ (x_1,x_2)~:~\varepsilon i-\varepsilon_1 <x_1 <\varepsilon i+\varepsilon_1,~
\varepsilon j-\varepsilon_1 < x_2 <\varepsilon j-\varepsilon_1 \Big\} \cap \mathcal{A}_m.
\end{eqnarray}
Note that if $\varepsilon_1 \geq  \varepsilon$ then $\varepsilon i-\varepsilon_1$ is an open covering of $\mathcal{A}_m$, i.e.
\begin{eqnarray}\label{3Cover of variable space}
\mathcal{A}_m \subset \bigcup_{-N_{\varepsilon}\leq i,j\leq N_{\varepsilon} } \mathcal{R}_{m,i,j,\varepsilon,\varepsilon_1} ~~\text{with}~~N_{\varepsilon}=O( \varepsilon^{-1} \theta).
\end{eqnarray}
For each rectangle we define the representative outward normal
\begin{eqnarray*}
  n_{m,i,j,\varepsilon}=\frac{(\partial_1 \eta_m(c_{m,i,j,\varepsilon}),\partial_2 \eta_m(c_{m,i,j,\varepsilon}),-1)}{ \sqrt{1+|\partial_1 \eta_m(c_{m,i,j,\varepsilon})|^2+|\partial_2 \eta_m(c_{m,i,j,\varepsilon})|^2}}.
 \end{eqnarray*}
Let $(\hat{x}_{1,m,i,j,\varepsilon}, ~ \hat{x}_{2,m,i,j,\varepsilon}) \in \mathbb{S}^2$ be an orthonormal basis of the tangent space of $\partial\Omega$
 at $(c_{m,i,j,\varepsilon}, \eta_m(c_{m,i,j,\varepsilon})$. Remark that the three vectors $(\hat{x}_{1,m,i,j,\varepsilon}, ~ \hat{x}_{2,m,i,j,\varepsilon}, n_{m,i,j,\varepsilon}) $ is an orthonormal basis of $\mathbb{R}^3$  for each $m, i, j, \varepsilon$.\\

{\it Step 2. Decomposition of $ \Omega \times V$ }  We split the tangent velocity space at $(c_{m,i,j,\varepsilon},~\eta_{m}(c_{m,i,j,\varepsilon}))\in \partial \Omega$ as
\begin{eqnarray*}
\{v \in V: v\cdot n_{m,i,j,\varepsilon}=0\} \subseteq \bigcup_{l=0}^{L_{\varepsilon}} \Theta_{m,i,j,\varepsilon,\varepsilon_1,l}, ~~\text{with}~~L_{\varepsilon}=O(\frac{1}{\varepsilon}),
\end{eqnarray*}
where
\begin{eqnarray} \label{3Decomposition of velocity space}
 && \Theta_{m,i,j,\varepsilon,\varepsilon_1,l}
 : =  \Big\{ r_v \cos\theta_v \cos \phi_v \hat{x}_{1,m,i,j,\varepsilon}+  r_v \sin\theta_v \cos \phi_v \hat{x}_{2,m,i,j,\varepsilon}+ r_v \sin \phi_v n_{m,i,j,\varepsilon} \in V: \notag \\
&& \hspace{3cm}|r_v\sin\phi_v|\leq 8C_{\eta} \varepsilon_1 \max\{r_v,1\},~ ~|\theta_v -\varepsilon l|\leq \varepsilon_1 ~~\text{for}~~r_v \geq 0
 \Big\},
\end{eqnarray}
with the constant $C_{\eta}>0$ form (\ref{3Error of boundary function}).

Remark that for $\varepsilon_1 \geq \varepsilon, $
\begin{eqnarray}\label{3Joint set of the angular}
 \bigcup_{l=0}^{L_{\varepsilon}} \Theta_{m,i,j,\varepsilon,\varepsilon_1,l}=
 \Big\{v\in V: |v\cdot n_{m,i,j,\varepsilon}|\leq 8C_{\eta}\varepsilon_1  \max\{r_v,1\}\Big\}.
\end{eqnarray}

Now, we are ready to construct the following open sets corresponding to $\mathcal{R}_{m,i,j,\varepsilon,\varepsilon_1} \times \Theta_{m,i,j,\varepsilon,\varepsilon_1,l}$ as
\begin{eqnarray}\label{3Small pieces in velocity space}
\mathcal{O}_{m,i,j,\varepsilon,\varepsilon_1,l}:= \Big[ \bigcup_{x\in \mathcal{X}_{m,i,j,\varepsilon,\varepsilon_1,l}}B_{\mathbb{R}^3}(x;\varepsilon_1) \Big]\times \Theta_{m,i,j,\varepsilon,\varepsilon_1,l},
\end{eqnarray}
where the index set is defined as
\begin{eqnarray}\label{3Small pieces in position space}
&&\mathcal{X}_{m,i,j,\varepsilon,\varepsilon_1,l}:=\Big\{ (x_1,x_2,\eta_m(x_1,x_2))+\tau[\cos\theta \hat{x}_{1,m,i,j,\varepsilon}+\sin\theta \hat{x}_{2,m,i,j,\varepsilon} ] +s n_{m,i,j,\varepsilon} \in \mathbb{R}^3:\nonumber\\
& &\hspace{3cm}(x_1,x_2) \in \mathcal{R}_{m,i,j,\varepsilon,\varepsilon_1},\theta\in (\varepsilon l-\varepsilon_1,\varepsilon
l +\varepsilon_1), ~s\in (-\varepsilon_1, \varepsilon_1) \\
&&\hspace{3cm}\tau \in [0,t_f\big((x_1,x_2,\eta_m(x_1,x_2)),\cos\theta\hat{x}_{1,m,i,j,\varepsilon}+\sin\theta \hat{x}_{2,m,i,j,\varepsilon}  \big)]  \Big\}.\nonumber
\end{eqnarray}
We denote that $\mathcal{O}_{m,i,j,\varepsilon,\varepsilon_1,l}$ is an infinite union of open sets and hence is an open set.
Finally, we define
\begin{eqnarray}\label{3Definition of the neighborhood}
\mathcal{O}_{\varepsilon,\varepsilon_1}:= \bigcup_{m,i,j,l}\mathcal{O}_{m,i,j,\varepsilon,\varepsilon_1,l}~\bigcup ~\left\{\mathbb{R}^3 \times B_{\mathbb{R}^3}(0,\varepsilon_1)\right\} ,
\end{eqnarray}
where $1\leq m\leq M_{\Omega,\delta}=O(\theta^{-2})$, $-N_{\varepsilon}\leq i,j \leq N_{\varepsilon}=O( \theta \varepsilon^{-1})$ and
$0\leq l\leq L_{\varepsilon}=O(\varepsilon^{-1})$. Since $\mathcal{O}_{\varepsilon,\varepsilon_1}$ is a union of open sets, it is also an open set.\\

With the covering set $\mathcal{O}_{\varepsilon,\varepsilon_1}$ on hand, we now prove the properties.\\

 \underline{\it Proof of (\ref{3Neighborhood of the singular set})}.  Suppose there exists $(x,v) \in \mathfrak{S}_{B} $, by the definition of $\mathfrak{S}_B$ in (\ref{1singular set} ), there exists $y=x_b(x,v) \in \partial \Omega$ such that $x=y+t_b(x,v)v$ and $v\cdot n(y)=0$ from (\ref{1Backward exit time}) and (\ref{1Backward exit P}). Then $y\in \mathcal{U}_m$ for some $m$, that is, $y=(y_1,y_2,\eta_m(y_1,y_2))$ and $(y_1,y_2) \in \mathcal{R}_{m,i,j,\varepsilon,\varepsilon_1}$ for some $i,j$.

Firstly, for any $|v|\geq 1 $,  we check that
%%%%%%%%%%%%%%%%%%%%%%%%%%%%%%%%%%%%%%%%%%%%%%%%%%%%%%%%%%%%%%%%%%%%%%%%%%%%%%%%%%%%%%%%%%%%%%%%%%%
% Here we should decompose $r_v <1$ and $r_v>1$. Else, $|n_{m,i,j,\epsilon}\cdot v| $ be bounded by some constant times 8C_{\eta}\epsilon, which
% is not included in $\Theta_{m,i,j,\epsilon,\epsilon_1,l}$.
%%%%%%%%%%%%%%%%%%%%%%%%%%%%%%%%%%%%%%%%%%%%%%%%%%%%%%%%%%%%%%%%%%%%%%%%%%%%%%%%%%%%%%%%%%%%%%%%%%%
\begin{eqnarray}\label{3Estimate of the singular vector}
&& |n_{m,i,j,\varepsilon} \cdot \frac{v}{|v|}|  = \Big|[n_{m,i,j,\varepsilon}-n(y_1,y_2,\eta_m(y_1,y_2))]\cdot v+ n(y_1,y_2,\eta_m(y_1,y_2))\cdot \frac{v}{|v|}| \notag\\
& & = \Big|[n_{m,i,j,\varepsilon}-n(y_1,y_2,\eta_m(y_1,y_2))]\cdot \frac{v}{|v|} \Big|\notag\\
& & \leq \frac{1}{\sqrt{1+|\nabla \eta_m(c_{m,i,j,\varepsilon})|^2} }|\nabla \eta_m(c_{m,i,j,\varepsilon})-\nabla \eta_m(y_1,y_2),0)| \notag\\
& & \hspace{3mm}+ \frac{\Big|\sqrt{1+|\nabla\eta_m(y_1,y_2)|^2}-\sqrt{1+|\nabla \eta_m(c_{m,i,j,\varepsilon})|^2}\Big| }{\sqrt{1+|\nabla \eta_m(c_{m,i,j,\varepsilon})|^2} \sqrt{1+|\nabla \eta_m(y_1,y_2})|^2}|\nabla\eta_m(y_1,y_2),-1|\\
& & \leq  |\nabla \eta_m(c_{m,i,j,\varepsilon})-\nabla \eta_m(y_1,y_2)|+\frac{\Big|\sqrt{1+|\nabla \eta_m(y_1,y_2)|^2}-\sqrt{1+|\nabla \eta_m(c_{m,i,j,\varepsilon})|^2}\Big| }{\sqrt{1+|\nabla \eta_m(c_{m,i,j,\varepsilon})|^2}}\notag\\
&&\leq 2 |\nabla \eta_m(c_{m,i,j,\varepsilon})-\nabla \eta_m(y_1,y_2)|,\nonumber
\end{eqnarray}
where we denoted $\nabla \eta_m(y_1,y_2) =(\partial_{y_1} \eta_m(y_1,y_2), \partial_2 \eta_m(y_1,y_2)) $ and used the fact that
\begin{eqnarray*}
&&\frac{\Big|\sqrt{1+|\nabla \eta_m(y_1,y_2)|^2}-\sqrt{1+|\nabla \eta_m(c_{m,i,j,\varepsilon})|^2}\Big| }{\sqrt{1+|\nabla \eta_m(c_{m,i,j,\varepsilon})|^2}} \\
 &&\leq \frac{|\nabla \eta_m(y_1,y_2)-\nabla \eta_m(c_{m,i,j,\varepsilon})|(|\nabla \eta_m(y_1,y_2)|+|\nabla \eta_m(c_{m,i,j,\varepsilon})) }{\sqrt{1+|\nabla \eta_m(c_{m,i,j,\varepsilon})|^2}(\sqrt{1+|\nabla \eta_m(y_1,y_2)|^2}+\sqrt{1+|\nabla \eta_m(c_{m,i,j,\varepsilon})|^2})}\\
&& \leq |\nabla \eta_m(y_1,y_2)-\nabla \eta_m(c_{m,i,j,\varepsilon})|.
\end{eqnarray*}
Using (\ref{3Error of boundary function}), for $(y_1,y_2) \in \mathcal{R}_{m,i,j,\varepsilon,\varepsilon_1}$, we have
\begin{eqnarray*}
 \Big|n_{m,i,j,\varepsilon}\cdot \frac{v}{|v|}\Big| &\leq &  2 ||\eta_m||_{C^2(\mathcal{R}_{m,i,j,\varepsilon,\varepsilon_1})} |c_{m,i,j,\varepsilon}- (y_1,y_2)| \\
&\leq & 8\varepsilon_1 ||\eta_m||_{C^2(\mathcal{R}_{m,i,j,\varepsilon,\varepsilon_1})}\leq  8\varepsilon_1 ||\eta_m||_{C^2(\mathcal{A}_{m})} \leq 8 C_{\eta} \varepsilon_1.
\end{eqnarray*}
Secondly, we consider the case $|v|\leq 1$. From (\ref{3Error of boundary function}) and the similar estimates of $|v| \geq 1$ case , we have
\begin{eqnarray*}
 \Big|n_{m,i,j,\varepsilon}\cdot v \Big| &\leq & |n(y)\cdot v|+|(n(y_1,y_2)-n_{m,i,j,\varepsilon})\cdot v| \\
 &\leq & 2|\nabla \eta_m(c_{m,i,j,\varepsilon})-\nabla \eta_m(y_1,y_2)| \leq  8 C_{\eta} \varepsilon_1.
\end{eqnarray*}
By (\ref{3Joint set of the angular}), we conclude  that $v\in \bigcup_{l=0}^{L_{\varepsilon}}\Theta_{m,i,j,\varepsilon,l}$.
 Since $(y_1,y_2) \in \mathcal{R}_{m,i,j,\varepsilon,\varepsilon_1} \in \mathcal{A}_m$, the distance $s$ in the direction $n_{m,i,j,\varepsilon}$ is less than the height $\sup_{\mathcal{A}_m} |\eta_m|$. From (\ref{3Error of boundary function}), we know that $ |\eta_m(x_1,x_2)|\leq |\nabla \eta_m| |(x_1,x_2)| \leq \varepsilon_1 $, i.e. $ |s| \leq \varepsilon_1$, and hence $(x,v) \in \mathcal{O}_{\varepsilon,\varepsilon_1}$.\\

\underline{\it Proof of (\ref{3Neighborhood of the closed set})}. It suffices to show that there exists a constant $C_* \gg 1$ such that if $(x,v) \in \overline{\mathcal{O}_{\varepsilon,\varepsilon_1}}$ then $(x,v) \in \mathcal{O}_{\varepsilon,C_*\varepsilon_1}$.
Since in the definition (\ref{3Definition of the neighborhood}) the union on $m,i,j,l$ is finite, we have
\begin{eqnarray*}
 \overline{\mathcal{O}_{\varepsilon,\varepsilon_1}}&=&\bigcup_{m,i,j,l} \overline{\mathcal{O}_{m,i,j,\varepsilon,\varepsilon_1,l}} ~\cup~\left\{\mathbb{R}^3 \times  B_{\mathbb{R}^3}(0;\varepsilon_1)\right\}\\
 & = & \bigcup_{m,i,j,l} \Big[~ \overline{\bigcup_{x \in \mathcal{X}_{m,i,j,\varepsilon,\varepsilon_1,l}}B_{\mathbb{R}^3}(x;\varepsilon_1)}~\times ~ \overline{\Theta_{m,i,j,\varepsilon,\varepsilon_1,l}}\Big] ~\cup~\left\{\mathbb{R}^3 \times B_{\mathbb{R}^3}(0;\varepsilon_1)\right\} .
\end{eqnarray*}
Let $z \in \overline{\bigcup_{x \in \mathcal{X}_{m,i,j,\varepsilon,\varepsilon_1,l}}B_{\mathbb{R}^3}(x;\varepsilon_1)} $. From the definition of closed set, we know that, for given $\varepsilon_1 $, there exists some $y \in \bigcup_{x \in \mathcal{X}_{m,i,j,\varepsilon,\varepsilon_1,l}}B_{\mathbb{R}^3}(x;\varepsilon_1)$ such that $|z-y| \leq \varepsilon_1$. Furthermore, we know that there exists some $x \in \mathcal{X}_{m,i,j,\varepsilon,\varepsilon_1,l}$ satisfies
$|y-x| \leq\varepsilon_1 $. So, we derive that $|z-x| \leq |z-y|+|y-x| \leq 2 \varepsilon_1 \leq C_*\varepsilon_1$ for sufficiently large $C_*\gg 1$. That is,
\begin{eqnarray}\label{3Claim for Closed Set}
\overline{\bigcup_{x \in \mathcal{X}_{m,i,j,\varepsilon,\varepsilon_1,l}}B_{\mathbb{R}^3}(x;\varepsilon_1)} ~ \subset ~\bigcup_{x \in \mathcal{X}_{m,i,j,\varepsilon,\varepsilon_1,l}} B_{\mathbb{R}^3}(x; C_*\varepsilon_1).
\end{eqnarray}
\\
On the other hand, from (\ref{3Decomposition of velocity space}), $C_* \gg 1$ and the fact that the vectors $\hat{x}_{1,m,i,j,\varepsilon},~\hat{x}_{2,m,i,j,\varepsilon}$, and $n_{m,i,j,\varepsilon}$ are fixed for given $m,i,j$,
\begin{eqnarray}\label{3Set for Closed Set 2}
\overline{\Theta_{m,i,j,\varepsilon,\varepsilon_1,l}} &=& \Big\{v = r_v \cos\theta_v \cos\phi_v \hat{x}_{1,m,i,j,\varepsilon} + r_v\sin\theta_v \cos\phi_v~\hat{x}_{2,m,i,j,\varepsilon}+r_v \sin\phi_v n_{m,i,j,\varepsilon} \in V : \notag\\
 & & \hspace{5pt} |r_v \sin\phi_v| \leq 8C_{\eta}\varepsilon_1 \max\{r_v,1\} ~ |\theta_v -\varepsilon l|\leq \varepsilon_1 ~\text{for}~r_v \geq 0
 \Big\}\notag\\
 & \subset & \Big\{v = r_v \cos\theta_v \cos\phi_v \hat{x}_{1,m,i,j,\varepsilon} + r_v\sin\theta_v \cos\phi_v~\hat{x}_{2,m,i,j,\varepsilon}+r_v \sin\phi_v n_{m,i,j,\varepsilon} \in V : \notag\\
 & & \hspace{5pt} |r_v \sin\phi_v| \leq 8C_{\eta}C_*\varepsilon_1 \max\{r_v,1\}, ~ |\theta_v -\varepsilon l|\leq C_* \varepsilon_1 ~\text{for}~r_v \geq 0 \Big\}\notag\\
 &=& \Theta_{m,i,j,\varepsilon,C_*\varepsilon_1,l}.
\end{eqnarray}
Finally, we conclude (\ref{3Neighborhood of the closed set}) from (\ref{3Claim for Closed Set})-(\ref{3Set for Closed Set 2}),
\begin{eqnarray*}
\overline{\mathcal{O}_{\varepsilon,\varepsilon_1}}&\subset& \bigcup_{m,i,j}\Big[ \bigcup_{x\in\mathcal{X}_{m,i,j,\varepsilon,C_*\varepsilon_1}}B_{\mathbb{R}^3(x;C_* \varepsilon_1)} \times \Theta_{m,i,j,\varepsilon,C_*\varepsilon_1,l}\Big]~\cup~[\mathbb{R}^3 \times B_{\mathbb{R}^3}(0;C_*\varepsilon_1)]\\
& = &\mathcal{O}_{\varepsilon,C_*\varepsilon_1}.
\end{eqnarray*}

\underline{\it Proof of (\ref{3Neighborhood set measure})}.  From the definition of $\mathcal{O}_{\varepsilon,\varepsilon_1}$, we deduce that
\begin{eqnarray}\label{3Main Estimates 1}
 \iint_{\Omega\times V} \mathbf{1}_{\mathcal{O}_{\varepsilon,C_*\varepsilon}} dv dx
& \leq &\sum_{m,i,j,l} \iint_{\Omega \times V}   \mathbf{1}_{\mathcal{O}_{m,i,j,\varepsilon,C_*\varepsilon,l}} dv dx +m_3(\Omega) O(|\varepsilon|^3) \notag\\
 & \leq &  M_{\Omega,\delta}(2N_{\varepsilon})^2L_{\varepsilon} \times \sup_{m,i,j,l} \iint_{\Omega\times V}\mathbf{1}_{\mathcal{O}_{m,i,j,\varepsilon,C_*\varepsilon,l}} dv dx+m_3(\Omega) O(|\varepsilon|^3)\\
 & &\leq  O(\frac{1}{\varepsilon^3}) \times \sup_{m,i,j,l} \iint_{\Omega\times V}\mathbf{1}_{\mathcal{O}_{m,i,j,\varepsilon,C_*\varepsilon,l}} dv dx+m_3(\Omega) O(|\varepsilon|^3).\notag
\end{eqnarray}
On the one hand, there is some constant $C_V$ such that $\max\{r_v,1\} \leq C_V$ because $V$ is a bounded domain. Then, $ |r_v \sin\phi_v| \leq 8 C_V C_{\eta} \varepsilon_1 $ if $ v \in \Theta_{m,i,j,\varepsilon,\varepsilon_1,l}$. So, it holds
\begin{eqnarray*}
\int_{V}\mathbf{1}_{\Theta_{m,i,j,\varepsilon,C_*\varepsilon,l}} dv
\leq  \int_{r_v}\int_{|\theta_v -l\varepsilon| \leq C_*\varepsilon}\int_{|r_v \sin\phi_v| \leq 8 C_V  C_{\eta} C_* \varepsilon } r_v^2 \sin \phi_v  d\phi_v d\theta_v dr_v \lesssim_{\Omega,V} \varepsilon^2.
\end{eqnarray*}
On the other hand, we claim that for $\varepsilon_1\geq \varepsilon$,
\begin{eqnarray}\label{3Measure of sets}
 m_3 \Big(\bigcup_{x\in \mathcal{X}_{m,i,j,\varepsilon,\varepsilon_1,l}}B_{\mathbb{R}^3}(x;\varepsilon_1)\Big)\leq_{\Omega} \varepsilon_1^2.
\end{eqnarray}
Without loss of generality, we assume that $i,j,l=0$. Therefore, $c_{m,i,j,\varepsilon}=0$. For simplicity, we denote the  $\mathbf{e}_1 = \hat{x}_{1,m,0,0,\varepsilon},~\mathbf{e}_2 = \hat{x}_{2,m,0,0,\varepsilon}$, and $\mathbf{e}_3 = n_{m,0,0,\varepsilon}$. Then
\begin{eqnarray}
 \mathcal{X}_{m,i,j,\varepsilon,\varepsilon_1,l}&\subset & \Big\{ (x_1,x_2,\eta_m(x_1,x_2))+\tau[\cos\theta \mathbf{e}_1 +\sin\theta \mathbf{e}_2] +s \mathbf{e}_3 \in \mathbb{R}^3:\nonumber\\
& &~~~~~~~~~~(x_1,x_2) \in \mathcal{R}_{m,0,0,\varepsilon,\varepsilon_1},\theta\in (-\varepsilon_1, \varepsilon_1), \\
&&~~~~~~~~~~ \tau \in [0,t_f\big((x_1,x_2,\eta_m(x_1,x_2)),\cos \theta\mathbf{e}_1+\sin\theta \mathbf{e}_2\big)],~s\in (-\varepsilon_1, \varepsilon_1)  \Big\}.\nonumber
\end{eqnarray}
Let $diam(\Omega)= \sup_{x,y\in \Omega}|x-y|<+\infty$. Since $||\cos \theta\mathbf{e}_1+\sin\theta \mathbf{e}_2||=1 $, the exit time $t_f$ satisfies
\begin{eqnarray*}
t_f\big((x_1,x_2,\eta_m(x_1,x_2)),\cos\theta \mathbf{e}_1 +\sin\theta \mathbf{e}_2\big) \leq diam(\Omega).
\end{eqnarray*}
From the definition of $ \bigcup_{x\in \mathcal{X}_{m,i,j,\varepsilon,\varepsilon_1,l}}B_{\mathbb{R}^3}(x;\varepsilon_1)$, we can check that
it is included in the truncated cone with height $diam(\Omega)$, top radius $[10+||\eta||_{C^1(\mathcal{A}_m)}]\varepsilon_1$
and the bottom radius $[10+||\eta||_{C^1(\mathcal{A}_m)}+ diam(\Omega)
 ||\eta||_{C^2(\mathcal{A}_m)}]\varepsilon_1$. More precisely, it holds that
\begin{eqnarray*}
 \bigcup_{x\in \mathcal{X}_{m,i,j,\varepsilon,\varepsilon_1,l}}B_{\mathbb{R}^3}(x;\varepsilon_1) \subset
 \bigcup_{\tau=0}^{2diam(\Omega)}B_{\mathbb{R}^3}(\tau \mathbf{e}_1;[10+||\eta||_{C^1(\mathcal{A}_m)}+ \tau
 ||\eta||_{C^2(\mathcal{A}_m)}]\varepsilon_1).
\end{eqnarray*}
 Therefore, using (\ref{3Error of boundary function}), we conclude (\ref{3Measure of sets})
 \begin{eqnarray*}
  m_3 \Big(\bigcup_{x\in \mathcal{X}_{m,i,j,\varepsilon, \varepsilon_1, l}}B_{\mathbb{R}^3}(x;\varepsilon_1)\Big)
  &\leq& m_3\Big(\bigcup_{\tau=0}^{2diam(\Omega)}B_{\mathbb{R}^3}(\tau \mathbf{e}_1;[10+||\eta||_{C^1(\mathcal{A}_m)}+ \tau
 ||\eta||_{C^2(\mathcal{A}_m)}]\varepsilon_1) \Big) \\
 &\leq & 3 diam(\Omega) \Big[10+||\eta||_{C^1(\mathcal{A}_m)}+ \tau
 ||\eta||_{C^2(\mathcal{A}_m)} \Big]^2 \times (\varepsilon_1)^2 \\
 &\leq & 3 diam(\Omega) \Big[10+\frac{1}{8}+ C_{\eta}diam(\Omega)\Big]^2 (\varepsilon_1)^2\\
 &\leq_{\Omega} & \varepsilon_1^2.
 \end{eqnarray*}
Finally, by selecting $\varepsilon_1=C_* \varepsilon $ in (\ref{3Measure of sets}),  we conclude (\ref{3Neighborhood set measure}) as
\begin{eqnarray*}
\iint_{\Omega \times V} \mathbf{1}_{\mathcal{O}_{\varepsilon,C_*\varepsilon}} dvdx \lesssim O(\frac{1}{\varepsilon^3})  m_3 \Big(\bigcup_{x\in \mathcal{X}_{m,i,j,\varepsilon,C_*\varepsilon,l}}B_{\mathbb{R}^3}(x;\varepsilon_1)\Big)
  \int_{V}\mathbf{1}_{\Theta_{m,i,j,\varepsilon,C_*\varepsilon,l}} dv+m_3(\Omega) O(|\varepsilon|^3)
   \lesssim   \varepsilon.
\end{eqnarray*}

\underline{\it Proof of (\ref{3Neighborhood distance})}. Since (\ref{3Neighborhood of the singular set}) holds for all $ \varepsilon \leq \varepsilon_1 $, it suffices to show there exists $C_2=C_2(C_*)>0$ such that
\begin{eqnarray}\label{3Main result of Distence}
dist(\overline{\Omega} \times \mathbb{R}^3 \ \mathcal{O}_{\varepsilon,C_*\varepsilon}, ~~\overline{\mathcal{O}_{\varepsilon,\varepsilon}})> C_2 \varepsilon.
\end{eqnarray}

By the definition of $\mathcal{O}_{\varepsilon,\varepsilon}$ in (\ref{3Definition of the neighborhood}),
\begin{eqnarray}
 && \text{dist}(\overline{\Omega} \times V \setminus \mathcal{O}_{\varepsilon,C_*\varepsilon}, ~~\overline{\mathcal{O}_{\varepsilon,\varepsilon}}) \nonumber\\
% &&= \inf\{|(x,v)-(y,u)|: (x,v)\in (\mathcal{O}_{\varepsilon,C_*\varepsilon})^c, ~(y,u)\in % \overline{\mathcal{O}_{\varepsilon,\varepsilon}}\}\nonumber\\
 %&&= \inf_{m,i,j,l} \{|(x,v)-(y,u)|: (x,v)\in (\mathcal{O}_{\varepsilon,C_*\varepsilon})^c, ~(y,u)\in %\overline{\mathcal{O}_{m,i,j,\varepsilon,\varepsilon,l}} \cup [\mathbb{R}^3 \times \overline{B_{\mathbb{R}^3}(0,\varepsilon)}]\}\nonumber\\
 &&\geq \inf_{m,i,j,l} \Big\{|(x,v)-(y,u)|: (x,v)\in (\mathcal{O}_{m,i,j,\varepsilon,\varepsilon,l})^c \cap [\mathbb{R}^3 \times (B_{\mathbb{R}^3}(0,C_*\varepsilon)]^c , \nonumber\\
 && \qquad\qquad\qquad~(y,u)\in \overline{\mathcal{O}_{m,i,j,\varepsilon,\varepsilon,l}} \cup [\mathbb{R}^3 \times B_{\mathbb{R}^3}(0,\varepsilon)]\Big\}\nonumber\\
  & & = \min \Big\{\inf_{m,i,j,l}\{ |(x,v)-(y,u)|: (x,v)\in (\mathcal{O}_{m,i,j,\varepsilon,\varepsilon,l})^c \cap [\mathbb{R}^3 \times (B_{\mathbb{R}^3}(0,C_*\varepsilon))^c], \nonumber\\
 && \qquad\qquad\qquad~(y,u)\in \mathbb{R}^3 \times B_{\mathbb{R}^3}(0;\varepsilon)\}, \label{3Dist 1}\\
 && \hspace{1cm} \inf_{m,i,j,l}\{ |(x,v)-(y,u)|: (x,v)\in (\mathcal{O}_{m,i,j,\varepsilon,\varepsilon,l})^c \cap [\mathbb{R}^3 \times (B_{\mathbb{R}^3}(0,C_*\varepsilon))^c], \nonumber\\
 && \qquad\qquad\qquad~(y,u)\in \overline{\mathcal{O}_{m,i,j,\varepsilon,\varepsilon,l}} \cap [\mathbb{R}^3 \times B_{\mathbb{R}^3}(0;\varepsilon)]^c \}\Big\}.\label{3Dist 2}
  \end{eqnarray}
Clearly, we have
\begin{eqnarray*}
(\ref{3Dist 1}) &\geq& \inf\{ |(x,v)-(y,u)|: (x,v)\in (\mathbb{R}^3 \times (B_{\mathbb{R}^3}(0,C_* \varepsilon))^c,
 ~(y,u)\in \mathbb{R}^3 \times B_{\mathbb{R}^3}(0;\varepsilon)\}\\
 &\geq & \inf\{|v-u|: v\in B_{\mathbb{R}^3}(0;C_*\varepsilon)^c,~u\in B_{\mathbb{R}^3}(0;\varepsilon)  \} \\
 &=& (C_*-1)\varepsilon.
\end{eqnarray*}

Now, we consider the lower bound of (\ref{3Dist 2}). Firstly, from the definition of $\mathcal{O}_{m,i,j,\varepsilon,C_*\varepsilon,l}$ in (\ref{3Small pieces in position space}), we divide $\{(x,v)\in (\mathcal{O}_{m,i,j,\varepsilon,C_*\varepsilon,l})^c\}$ in (\ref{3Dist 2}) into two parts,  we deduce that
\begin{eqnarray*}
(\mathcal{O}_{m,i,j,\varepsilon,C_*\varepsilon,l})^c &=& \Big[ \bigcup_{x\in \mathcal{X}_{m,i,j,\varepsilon,C_*\varepsilon,l}}B_{\mathbb{R}^3}(x;C_*\varepsilon)\Big] \times (\Theta_{m,i,j,\varepsilon,C_*\varepsilon,l})^c \\
&& ~\bigcup ~ \Big[ \bigcap_{x\in \mathcal{X}_{m,i,j,\varepsilon,C_*\varepsilon,l}}\big(B_{\mathbb{R}^3}(x;C_*\varepsilon)\big)^c\Big] \times V.
 \end{eqnarray*}
Therefore, (\ref{3Dist 2}) is bounded below by the minimum of the following two numbers
\begin{eqnarray*}
 && \inf\{|(x,v)-(y,u)|: (x,v) \in \Big[\bigcup_{x\in \mathcal{X}_{m,i,j,\varepsilon,C_*\varepsilon,l}}B_{\mathbb{R}^3}(x,C_*\varepsilon)\Big]\times [(\Theta_{m,i,j,\varepsilon,C_*\varepsilon,l})^c \setminus  B_{\mathbb{R}^3}(0;C_*\varepsilon) ],\nonumber\\
&& \hspace{3cm}(y,u) \in \overline{\mathcal{O}_{m,i,j,\varepsilon,C_*\varepsilon,l}} \cap [\mathbb{R}^3\times B_{\mathbb{R}^3}(0;\varepsilon)^c  ] \},\label{Dist 5}\vspace{3pt}
\\
 && \inf\{|(x,v)-(y,u)|: (x,v) \in \Big[\bigcap_{x\in \mathcal{X}_{m,i,j,\varepsilon,C_*\varepsilon,l}}(B_{\mathbb{R}^3}(x,C_*\varepsilon))^c\Big]\times [\mathbb{R}^3 \setminus  B_{\mathbb{R}^3}(0; C_*\varepsilon) ],\nonumber\vspace{3pt}\\
&& \hspace{3cm}(y,u) \in \overline{\mathcal{O}_{m,i,j,\varepsilon,C_*\varepsilon,l}} \cap [\mathbb{R}^3\times B_{\mathbb{R}^3}(0;\varepsilon)^c  ] \}.\label{Dist 6}
\end{eqnarray*}

Secondly, we consider $\{(y,u) \in \overline{\mathcal{O}_{m,i,j,\varepsilon,C_*\varepsilon,l}}\}$. From (\ref{3Claim for Closed Set}) with $\varepsilon_1=\varepsilon$ for some $\varsigma = \varsigma(\varepsilon,C_*)>0$ such that
\begin{eqnarray*}
\overline{\mathcal{O}_{m,i,j,\varepsilon,\varepsilon,l}} &=& \overline{\bigcup_{x\in\mathcal{X}_{m,i,j,\varepsilon,\varepsilon,l}}B_{\mathbb{R}^3}(x;\varepsilon)}\times \overline{\Theta_{m,i,j,\varepsilon,\varepsilon,l}}\\
&\subset& \Big[ \bigcup_{x\in \mathcal{X}_{m,i,j,\varepsilon,\frac{C_* \varepsilon}{2},l}} B_{\mathbb{R}^3}(x;\frac{C_*\varepsilon}{2})  \Big]\times \overline{\Theta_{m,i,j,\varepsilon,\varepsilon,l}}.
\end{eqnarray*}
So, we conclude that (\ref{3Dist 2}) is bounded below by the minimum of (A) and (B):
\begin{eqnarray*}
(A)&=& \inf\{|(x,v)-(y,u)|: (x,v) \in \bigcup_{x\in \mathcal{X}_{m,i,j,\varepsilon,C_*\varepsilon,l}}B_{\mathbb{R}^3}(x,C_*\varepsilon)\times [(\Theta_{m,i,j,\varepsilon,C_*\varepsilon,l})^c\setminus  B_{\mathbb{R}^3}(0;C_*\varepsilon) ],\nonumber\\
&& \hspace{2.5cm}(y,u) \in \Big[\bigcup_{x\in \mathcal{X}_{m,i,j,\varepsilon,\frac{C_*\varepsilon}{2},l}}B_{\mathbb{R}^3}(x,\frac{C_*\varepsilon}{2})\Big]\times [\overline{\Theta_{m,i,j,\varepsilon,C_*\varepsilon,l}} \setminus B_{\mathbb{R}^3}(0;\varepsilon)  ] \},\label{3Dist 3}
\\
(B)&=& \inf\{|(x,v)-(y,u)|: (x,v) \in \Big[\bigcap_{x\in \mathcal{X}_{m,i,j,\varepsilon,C_*\varepsilon,l}}(B_{\mathbb{R}^3}(x,C_*\varepsilon))^c\Big]\times [\mathbb{R}^3\setminus   B_{\mathbb{R}^3}(0; C_*\varepsilon) ],\nonumber\\
&& \hspace{2.5cm}(y,u) \in \Big[\bigcup_{x\in \mathcal{X}_{m,i,j,\varepsilon,\frac{C_*\varepsilon}{2},l}}B_{\mathbb{R}^3}(x,\frac{C_*\varepsilon}{2})\Big]\times [\overline{\Theta_{m,i,j,\varepsilon,C_*\varepsilon,l}} \setminus  B_{\mathbb{R}^3}(0;\varepsilon)  ] \}.\label{3Dist 4}
\end{eqnarray*}
In the following, we firstly prove that $(A)\geq \varepsilon$. Let $v\in (\Theta_{m,i,j,\varepsilon,C_*\varepsilon,l})^c \setminus B_{\mathbb{R}^3}(0; C_*\varepsilon)$. By (\ref{3Decomposition of velocity space}), it could be rewritten as
\begin{eqnarray*}
v = r_v \cos\theta_v \cos\phi_v \hat{x}_{1,m,i,j,\varepsilon} +  r_v \sin\theta_v \cos\phi_v \hat{x}_{2,m,i,j,\varepsilon}+
 r_v  \sin\phi_v n_{m,i,j,\varepsilon},
\end{eqnarray*}
where
\begin{eqnarray}\label{3Dist 7}
 \min\{r_v,1\}|\sin\phi_v| \geq 8C_{\eta}C_* \varepsilon, ~~\text{or}~~ |\theta_v-l\varepsilon| \geq C_*\varepsilon.
\end{eqnarray}
Let $u \in \overline{\Theta_{m,i,j,\varepsilon,\varepsilon,l}}\setminus B_{\mathbb{R}^3 }(0;\varepsilon) $. Again from (\ref{3Decomposition of velocity space}), we have
\begin{eqnarray*}
u = r_u \cos\theta_u \cos\phi_u \hat{x}_{1,m,i,j,\varepsilon} +  r_u \sin\theta_u \cos\phi_u \hat{x}_{2,m,i,j,\varepsilon}+
 r_u  \sin\phi_u n_{m,i,j,\varepsilon},
\end{eqnarray*}
where
\begin{eqnarray}\label{3Dist 8}
 \min\{1,r_v\}| \sin\phi_u| \leq 8C_{\eta} \varepsilon, ~~\text{and}~~ |\theta_u -l\varepsilon| \leq \varepsilon.
\end{eqnarray}

We discuss $(A)\geq \varepsilon$ in the following cases.

1) If $|\theta_v -l\varepsilon| \geq C_* \varepsilon$ for $C_*\gg 1$, then clearly $|v-u| \gtrsim \varepsilon$ since $|\theta_u -l\varepsilon| \leq \varepsilon$.

2) For the case $|\theta_v -l\varepsilon| \leq C_*\varepsilon$, it would be divided into several cases.

 (a) If $|r_v|,~|r_u| \leq 1$, then $|r_v \sin\phi_v| \geq 8 C_{\eta} C_* \varepsilon$ from (\ref{3Dist 7}) and $|r_u \sin\phi_u| \leq 8 C_{\eta} \varepsilon$ from (\ref{3Dist 8}). So,
\begin{eqnarray*}
|v-u| &\geq & |(v-u) \cdot n_{m,i,j,\varepsilon}| \geq |v\cdot n_{m,i,j,\varepsilon}|-|u\cdot n_{m,i,j,\varepsilon}| \\
&\geq &  |r_v \sin\phi_v|-|r_u\sin\phi_u| \geq 8C_{\eta}C_* \varepsilon -8C_{\eta} \varepsilon\\
& \gtrsim & \varepsilon.
\end{eqnarray*}

(b) If $|r_v| \geq 1$ and $|r_u|\leq 1$, then
 $|\sin\phi_v| \geq 8 C_{\eta} C_* \varepsilon$ from (\ref{3Dist 7}) and $|r_u \sin\phi_u| \leq 8 C_{\eta} \varepsilon$ from (\ref{3Dist 8}). So
\begin{eqnarray*}
|v-u| &\geq & |(v-u) \cdot n_{m,i,j,\varepsilon}| \geq   |r_v \sin\phi_v|-|r_u\sin\phi_u| \\
&\geq & |\sin\phi_v |-|r_u\sin\phi_u| \geq  8C_{\eta}C_* \varepsilon -8C_{\eta} \varepsilon\\
& \gtrsim & \varepsilon.
\end{eqnarray*}

(c) If $|r_v| \leq 1$ and $|r_u| \geq 1$. Then $|r_v \sin\phi_v| \geq 8 C_{\eta} C_* \varepsilon$ from (\ref{3Dist 7}) and $| \sin\phi_u| \leq 8 C_{\eta} \varepsilon$ from (\ref{3Dist 8}). We will discuss it in the following subcases.

 Fix $0<c_* \ll 1\ll C_*$. When $|r_u|\leq C_*-c_*$, then
\begin{eqnarray*}
|v-u| &\geq & |(v-u) \cdot n_{m,i,j,\varepsilon}| \geq |v\cdot n_{m,i,j,\varepsilon}|-|u\cdot n_{m,i,j,\varepsilon}| \\
&\geq &  |r_v \sin\phi_v|-|r_u\sin\phi_u| \geq 8C_{\eta}C_* \varepsilon -|r_u| \times 8C_{\eta} \varepsilon\\
& = & 8C_{\eta}(C_*-|r_u|)\varepsilon \geq 8C_{\eta}c_* \varepsilon.
\end{eqnarray*}
When $|r_u| \geq C_*-c_*$, then
\begin{eqnarray*}
|v-u| &\geq & |[u-(u \cdot n_{m,i,j,\varepsilon})n_{m,i,j,\varepsilon}]-[v-(v \cdot n_{m,i,j,\varepsilon})n_{m,i,j,\varepsilon}]|\\ &\geq & |r_u \cos \phi_u|-|r_v|\,|\cos\phi_v| \geq |r_u|\sqrt{1-64(C_{\eta})^2\varepsilon^2}-|\cos\phi_v| \\
& \geq & (C_*-c_*) \sqrt{1-64(C_{\eta})^2\varepsilon^2}-1 \\
&\gtrsim & 1 .
\end{eqnarray*}

(d) If $|r_v|\geq 1$ and $|r_u|\geq 1$. In this case, $|\sin\phi_v|\geq 8C_{\eta}C_* \varepsilon$ from (\ref{3Dist 7}) and
$|\sin\phi_u|\leq 8C_{\eta}\varepsilon$  from (\ref{3Dist 8}). Let $|r_u|=k|r_v|$. We also introduce $0<c_*\ll 1 \ll C_*$.
When $k \leq C_*-c_*$, then
\begin{eqnarray*}
|v-u| \geq  |r_v | |\sin\phi_v -k\sin\phi_u| \geq |\sin\phi_v -k\sin\phi_u| \geq 8C_{\eta} c_* \varepsilon.
\end{eqnarray*}
When $k \geq C_*-c_*$, one has
\begin{eqnarray*}
 |v-u| &\geq & |[u-(u \cdot n_{m,i,j,\varepsilon})n_{m,i,j,\varepsilon}]-[v-(v \cdot n_{m,i,j,\varepsilon})n_{m,i,j,\varepsilon}]|\\ &\geq & |r_u \cos \phi_u|-|r_v||\cos\phi_v| \geq k \sqrt{1-64(C_{\eta})^2\varepsilon^2}-|\cos\phi_v| \\
& \geq & (C_*-c_*) \sqrt{1-64(C_{\eta})^2\varepsilon^2}-1 \\
&\gtrsim & 1.
\end{eqnarray*}
Combing all the cases, we deduce $(A) \gtrsim \varepsilon$ for $\varepsilon$ small enough.\\

Secondly, we prove $(B) \gtrsim \varepsilon$. It is true due to
\begin{eqnarray*}
(B) &\geq& \inf\Big\{|x-y|: x \in \bigcap_{z\in \mathcal{X}_{m,i,j,\varepsilon,C_*\varepsilon,l}}(B_{\mathbb{R}^3}(z,C_*\varepsilon))^c,~y \in \bigcup_{z\in \mathcal{X}_{m,i,j,\varepsilon,\frac{C_*\varepsilon}{2},l}}B_{\mathbb{R}^3}(z,\frac{C_*\varepsilon}{2})\Big\}\\
&\geq & \inf\Big\{|x-y|: x \in \bigcap_{z\in \mathcal{X}_{m,i,j,\varepsilon,\frac{C_*\varepsilon}{2},l}}(B_{\mathbb{R}^3}(z,C_*\varepsilon))^c ,~y \in \bigcup_{z\in \mathcal{X}_{m,i,j,\varepsilon,\frac{C_*\varepsilon}{2},l}}B_{\mathbb{R}^3}(z,\frac{C_*\varepsilon}{2})\Big\}\\
&\geq & \inf_{z\in \mathcal{X}_{m,i,j,\varepsilon,C_*\varepsilon,l}} \Big\{|x-y|: x \in \bigcap(B_{\mathbb{R}^3}(z,C_*\varepsilon))^c,~y \in B_{\mathbb{R}^3}(z,\frac{C_*\varepsilon}{2})\Big\}\\
&\geq & \frac{C_* \varepsilon}{2}.
\end{eqnarray*}
So, the estimate of (\ref{3Neighborhood distance}) from $ (A) \gtrsim \varepsilon,~~\text{and}~~ (B) \gtrsim \varepsilon.
$ immediately. The proof of Lemma \ref{3Neighborhood of the singular set properties} is completed. $\hfill \square$

\subsection{Construction of cut-off function}

Recall the standard mollifier $\psi: \mathbb{R}^3\times \mathbb{R}^3 \rightarrow [0,\infty)$,
\begin{eqnarray*}
 \psi(x,v)=\left\{\begin{array}{lll} C\exp\left(\frac{1}{|x|^2+|v|^2-1}\right),& \text{for}~~|x|^2+|v|^2 < 1,\\
 0,&  \text{for}~~|x|^2+|v|^2 \geq 1 \end{array}\right.
\end{eqnarray*}
where the constant $C>0$ is selected so that $ \iint_{\mathbb{R}^3\times \mathbb{R}^3} \psi(x,v)dxdv=1$.

For each $\varepsilon >0$, set
\begin{eqnarray}\label{3Definition of Mollifier1}
\psi_{\varepsilon}(x,v) = (\frac{\varepsilon}{\tilde{C}})^6 \psi(\frac{\sqrt{|x|^2+|v|^2}}{\varepsilon/\tilde{C}}),
~~\text{with}~~\tilde{C} \gg C_* \gg 1.
\end{eqnarray}
Clearly, $\psi_{\varepsilon}$ is smooth and bounded and satisfies
\begin{eqnarray*}
\iint_{\mathbb{R}^3\times \mathbb{R}^3} \psi_{\varepsilon}(x,v)dxdv=1,~~~~\text{spt}(\psi_{\varepsilon})\subset B_{\mathbb{R}^3\times \mathbb{R}^3}(0;\varepsilon/\tilde{C}).
\end{eqnarray*}

\begin{defi} \label{3Definition of Mollifier 2}
We define a smooth cut-off function $\chi_{\varepsilon}: \overline{\Omega} \times V \rightarrow [0,2]$ as
\begin{eqnarray}\label{3Definition of Mollifier3}
\begin{array}{lll} \chi_{\varepsilon}(x,v)&:=& \mathbf{1}_{\overline{\Omega} \times V \setminus \mathcal{O}_{\varepsilon,C_*\varepsilon}} * \psi_{\varepsilon}(x,v) \\& \,= & \displaystyle{\iint_{\mathbb{R}^3\times\mathbb{R}^3} \mathbf{1}_{\overline{\Omega} \times V \setminus \mathcal{O}_{\varepsilon,C_*\varepsilon}}(y,u)\psi_{\varepsilon}(x-y,v-u)dydu }.
\end{array}\end{eqnarray}
\end{defi}
The following properties of the cut-off function are crucial for our analysis.

\begin{pro} \label{3Modified cover function} There exist $\tilde{C}\gg C_*\gg 1$ in (\ref{3Definition of Mollifier1}) and
(\ref{3Definition of Mollifier3}) and $\varepsilon_0=\varepsilon_0(\Omega,V)>0$ such that if $0\leq \varepsilon\leq \varepsilon_0$, then
\begin{eqnarray}\label{3Set Included Singular Poits}
\mathfrak{S}_B \subset \{ (x,v) \in \overline{\Omega} \times V: \chi_{\varepsilon}(x,v)=0\},
\end{eqnarray}
and, for either $\partial = \nabla_x $ or $\partial=\nabla_v$, it holds that
\begin{eqnarray}
\iint_{\Omega\times V}[1-\chi_{\varepsilon}(x,v)]dvdx &\lesssim & \varepsilon,\label{3Measure for Singualr set by cover}\\
 \iint_{\Omega\times V}|\partial \chi_{\varepsilon}(x,v)| dvdx &\lesssim & 1\label{3Measure for Singualr Differetnial}.
\end{eqnarray}
\end{pro}

\prof This Proposition will be proved by the definition of the cut-off function $\chi_{\varepsilon}$ directly.

\underline{\it Proof of (\ref{3Set Included Singular Poits})}. Let $(x,v) \in \mathfrak{S}_B$. Due to (\ref{3Definition of Mollifier1}) if $|(x,v)-(y,u)|\geq \varepsilon/\tilde{C}$ then $\psi_{\varepsilon}(x,v)=0$. Therefore, (\ref{3Definition of Mollifier3}) can be rewritten as
\begin{eqnarray*}
\chi_{\varepsilon}(x,v)=  \iint_{B_{\mathbb{R}^6}((x,v);\varepsilon/\tilde{C})} \mathbf{1}_{\overline{\Omega} \times V \setminus \mathcal{O}_{\varepsilon,C_*\varepsilon}}(y,u)\psi_{\varepsilon}(x-y,v-u)dydu
\end{eqnarray*}
On the other hand, if $(y,u) \in B_{\mathbb{R}^6}((x,v);\varepsilon/\tilde{C})$, duo to (\ref{3Neighborhood of the singular set})-(\ref{3Neighborhood of the closed set}) with $\varepsilon_1 =\varepsilon$ and $\tilde{C} \gg C_*\gg 1$, we have $(y,u) \in \overline{\mathcal{O}_{\varepsilon,\varepsilon}} \subset
\mathcal{O}_{\varepsilon,C_*\varepsilon}$ and
\begin{eqnarray*}
 \mathbf{1}_{\overline{\Omega} \times V \setminus \mathcal{O}_{\varepsilon,C_*\varepsilon} }(y,u)=0.
\end{eqnarray*}
Therefore, we conclude that $\chi_{\varepsilon}(x,v)=0$ for all $\mathfrak{S}_B$ and (\ref{3Set Included Singular Poits}) is true.

\underline{\it Proof of (\ref{3Measure for Singualr set by cover})}. We use (\ref{3Neighborhood set measure}) with $\varepsilon_1=\varepsilon$ to have
\begin{eqnarray*}
&&\iint_{\Omega\times V}\iint_{\mathbb{R}^3\times \mathbb{R}^3}[1-\mathbf{1}_{\overline{\Omega} \times V \setminus \mathcal{O}_{\varepsilon,C_*\varepsilon}}(y,u)]\psi_{\varepsilon}(x-y,v-u)dudy dvdx\\
&&\leq  \iint_{\Omega\times V}\mathbf{1}_{\overline{\Omega} \times V\setminus \mathcal{O}_{\varepsilon,C_*\varepsilon}}(y,u)dudy\iint_{\mathbb{R}^3\times \mathbb{R}^3} \psi_{\varepsilon}(x-y,v-u)dvdx\\
&& \leq \frac{C_1 \varepsilon}{2} \iint_{B_{\mathbb{R}^3}(0,\varepsilon/\tilde{C})}\psi_{\varepsilon}(x,v) dvdx\\
&&\leq \varepsilon
\end{eqnarray*}

\underline{\it Proof of (\ref{3Measure for Singualr Differetnial})}. Note that from a standard scaling argument and (\ref{3Definition of Mollifier1}), one has
\begin{eqnarray*}
|\partial \psi_{\varepsilon}(x,v) | \leq \frac{\tilde{C}^6}{\varepsilon^7} \mathbf{1}_{B_{\mathbb{R}^6}(0,\varepsilon/\tilde{C})}(x,v).
\end{eqnarray*}
We also note that $\partial \chi_{\varepsilon} =- \partial[1-\chi_{\varepsilon}]$. Therefore, by Lemma 1,
\begin{eqnarray*}
&&\iint_{\Omega\times V} |\partial \chi_{\varepsilon}(x,v)|dvdx \\
&&=\iint_{\Omega\times V}\iint_{\mathbb{R}^3\times \mathbb{R}^3}[1-\mathbf{1}_{\overline{\Omega} \times \mathbb{R}^3 \setminus \mathcal{O}_{\varepsilon,C_*\varepsilon}}(y,u)]|\partial \psi_{\varepsilon}(x-y,v-u)|dudy dvdx\\
&&\leq  \iint_{\Omega\times V}\mathbf{1}_{\overline{\Omega} \times \mathbb{R}^3 \setminus \mathcal{O}_{\varepsilon,C_*\varepsilon}}(y,u)dudy\iint_{\mathbb{R}^3\times \mathbb{R}^3} O(\varepsilon^{-7} \tilde{C}^6) \mathbf{1}_{B_{\mathbb{R}^3}(0,\varepsilon/\tilde{C})}dvdx\\
&& \leq O(\varepsilon) \times O(\varepsilon^{-1})\\
&&\leq 1.
\end{eqnarray*}
This completed the proof of Proposition \ref{3Modified cover function}. $\hfill\square$

\begin{pro} \label{3Measure for singualr set on boundary} With the same constant $\tilde{C} \gg C_* \gg 1$ as in Proposition \ref{3Modified cover function} and $0<\varepsilon \leq \varepsilon_0$, then
\begin{eqnarray}
 \mathfrak{S}_B \cap [\partial \Omega \times V] \subset \{(x,v)\in \partial \Omega \times V : \chi_{\varepsilon}(x,v)=0\}. \label{3Set Included Boundary Singular Point}
\end{eqnarray}
Moreover if $|y,u| \leq \varepsilon/\tilde{C}$ for $\tilde{C} \gg C_* \gg 1$,
\begin{eqnarray}
\int_{\partial \Omega}\int_{n(x)\cdot v <0}\mathbf{1}_{\varepsilon,C_*\varepsilon}(x-y,v-u)  |n(x-y)\cdot(v-u)|du dS_{x} \lesssim \varepsilon, \label{3Measure on Boundary}
\end{eqnarray}
and
\begin{eqnarray}
\int_{\gamma_-}[1-\chi_{\varepsilon}(x,v)] d\gamma &\lesssim_{\Omega,V} & \varepsilon, \label{3Measure on Outing Boundary}\\
\int_{\gamma_-}|\partial\chi_{\varepsilon}(x,v)|  d\gamma &\lesssim_{\Omega,V} & 1. \label{3Measure on Outing Boundary Partial}
\end{eqnarray}
\end{pro}

The following fact is crucial to prove Proposition \ref{3Measure for singualr set on boundary} and especially (\ref{3Measure on Boundary}). The proof is similar to that in \cite{[GKTT2]}.

\begin{lem}\label{3Classification of velocity space} We fix $m_0=1,2,\cdots, M_{\Omega,\theta} $ in (\ref{3Composition of Boundary}). We may assume (up to rotations and translations) there exists a $C^2$ function $\eta_{m_0}: [-\theta,\theta]\times [-\theta,\theta] \rightarrow \mathbb{R}$, whose graph is the boundary $\mathcal{U} \cap \partial \Omega$. Let $(x_1,x_2) \in \mathcal{A}_{m_0} \cap [-\theta,\theta]\times [-\theta,\theta] $ and $(x_1,x_2)\in \mathcal{R}_{m_0,i_0,j_0,\varepsilon,C_*\varepsilon}$ for $|i_0|, ~|j_0| \leq N_{\varepsilon}$. Furthermore, we suppose that
\begin{itemize}
\item $|y| \leq \frac{\varepsilon}{\tilde{C}}$ and
\begin{eqnarray}\label{3Almost Boundary}
(x_1,x_2, \eta_{m_0}(x_1,x_2)-y,v) \in \mathcal{O}_{\varepsilon,C_*\varepsilon},
\end{eqnarray}
\item For large but fixed $s_* \gg 1$,
\begin{eqnarray}\label{3Almost parallel set}
-1\leq n_{m_0}(0,0) \cdot \frac{v}{|v|} \leq -s_* C_2 \sqrt{\varepsilon},~~\text{with}~~C_2:= \sqrt{\frac{8C_*}{3}[1+||\eta_{m_0}||_{C^2(\mathcal{A}_{m_0})}]}.
\end{eqnarray}
\end{itemize}
Then either $|v| \leq \varepsilon^{1/3}$ or there exists $(i,j) \in [-N_1+i_0, N_1+i_0] \times [-N_1+j_0,N_1+j_0]$ with
\begin{eqnarray}\label{3Range of ij}
N_1:= [\frac{8C_3}{\sqrt{\varepsilon}}],~~~C_3 := \frac{4C_*+8C_* [1+||\eta_{m_0}||_{C^2(\mathcal{A}_{m_0})}]^{1/2}+2/\tilde{C} }{s_* C_2},
\end{eqnarray}
such that
\begin{eqnarray*}
 \big((x_1,x_2,\eta_{m_0}(x_1,x_2)-y,v)\big) \in \bigcup_{0\leq l\leq L_{\varepsilon}} \mathcal{O}_{m_0,i,j,\varepsilon,C_*\varepsilon,l} \cap [\overline{\Omega}\times \{v\in \mathbb{R}^3: |v|\geq \varepsilon^{1/3}],
\end{eqnarray*}
and
\begin{eqnarray}
|n_{m_{0}}(0,0)\cdot \frac{v}{|v|}| \leq C_4 \sqrt{\varepsilon} ~\quad~~~\text{with}~~\quad~~C_4=C_3 [1+||\eta_{m_0}||_{C^2(\mathcal{A}_{m_0})}].\label{3Estimate of angle}
\end{eqnarray}
\end{lem}

{\it Proof of Lemma \ref{3Classification of velocity space}}. Without loss of generality (up to rotations and translations), we may assume
\begin{eqnarray}\label{3Simplity for Boundary function}
 (i_0,j_0)=(0,0) \quad \text{and} \quad \eta_{m_0}(0,0)=0 \quad \text{and} \quad \nabla \eta_{m_0}(0,0)=0.
\end{eqnarray}
Consider the case of $|v| \geq \varepsilon^{1/3}$. Since $ \big((x_1,x_2,\eta_{m_0}(x_1,x_2)-y,v)\big)  \in \mathcal{O}_{\varepsilon,C_*\varepsilon}$ we use the definition of $\mathcal{O}_{\varepsilon,C_*\varepsilon} $ in (\ref{3Definition of the neighborhood}) to have
\begin{eqnarray}\label{3Possible choosing}
\text{either} \quad |v| \leq C_*\varepsilon \quad \text{or}\quad (x-y,v) \in \bigcup_{m,i,j,l}\mathcal{O}_{m,i,j,\varepsilon,C_*\varepsilon,l}.
\end{eqnarray}
For small $0<\varepsilon \ll 1$, we can exclude the first case of $|v|\leq C_* \varepsilon$ since $|v|\geq \varepsilon^{1/3} \gg C_* \varepsilon$.

Now we consider the latter case in (\ref{3Possible choosing}). In this case, we claim that
\begin{eqnarray}\label{3Range of m0}
  \big((x_1,x_2,\eta_{m_0}(x_1,x_2)-y,v)\big) \in \bigcup_{i,j,l}\mathcal{O}_{m_0,i,j,\varepsilon,C_*\varepsilon,l}.
\end{eqnarray}
From (\ref{3Possible choosing}) and the definition of $\mathcal{O}_{m,i,j,\varepsilon,C_*\varepsilon,l} $ in (\ref{3Small pieces in velocity space}), there exist $m,i,j,l$ such that
\begin{eqnarray*}
  \big((x_1,x_2,\eta_{m_0}(x_1,x_2)-y,v)\big) \in \Big[\bigcup_{p \in \mathcal{X}_{m,i,j,\varepsilon,C_*\varepsilon,l}}B_{\mathbb{R}^3}(p;C_*,\varepsilon)\Big] \times \Theta_{m,i,j,\varepsilon,C_*\varepsilon,l}.
\end{eqnarray*}
In particular, there exists $p \in \mathcal{X}_{m,i,j,\varepsilon,C_*\varepsilon,l}$ satisfying $ |p- (x_1,x_2,\eta_{m_0}(x_1,x_2)-y )| \leq C_*\varepsilon.$ By the definition of $ \mathcal{X}_{m,i,j,\varepsilon,C_*\varepsilon,l}$ in (\ref{3Small pieces in position space}), one has
\begin{eqnarray*}
p=\big(\overline{p}_1,\overline{p}_2,\eta_{m}(\overline{p}_1,\overline{p}_2) \big)+\overline{\tau}[\cos\overline{\theta} \hat{x}_{1,m,i,j,\varepsilon}+\sin\overline{\theta} \hat{x}_{2,m,i,j,\varepsilon}] +\overline{s} n_{m,i,j,\varepsilon},
\end{eqnarray*}
for some
\begin{eqnarray*}
(\overline{p}_1,\overline{p}_2) &\in & \mathcal{R}_{m,i,j,\varepsilon,C_*\varepsilon},\\
\overline{\theta} &\in & (l\varepsilon-C_*\varepsilon,l\varepsilon+C_*\varepsilon),\\
 \overline{\tau} &\in& \big[0,t_{\mathbf{f}} \big( \overline{p}_1,\overline{p}_2,\eta_{m}(\overline{p}_1,\overline{p}_2)    ,(\cos\overline{\theta} \hat{x}_{1,m,i,j,\varepsilon} +\sin\overline{\theta} \hat{x}_{2,m,i,j,\varepsilon})  \big)\big],\\
 \overline{s} &\in & [-C_*\varepsilon, C_*\varepsilon].
\end{eqnarray*}
By the definition of $t_{\mathbf{f}}$ in (\ref{1Forward exit time}),
\begin{eqnarray*}
 z:= p-\overline{s} n_{m,i,j,\varepsilon} = \big(\overline{p}_1,\overline{p}_2,\eta_{m}(\overline{p}_1,\overline{p}_2)\big)+
\overline{\tau}[\cos\overline{\theta} \hat{x}_{1,m,i,j,\varepsilon}+\sin\overline{\theta} \hat{x}_{2,m,i,j,\varepsilon}]
 \in \Omega.\end{eqnarray*}
Then, we have
\begin{eqnarray} \label{3Distance of z}
&& |z-\big((x_1,x_2,\eta_{m_0}(x_1,x_2)-y)\big)| \notag\\
&&\leq  |z-p|+|p-\big((x_1,x_2,\eta_{m_0}(x_1,x_2)-y)\big)| \\
&&\leq 2 C_*\varepsilon.
\notag
\end{eqnarray}
From (\ref{3Simplity for Boundary function}) and (\ref{3Distance of z}), and $|y| \leq \varepsilon/\tilde{C}$, we deduce
\begin{eqnarray*}
&& |z-\big(0,0,\eta_{m_0}(0,0)\big)|\\
&&\leq |z-\big((x_1,x_2,\eta_{m_0}(x_1,x_2)-y)\big)| +|\big((x_1,x_2,\eta_{m_0}(x_1,x_2))\big)-\big(0,0,\eta_{m_0}(0,0)\big)|+|y|  \\
&&\leq 2C_* \varepsilon +4C_*\varepsilon (1+||\eta_{m_0}||_{C^1(\mathcal{A}_{m_0})})+\varepsilon/\tilde{C}.
\end{eqnarray*}
Denote $(\overline{z}_1,\overline{z}_2)= (\overline{p}_1,\overline{p}_2)$. By the definition of $t_{\mathbf{b}}$ and $t_{\mathbf{f}}$ in (\ref{1Backward exit time}) and (\ref{1Forward exit time}),
\begin{eqnarray}\label{3Position of Backward}
x_{\mathbf{b}}(z, \cos\overline{\theta} \hat{x}_{1,m,i,j,\varepsilon}+\sin\overline{\theta} \hat{x}_{2,m,i,j,\varepsilon}+0n_{m,i,j,\varepsilon}  )= (\overline{z}_1,\overline{z}_2,\eta_{m_0}(\overline{z}_1,\overline{z}_2)).
\end{eqnarray}
On the other hand, by the definition of $\Theta_{m,i,j,\varepsilon,C_*\varepsilon,l}$ in (\ref{3Decomposition of velocity space}),
\begin{eqnarray}\label{3Expression of v}
\frac{v}{|v|}= \cos \theta_v \cos\phi_v \hat{x}_{1,m,i,j,\varepsilon}+\sin \theta_v \cos\phi_v \hat{x}_{2,m,i,j,\varepsilon}+ \sin\phi_vn_{m,i,j,\varepsilon},
\end{eqnarray}
with
\begin{eqnarray*}
&&|\theta_v -l\varepsilon| \leq C_*\varepsilon,\\
&& |v\cdot n_{m,i,j,\varepsilon}|\leq 8C_{\eta}C_*\varepsilon \,\,\text{for}\,\,\varepsilon^{1/3}\leq |v|\leq 1,\\
&&|\frac{v}{|v|}\cdot n_{m,i,j,\varepsilon}|\leq 8C_{\eta}C_*\varepsilon \,\,\text{for}\,\,\varepsilon^{1/3}\leq |v|\leq 1.
\end{eqnarray*}
Therefore, for $0<\varepsilon\ll 1$,
\begin{eqnarray}\label{3Estimate of Angle1}
 |\frac{v}{|v|}\cdot n_{m,i,j,\varepsilon}|=|\sin\phi_v| \leq \max\{8C_{\eta}C_*\varepsilon^{2/3},8C_{\eta}C_*\varepsilon\} \leq 16 C_{\eta}C_*\varepsilon^{2/3}.
\end{eqnarray}
Now we estimate as
\begin{eqnarray*}
&&n_{m_0}(0,0) \cdot (\cos\overline{\theta} \hat{x}_{1,m,i,j,\varepsilon}+\sin\overline{\theta} \hat{x}_{2,m,i,j,\varepsilon}+0n_{m,i,j,\varepsilon} ) \\
&& \leq n_{m_0}(0,0) \cdot\frac{v}{|v|}+ n_{m_0}(0,0) \cdot \Big(\frac{v}{|v|}- (\cos\overline{\theta} \hat{x}_{1,m,i,j,\varepsilon}+\sin\overline{\theta} \hat{x}_{2,m,i,j,\varepsilon}+0n_{m,i,j,\varepsilon})\Big).
\end{eqnarray*}
We use (\ref{3Expression of v})-(\ref{3Estimate of Angle1}), and $\overline{\theta} \in (l\varepsilon-C_*\varepsilon, l\varepsilon+C_*\varepsilon)$ to conclude that, for $0<\varepsilon\ll 1$,
\begin{eqnarray*}
&&\Big| \frac{v}{|v|}- \cos\overline{\theta} \hat{x}_{1,m,i,j,\varepsilon}+\sin\overline{\theta} \hat{x}_{2,m,i,j,\varepsilon}+0n_{m,i,j,\varepsilon}\Big| \\
&& \leq 2\{|\cos\theta_v-\cos\overline{\theta}|+ |\cos\theta_v||\cos\phi_v-1|+|\sin\theta_v-\sin\overline{\theta}|
+|\sin\theta_v||\cos\theta_v-1|+|\sin\theta_v| \}\\
&&\leq 2\{4C_*\varepsilon+16 C_{\eta}C_* \varepsilon^{2/3}+2(16C_{\eta}C_*)^2\varepsilon^{4/3} \}\\
&& \leq 200C_{\eta}C_* \varepsilon^{2/3}.
\end{eqnarray*}

Finally from (\ref{3Almost parallel set}), for $0<\varepsilon\ll 1$,
\begin{eqnarray}\label{3Estimate of Almost parallel set}
-1&\leq& n_{m_0}(0,0)\cdot \big( \cos\overline{\theta} \hat{x}_{1,m,i,j,\varepsilon}+\sin\overline{\theta} \hat{x}_{2,m,i,j,\varepsilon}+0n_{m,i,j,\varepsilon}\big)\notag\\
& \leq &-s_*\times C_2 \sqrt{\varepsilon} +400C_{\eta}C_*\varepsilon^{2/3}\\
& \leq &-\frac{s_*\times C_2}{2} \sqrt{\varepsilon}.\notag
\end{eqnarray}
Now we are ready to prove the claim (\ref{3Range of m0}). Denote
\begin{eqnarray*}
 \hat{u}:= cos\overline{\theta} \hat{x}_{1,m,i,j,\varepsilon}+\sin\overline{\theta} \hat{x}_{2,m,i,j,\varepsilon}.
\end{eqnarray*}
Recall that $|z| \leq (2C_*+4C_*(1+||\eta_{m_0}||_{C^1(\mathcal{A}_{m_0})})+1/\tilde{C} )\varepsilon $ and $z\in \Omega$.
Therefore for $0<\varepsilon \ll 1$, the function $\eta_{m_0} $ is defined around $(z_1,z_2)$ and $z_3 \geq \eta_{m_0}(z_1,z_2)$.

We define, for $|\tau| \ll 1$,
\begin{eqnarray}\label{3Auxilly Function}
\Phi(\tau) =z_3-\hat{u}_3 \tau -\eta_{m_0}(z-\hat{u}_1\tau,z_2-\hat{u}_2\tau), \quad\Phi(0)>0.
\end{eqnarray}
Expanding $\Phi(\tau)$ in $\tau$, form $-\hat{u}_3 =n_{m_0}(0,0)\cdot \hat{u} $ and
(\ref{3Estimate of Almost parallel set}), we have
\begin{eqnarray*}
\Phi(\tau) &\leq & -\hat{u}_3 +|z_3|+|\eta_{m_0}(z_1-\hat{u}_1\tau,z_2-\hat{u}_2\tau)|\\
&\leq & -\frac{s_*C_2}{2} \sqrt{\varepsilon}\tau + (2C_*+4C_*(1+||\eta_{m_0}||_{C^1(\mathcal{A}_{m_0})})+1/\tilde{C} )\varepsilon\\
& & + ||\eta_{m_0}||_{C^2(\mathcal{A}_{m_0})}(2C_*+4C_*(1+||\eta_{m_0}||_{C^1(\mathcal{A}_{m_0})})+1/\tilde{C} )^2\varepsilon^2 \\
&& +||\eta_{m_0}||_{C^2(\mathcal{A}_{m_0})}|\tau|^2,
\end{eqnarray*}
where we have used the fact that
\begin{eqnarray*}
&&\eta_{m_0}(z_1-\hat{u}_1\tau,z_2-\hat{u}_2\tau)\\
&& = \eta_{m_0}(z_1,z_2) +\int_0^{\tau} \frac{d}{ds} \eta_{m_0}(z_1-\hat{u}_1s,z_2-\hat{u}_2s )ds\\
&&= \eta_{m_0}(z_1,z_2)-(\hat{u}_1,\hat{u}_2)\cdot\nabla \eta_{m_0}\tau  +\int_0^{\tau}\int_0^s \frac{d^2}{ds^2} \eta_{m_0}(z_1-\hat{u}_1s_1,z_2-\hat{u}_2s_1 )ds_1ds\\
&&\leq ||\eta_{m_0}||_{C^2(\mathcal{A}_{m_0})}\frac{|z|^2}{2}+ |(\hat{u}_1,\hat{u}_2)\cdot\nabla \eta_{m_0}(0,0)|\,|\tau|
+ ||\eta_{m_0}||_{C^2(\mathcal{A}_{m_0})}(|z||\tau|+\frac{|\tau|^2}{2}) \\
&& \leq ||\eta_{m_0}||_{C^2(\mathcal{A}_{m_0})}(|z|^2 +|\tau|^2).
\end{eqnarray*}
Now we plug $\tau=\frac{C_3 \sqrt{\varepsilon}}{s_*}$ with the constant $C_3$ in (\ref{3Range of ij}) to have, for $s_* \gg1$ and $0<\varepsilon \ll 1$,
\begin{eqnarray*}
 \Phi(\tau) &\leq& \Big[\frac{C_2C_3}{2}-(2C_*+4C_*[1+||\eta_{m_0}||_{C^1(\mathcal{A}_{m_0})}]+1/\tilde{C})
 -\frac{||\eta_{m_0}||_{C^2(\mathcal{A}_{m_0})}C_3^2}{(s_*)^2}\Big]\varepsilon + O(\varepsilon^2)\\
 &<&0.
\end{eqnarray*}
By the mean value theorem, there exists at least one $\tau \in (0,C_3 \sqrt{\varepsilon}]$ satisfying $\Phi(\tau)=0$.
We choose the smallest one if them and denote it as $\tau_0 \in (0,C_3 \sqrt{\varepsilon}]$. By this definition and (\ref{3Position of Backward}), for $0<\varepsilon\ll 1$,
\begin{eqnarray*}
x_{\mathbf{b}}(z,\hat{u}) &=&x_{\mathbf{b}}(z,\cos\overline{\theta} \hat{x}_{1,m,i,j,\varepsilon}+\sin\overline{\theta} \hat{x}_{2,m,i,j,\varepsilon})\\
&= &z-\tau_0 \hat{u} = (z_1-\tau_0 \hat{u}_1,
z_2-\tau_0 \hat{u}_2,z_3-\tau_0 \hat{u}_3).
\end{eqnarray*}
Therefore, $x_{\mathbf{b}}(z, \hat{u}) \in \partial \Omega \cap \mathcal{U}_{m_0}$ and this proves the claim (\ref{3Range of m0}).
For $0<\varepsilon\ll 1$,
\begin{eqnarray*}
 |(z_1-\tau_0 \hat{u}_1,
z_2-\tau_0 \hat{u}_2)| \leq |z|+\tau_0 |\hat{u}| \leq (2C_*+4C_*[1+||\eta_{m_0}||_{C^1(\mathcal{A}_{m_0})}]+1/\tilde{C})\varepsilon +C_3 \sqrt{\varepsilon} \leq 2C_3 \sqrt{\varepsilon}.
\end{eqnarray*}
 Moreover, for $|i-i_0|, |j-j_0| \leq (2C_3 \sqrt{\varepsilon})/\varepsilon \leq 2C_3 \times \frac{1}{\sqrt{\varepsilon}}\leq N_1$,
 \begin{eqnarray*}
  (z_1-\tau_0 \hat{u}_1,
z_2-\tau_0 \hat{u}_2) \in \mathcal{R}_{m_0,i,j,\varepsilon,C_*\varepsilon}.
 \end{eqnarray*}
Finally, we need to prove (\ref{3Estimate of angle}). From (\ref{3Estimate of Angle1}) and (\ref{3Range of ij})
\begin{eqnarray*}
\big|n_{m_0}(0,0) \cdot \frac{v}{|v|}\big| &\leq & \big|n_{m_0,i,j,\varepsilon,C_*\varepsilon} \cdot \frac{v}{|v|}\big|+
\big|(n_{m_0}(0,0)-n_{m_0,i,j,\varepsilon,C_*\varepsilon}) \cdot \frac{v}{|v|}\big|\\
&\leq & 16 C_{\eta}C_* \varepsilon^{2/3} +|| n_{m_0}||_{C^1(\mathcal{A}_{m_0})}[N_1\varepsilon+C_*\varepsilon]\\
&\leq  & 16 C_{\eta}C_* \varepsilon^{2/3} +|| n_{m_0}||_{C^1(\mathcal{A}_{m_0})}[2C_3\sqrt{\varepsilon}+C_*\varepsilon]\\
&\leq & 10C_3(1+||\eta_{m_0}||_{C^2(\mathcal{A}_{m_0})}) \sqrt{\varepsilon}\\
&\leq & C_4 \sqrt{\varepsilon},
\end{eqnarray*}
and (\ref{3Estimate of angle}) follows. $\hfill\square $\\

Now, we turn to the proof of Proposition \ref{3Classification of velocity space}.

{\it\bf Proof of Proposition \ref{3Classification of velocity space}}. The first statement (\ref{3Set Included Boundary Singular Point}) is clear from (\ref{3Set Included Singular Poits}).

\underline{\it Proof of (\ref{3Measure on Boundary})}. Let $(y,u)\leq \varepsilon/\tilde{C}$. We use (\ref{3Composition of Boundary}) to decompose  that
\begin{eqnarray*}
&&\int_{\partial \Omega}\int_{n(x)\cdot v <0}\mathbf{1}_{\varepsilon,C_*\varepsilon}(x-y,v-u) |n(x-y)\cdot(v-u)|du dS_{x} \\
&& \leq \sum_{m=1}^{M_{\Omega,\theta}} \int_{\mathcal{U}_m \cap \partial \Omega}\int_{n_m(x)\cdot v<0} \mathbf{1}_{\mathcal{O}_{\varepsilon,C_*\varepsilon}}(x-y,v-u) |n_m(x-y)\cdot (v-u)|dS_xdv\\
&& \leq M_{\Omega,\theta}\times \sup_{m} \int_{\mathcal{U}_m \cap \partial \Omega}\int_{n_m(x)\cdot v<0} \mathbf{1}_{\mathcal{O}_{\varepsilon,C_*\varepsilon}}(x-y,v-u)  |n_m(x-y)\cdot (v-u)|dS_xdv.
%\\
%&& \leq \frac{1}{\theta^2}\sup_{m}\int_{\mathcal{U}_m \cap \partial \Omega}\int_{n_m(x)\cdot v<0} %\mathbf{1}_{\mathcal{O}_{\varepsilon,C_*\varepsilon}}(x-y,v-u) |n_m(x-y)\cdot (v-u)|dS_xdv.
\end{eqnarray*}
For fixed $m=1,2,\cdots, M_{\Omega,\theta}$, we use (\ref{3Piece of boundary}) and (\ref{3Cover of variable space}) again to decompose
\begin{eqnarray*}
 &&\int_{\mathcal{U}_m \cap \partial \Omega}\int_{n_m(x)\cdot v<0} \mathbf{1}_{\mathcal{O}_{\varepsilon,C_*\varepsilon}}(x-y,v-u) |n_m(x-y)\cdot (v-u)|dS_xdv\\
 && = \int_{\mathcal{A}_m}\int_{n_m(x)\cdot v<0} \mathbf{1}_{\mathcal{O}_{\varepsilon,C_*\varepsilon}}(x_1-y_1,x_2-y_2,\eta_m(x_1,x_2)-y_3, v-u)\\
 && \hspace{3cm} \times  |n_m(x-y)\cdot (v-u)|\sqrt{1+|\nabla\eta_m(x_1,x_2)|}dx_1dx_2 dv \\
 && =\sum_{-N_{\varepsilon}\leq i,j\leq N_{\varepsilon}} \int_{\mathcal{R}_{m,i,j,\varepsilon,C_*\varepsilon}}\int_{n_m(x)\cdot v<0} \mathbf{1}_{\mathcal{O}_{\varepsilon,C_*\varepsilon}}(x_1-y_1,x_2-y_2,\eta_m(x_1,x_2)-y_3, v-u)\\
 && \hspace{3cm} \times  |n_m(x-y)\cdot (v-u)|\sqrt{1+|\nabla\eta_m(x_1,x_2)|}dx_1dx_2 dv \\
 && \leq  \frac{\theta^2}{\varepsilon^2} \sup_{-N_{\varepsilon}\leq i,j\leq N_{\varepsilon}} \int_{\mathcal{R}_{m,i,j,\varepsilon,C_*\varepsilon}}\int_{n_m(x)\cdot v<0} \mathbf{1}_{\mathcal{O}_{\varepsilon,C_*\varepsilon}}(x_1-y_1,x_2-y_2,\eta_m(x_1,x_2)-y_3, v-u)\\
 && \hspace{3cm} \times  |n_m(x-y)\cdot (v-u)|\sqrt{1+|\nabla\eta_m(x_1,x_2)|}dx_1dx_2 dv \\
 &&  \leq \frac{\theta^2}{\varepsilon^2} \sup_{-N_{\varepsilon}\leq i,j\leq N_{\varepsilon}} \int_{[-C_*\varepsilon,C_*\varepsilon]^2}\int_{n_m(x_1,x_2)\cdot (v+u)<0}  \mathbf{1}_{\mathcal{O}_{\varepsilon,C_*\varepsilon}}(x_1-y_1,x_2-y_2,\eta_m(x_1,x_2)-y_3, v) \notag\\
 && \hspace{3cm} \times |n_m(x-y)\cdot v|\sqrt{1+|\nabla\eta_m(x_1,x_2)|}dx_1dx_2 dv
\end{eqnarray*}
where $n_m(x_1,x_2) =\frac{1}{\sqrt{1+|\nabla \eta_m(x_1,x_2)|^2}}(\partial_1 \eta_m(x_1,x_2),\partial_2 \eta_m(x_1,x_2),-1)$.

We fix $i,j$. Without loss of generality (up to rotations and translations), we may assume
\begin{eqnarray*}
c_{m,i,j,\varepsilon}=(0,0),\, \partial_1 \eta_m(0,0)=\partial_2 \eta_m(0,0)=0,\, n_{m,i,j,\varepsilon}=(0,0,-1).
\end{eqnarray*}
For $(x_1,x_2) \in [-C_*\varepsilon,C_*\varepsilon]^2$, $|(y,u)|\leq \varepsilon/\tilde{C}$ and $n_m(x_1,x_2)\cdot(v+u)<0$, we deduce
\begin{eqnarray}
n_{m,i,j,\varepsilon}\cdot v &=& n_m(0,0)\cdot v =n_m(x_1,x_2)\cdot (v+u) + [n_m(0,0)\cdot v - n_m(x_1,x_2)\cdot (v+u)]\notag\\
& < & |n_m(x_1,x_2)\cdot u|+|n_m(x_1,x_2)\cdot v| \leq  \frac{\varepsilon}{\tilde{C}} +2C_* \varepsilon ||\eta_m||_{C^2 [-C_*\varepsilon,C_*\varepsilon]}|v|\\
& \leq & C_5 \varepsilon, \notag
\end{eqnarray}
where $C_5= \max\{ 1/\tilde{C}, 2C_* ||\eta_m||_{C^2 [-C_*\varepsilon,C_*\varepsilon]} diam{V}\}.$ Therefore, from Lemma \ref{3Classification of velocity space}, we decompose the domain as follows
\begin{eqnarray*}
 && \int_{[-C_*\varepsilon,C_*\varepsilon]^2 } \int_{n_m(0,0)\cdot v \leq C_5(1+|v|)\varepsilon } \mathbf{1}_{\mathcal{O}_{\varepsilon,C_*\varepsilon}}(x_1-y_1,x_2-y_2,\eta_m(x_1,x_2)-y_3, v) \notag\\
 && \hspace{3cm} \times  |n_m(x-y)\cdot v|\sqrt{1+|\nabla\eta_m(x_1,x_2)|}dx_1dx_2 dv \\
 && \leq \int_{[-C_*\varepsilon,C_*\varepsilon]^2 } \Big[\int_{-s_*C_2 \sqrt{\varepsilon} \leq n_m(0,0)\cdot\frac{v}{|v|} \leq C_5\frac{1+|v|}{|v|}\varepsilon } + \int_{n_m(0,0)\cdot\frac{v}{|v|}\leq -s_*C_2 \sqrt{\varepsilon} } \Big] \cdots\\
 &&:= (\mathbf{I}) +(\mathbf{II}).
\end{eqnarray*}
We consider $(\mathbf{I})$. If $-s_* C_2 \sqrt{\varepsilon}\leq n_m(0,0)\cdot\frac{v}{|v|} \leq 0 $, then, $0\leq v_3 =-n_m(0,0)\cdot v \leq s_* C_2 |v| \sqrt{\varepsilon}$ and
\begin{eqnarray*}
0\leq v_3 \leq 2s_*C_2 \sqrt{|v_1|^2+|v_2|^2}\sqrt{\varepsilon},\quad \text{for}\quad 0<\varepsilon\ll 1.
\end{eqnarray*}
Moreover,
\begin{eqnarray*}
|n_m(x-y)\cdot v| &\leq& |n_m(0,0)\cdot v|+||n_m||_{C^1[-C_*\varepsilon,C_* \varepsilon]^2} (C_*+1/\tilde{C})|v|\varepsilon \\
& \leq & s_* C_2 |v|\sqrt{\varepsilon} +4 ||\eta_m||_{C^2[-C_*\varepsilon,C_* \varepsilon]^2} (C_*+1/\tilde{C})|v|\varepsilon.
\end{eqnarray*}
If $n_m(0,0) \cdot\frac{v}{|v|} \leq C_5 \frac{1+|v|}{|v|} \varepsilon$ then for $0<\varepsilon\ll 1$,
\begin{eqnarray*}
|v_3|=|n_m(0,0)\cdot v| \leq 2 C_5 (1+\sqrt{|v_1|^2+|v_2|^2})\varepsilon.
\end{eqnarray*}
Therefore,
\begin{eqnarray}
(\mathbf{I}) &\leq &\int_{[-C_*\varepsilon,C_*\varepsilon]^2 } \int_{0\leq v_3 \leq 2s_*C_2 \sqrt{|v_1|^2+|v_2|^2} \sqrt{\varepsilon}}
 \Big\{s_* C_2 |v|\sqrt{\varepsilon} +4 ||\eta_m||_{C^2[-C_*\varepsilon,C_* \varepsilon]^2} (C_*+1/\tilde{C})|v|\varepsilon \Big\} \notag\\
& & +  \int_{[-C_*\varepsilon,C_*\varepsilon]^2 } \int_{|v_3| \leq 2C_5 (1+\sqrt{|v_1|^2+|v_2|^2}) \varepsilon} \cdots \notag\\
&\lesssim & m_2([-C_*\varepsilon,C_*\varepsilon]^2 )\times \sqrt{\varepsilon}\iint_{V'}\sqrt{|v_1|^2+|v_2|^2}dv_1dv_2 \int_0^{2s_*C_2\sqrt{|v_1|^2+|v_2|^2}\sqrt{\varepsilon}}dv_3 \notag\\
&& + m_2([-C_*\varepsilon,C_*\varepsilon]^2 )\times\iint_{V'}\sqrt{|v_1|^2+|v_2|^2} dv_1dv_2 \int_0^{2C_5(1+\sqrt{|v_1|^2+|v_2|^2})\varepsilon}dv_3\notag\\
&\lesssim & \varepsilon^3, \label{3Estimate of I}
\end{eqnarray}
where $V'$ is the projection of $V$ onto the space $(v_1,v_2) \in \mathbb{R}^2$, which is also bound.

Now we decompose $(\mathbf{II})$ according to Lemma \ref{3Classification of velocity space}:
\begin{eqnarray*}
 (\mathbf{II})= \int_{[-C_*\varepsilon,C_*\varepsilon]^2 } \int_{|v|\leq \varepsilon^{1/3}} +\int_{[-C_*\varepsilon,C_*\varepsilon]^2 } \int_{-1\leq n_m(0,0)\cdot\frac{v}{|v|}\leq -s_*C_2\sqrt{\varepsilon} ~\text{and}~|v|\geq \varepsilon^{1/3}}.
\end{eqnarray*}
The first term is clearly bounded by $O(1) \varepsilon^3$. For the second term, we use (\ref{3Estimate of angle}) to have
\begin{eqnarray*}
 \{-1\leq n_m(0,0)\cdot\frac{v}{|v|}\leq -s_*C_2\sqrt{\varepsilon} ~\text{and}~|v|\geq \varepsilon^{1/3}\} \subset
 \{|n_m(0,0)\cdot\frac{v}{|v|}| \leq C_4\sqrt{\varepsilon} ~\text{and}~|v|\geq \varepsilon^{1/3}\}.
\end{eqnarray*}
So, we follow the same proof for (\ref{3Estimate of I}) to obtain
\begin{eqnarray}
(\mathbf{II}) &\lesssim& \varepsilon^3 + \int_{[-C_*\varepsilon,C_*\varepsilon]^2 } \int_{|n_m(0,0)\cdot\frac{v}{|v|}| \leq C_4\sqrt{\varepsilon} } \{C_4 |v|\sqrt{\varepsilon}+4||\eta_m||_{C^2([-C_*\varepsilon,C_*\varepsilon]^2)}(C_*+1/\tilde{C})|v|\varepsilon \}\notag\\
&\lesssim & \varepsilon^3. \label{3Estimate of II}
\end{eqnarray}
We conclude the estimate of (\ref{3Measure on Boundary}) form (\ref{3Estimate of I}) and (\ref{3Estimate of II}).

\underline{\it Proof of (\ref{3Measure on Outing Boundary})}. Due to the properties of the standard mollifier (\ref{3Definition of Mollifier1}), we obtain
\begin{eqnarray*}
&&\iint_{x \in \partial \Omega,n(x)\cdot v <0} [1-\chi_{\varepsilon}(x,v)] |n(x)\cdot v|dS_x dv\\
&&=\iint_{x \in \partial \Omega,n(x)\cdot v <0}\iint_{\mathbb{R}^3\times \mathbb{R}^3} [1-\mathbf{1}_{\overline{\Omega}\times V \setminus \mathcal{O}_{\varepsilon,C_*\varepsilon}}(x-y,v-u)]\psi_{\varepsilon}(y,u) |n(x)\cdot v|dudy dS_x dv\\
&&\leq\iint_{\mathbb{R}^3\times \mathbb{R}^3}\psi_{\varepsilon}(y,u) dudy \iint_{x \in \partial \Omega,n(x)\cdot v <0} \mathbf{1}_{\mathcal{O}_{\varepsilon,C_*\varepsilon}}(x-y,v-u) |n(x)\cdot v| dS_x dv\\
&&\leq\iint_{B_{\mathbb{R}^6}(0;\varepsilon/\tilde{C})}\psi_{\varepsilon}(y,u) dudy
 \iint_{x \in \partial \Omega,n(x)\cdot v <0} \mathbf{1}_{\mathcal{O}_{\varepsilon,C_*\varepsilon}}(x-y,v-u) |n(x)\cdot v| dS_x dv.
\end{eqnarray*}
Since $\sqrt{|y|^2+|u|^2}\leq \varepsilon/\tilde{C}$ and $n(x)\cdot v \leq 0$, we have
\begin{eqnarray*}
n(x)\cdot v &=& n(x-y)\cdot(v-u)+ (n(x)-n(x-y)) \cdot v +n(x-y)\cdot u\\
&=& n(x-y)\cdot(v-u)+O(\varepsilon/\tilde{C})(1+|v|).
\end{eqnarray*}
Therefore, we use (\ref{3Measure on Boundary}) to bounded (\ref{3Measure on Outing Boundary}) further as
\begin{eqnarray*}
&& \int_{\gamma_-}[1-\chi_{\varepsilon}(x,v)]d\gamma = \iint_{x\in \partial \Omega, n(x)\cdot v <0}| 1- \chi_{\varepsilon}(x,v)|\, |n(x)\cdot v|dS_x dv \\
&& \leq   \iint_{B_{\mathbb{R}^6}(0;\varepsilon/\tilde{C})} \psi_{\varepsilon}(y,u)dudy  \\
&& \hspace{1cm} \times \Big[
  \iint_{x \in \partial \Omega,n(x)\cdot v <0}\mathbf{1}_{ \mathcal{O}_{\varepsilon,C_*\varepsilon}}(x-y,v-u) |n(x-y)\cdot (v-u)| dS_x dv \Big] \\
&& \hspace{1cm} + O(\frac{\varepsilon}{\tilde{C}})\times m_3(\partial \Omega)  \times \int_{V}(1+|v|) dv\\
& & \leq_{\Omega,V} \varepsilon.
\end{eqnarray*}

\underline{\it Proof of (\ref{3Measure on Outing Boundary Partial})}. Following the same proof of (\ref{3Measure on Outing Boundary}), we deduce that
\begin{eqnarray*}
&&\int_{\gamma_-}|\partial \chi_{\varepsilon}(x,v)|d\gamma  = \iint_{x\in \partial \Omega, n(x)\cdot v <0}|\partial \chi_{\varepsilon}(x,v)| |n(x)\cdot v|dS_x dv \\
&&= \iint_{x\in \partial \Omega, n(x)\cdot v <0}|\partial [1-\chi_{\varepsilon}(x,v)]|\,|n(x)\cdot v|dS_x dv \\
&&= \iint_{x\in \partial \Omega, n(x)\cdot v <0}|\iint_{\mathbb{R}^3\times\mathbb{R}^3}\mathbf{1}_{\mathcal{O}_{\varepsilon,C_*\varepsilon}}(y,u)  \partial \psi_{\varepsilon}(x-y,v-u)dudy|\,|n(x)\cdot v|dS_x dv \\
&&\leq
 \iint_{B_{\mathbb{R}^6}(0;\varepsilon/\tilde{C})}|\partial \psi_{\varepsilon}(y,u)| dudy \\
&& \hspace{1cm}\times \Big[\iint_{x \in \partial \Omega,n(x)\cdot v <0}\mathbf{1}_{ \mathcal{O}_{\varepsilon,C_*\varepsilon}}(x-y,v-u)|n(x-y)\cdot (v-u)| dS_x dv \\
&& \hspace{1cm}+  O(\frac{\varepsilon}{\tilde{C}}) \times m_3(\partial \Omega) \times \int_{V} dv\Big]\\
&& \lesssim \frac{1}{\varepsilon} \sup_{(y,u)\in B_{\mathbb{R}^6}(0;\varepsilon/\tilde{C})} \Big[\iint_{x \in \partial \Omega,n(x)\cdot v <0}\mathbf{1}_{ \mathcal{O}_{\varepsilon,C_*\varepsilon}}(x-y,v-u) |n(x-y)\cdot (v-u)| dS_x dv \\
&& \hspace{1cm}+  O(\frac{\varepsilon}{\tilde{C}}) \times m_3(\partial \Omega) \times \int_{V}(1+|v|)dv\Big]\\
&& \lesssim 1.
\end{eqnarray*}
The proof of Proposition \ref{3Measure for singualr set on boundary} is completed.
 $\hfill\square$

\subsection{New Trace Theorem via the Double Iteration}
In this section we prove the following geometric result. For the later purpose, we state the result for
the sequence of solutions.

\begin{pro} \label{33Estimate of sequence} Let $h_0 \in L^1(\Omega \times V)$. Let $(h^m)_{m\geq 0} \subset L^{\infty}([0,T]; L^1(\Omega \times V))\cap L^1([0,T]; L^1(\gamma_+, d\gamma))$
solve
\begin{eqnarray}
(\partial_t + v \cdot \nabla_x + \Sigma ) h^{m+1} = H^m,\qquad  h^{m+1}|_{t=0} = h_0, \label{33Equation of Series}
\end{eqnarray}
where $\Sigma = \Sigma(x, v) \geq 0$, and such that the following inequality holds for all $0\leq t\leq T$ and $(x, v) \in \gamma_-$,
\begin{eqnarray}\label{33Assumption on Boundary}
|h^{m+1}(t, x, v)| \leq  C_1  \Big(1+\frac{1}{|n(x)\cdot v|}\Big) \bigg[\int_{n(x)\cdot v'>0}|h^m(t,x,v')| \{n(x)\cdot v'\} dv'+ R^  m\bigg],
\end{eqnarray}
where $H^m \in L^1([0, T]; L^1(\Omega\times V))$ and $R^m \in L^1([0, T]; L^1(\partial\Omega\times V, dS_xdv))$.

Then for all $m \geq  1$, $h^{m+1}_{\gamma_-}
\in L^1([0, T]; L^1(\gamma_-,d\gamma))$ and satisfies, for $\tau, t \in [0, T]$ and $0 < \delta \ll 1$,
\begin{eqnarray}\label{33Estimate on In-Flow Boundary}
\int_{\t}^t |h^{m+1}(s)|_{\gamma_-,1} &\leq& O(\delta) \int_{\t}^t |h^{m-1}(s)|_{\gamma_+,1} + C_{\delta} ||h(\t)||_1 \notag\\
 && + C_{\delta} \max_{i=m,m-1} \int_{\t}^t\Big\{ ||h^i(s)||_1+|| R^i(s)||_1 + ||H^i(s)||_1\Big\}.
\end{eqnarray}
\end{pro}

The proof of this proposition requires the following lemma:

\begin{lem}\label{33Estimate of coordinate transform} Let $\Omega \subset \mathbb{R}^3$ be an open bounded set with a smooth boundary $\partial \Omega$.
For $k \in \mathds{N}$, consider the map
\begin{eqnarray*}
\Phi_k: \{(x, v) \in \gamma_+ :  n(x_{\mathbf{b}}(x, v)) \cdot v < -1/k\} &\rightarrow & \{(x_{\mathbf{b}}, v) \in \gamma_- : n(x_{\mathbf{b}}) \cdot v < -1/k \} ,\\
(x, v) &\rightarrow & \Phi_k(x, v) := (\tilde{x}, v) := (x_{\mathbf{b}}(x, v), v).
\end{eqnarray*}
Then $\Phi_k $is one-to-one and we have a change of variables formula for all $k \in \mathds{N}$ :
\begin{eqnarray*}
\mathbf{1}_{\{n(\tilde{x})\cdot v <-1/k \}}|n(\tilde{x})\cdot v| dS_{\tilde{x}}dv &=& \mathbf{1}_{\{n(x_{\mathbf{b}}(x,v)\cdot v <-1/k)\}}|n(x)\cdot v| dS_xdv.
\end{eqnarray*}
\end{lem}

{\it Proof of lemma \ref{33Estimate of coordinate transform}}: This lemma deals with the change of variable formula: $ (x,v) \rightarrow (x_{\mathbf{b}(x,v)}, v)$. The proof is the same as in \cite{[GKTT2]}. We omit it here. $\hfill\square\\$

{\bf Proof of Proposition \ref{33Estimate of sequence}}. We now prove the estimate (\ref{33Estimate on In-Flow Boundary}).
Using (\ref{33Assumption on Boundary}), we obtain
\begin{eqnarray*}
\int_{\tau}^t|h^{m+1}(s)|_{\gamma_-,1}:= \int_{\tau}^t\iint_{n(x)\cdot v <0} |h^{m+1}(s,x,v)||n(x)\cdot v|dS_xdv ds \lesssim (A)+(B),   \end{eqnarray*}
where
\begin{eqnarray*}
(A)&:=& \int_{\t}^t\iint_{n(x)\cdot v > 0}|h^{m}(s,x,v)|\, |n(x)\cdot v|dS_xdv ds ,\\
(B)&:=&  \int_{\t}^t\iint_{n(x)\cdot v < 0}|R^{m}(s,x,v)| \,[1+|n(x)\cdot v|]dS_xdv ds .
\end{eqnarray*}
Clearly the last term (B) is bounded by the RHS of (\ref{33Estimate on In-Flow Boundary}).

We focus on (A) in the following. We split the outgoing part as
$\gamma_+ = \gamma_+^{\delta} \cup (\gamma_+ \setminus \gamma_+^{\delta}), $  where the almost grazing set $\gamma_+^{\delta}$ is defined in (\ref{1Almost grazing set}) and the non-grazing set $\gamma_+ \setminus \gamma_+^{\delta}$ is defined in (\ref{1Non-grazing set}).
Due to Lemma \ref{Trace theorem}, the non-grazing part $\gamma_+ \setminus \gamma_+^{\delta}$ of the integral is bounded as
\begin{eqnarray}\label{33Estimate on non-grazing set}
\int_0^t \int_{\gamma_+ \setminus \gamma_+^{\delta}} |h^m(s)|d\gamma &\lesssim_{t,\delta,\Omega}& ||h^m(\t)||_1 +\int_{\t}^t \{||h^m(s)||_1+||[\partial_t+v\cdot \nabla_x +\Sigma] h^m(s) ||_1 \}ds \notag\\
& \lesssim_{t,\delta,\Omega}& ||h^m(\t)||_1+\int_{\t}^t||h^m(s)||_1 +\int_{\t}^t||H^{m-1}(s)||_1.
\end{eqnarray}
It is also bounded by the RHS of (\ref{33Estimate on In-Flow Boundary}).

Now, we deal with the almost grazing set $\gamma_+^{\delta} $. We claim that the following truncated term with a number $k \in \mathds{N}$ is uniformly bounded in $k$ as follows:
\begin{eqnarray}\label{33Claim for Almost Grazing set}
&&\int_{\t}^t\iint_{x\in\partial\Omega,n(x)\cdot v>0} \mathbf{1}_{\{(x,v) \in \gamma_+^{\delta}\}}\mathbf{1}_{\{1/k<|n(x_{\mathbf{b}}(x,v))\cdot v|\}} |h^m(s,x,v)| \,\{n(x)\cdot v\} dvdS_x ds \notag\\
&& \leq O(\delta) \int_{\t}^t |h^{m-1}(s)|_{\gamma_+,1}+C_{\delta}\left[||h_0||_1+\int_{\t}^t\Big(||h^{m-1}(s)||_1+ ||H^{m-1}(s)||_1 +t||R^{m-1}||_1\Big) \right].
\end{eqnarray}
In order to show (\ref{33Claim for Almost Grazing set}), we use the Duhamel formula of the equation (\ref{33Equation of Series})
together with (\ref{33Assumption on Boundary}): for $(x, v) \in \gamma_+^{\delta}$ and $1/k < |n(x_{\mathbf{b}}(x, v)) \cdot v|$
\begin{eqnarray*}
&&|h^m(s,x,v)|  \mathbf{1}_{\{(x,v) \in \gamma_+^{\delta}\}}\mathbf{1}_{\{1/k<|n(x_{\mathbf{b}}(x,v))\cdot v|\}}\\
&&\hspace{2mm} \leq  \mathbf{1}_{\{s-t_{\mathbf{b}}(x,v)<\t \}}\mathbf{1}_{\{1/k<|n(x_{\mathbf{b}}(x,v))\cdot v|\}}|h^m(\t,x-(s-\t)v,v)|\\
&& \hspace{5mm}+ \mathbf{1}_{\{1/k<|n(x_{\mathbf{b}}(x,v))\cdot v|\}}\int_{\max\{\t,s-t_{\mathbf{b}}(x,v)\}}^{s}|H^{m-1}(\tau',x-(s-\tau')v,v)|d\tau'\\
&&\hspace{5mm} +  \mathbf{1}_{\{s-t_{\mathbf{b}}(x,v)>\t \}}\mathbf{1}_{\{1/k<|n(x_{\mathbf{b}}(x,v))\cdot v|\}} C_1 \Big(1+\frac{1}{n(x_{\mathbf{b}}(x,v))\cdot v}\Big)\\
& & \hspace{2cm} \times \int_{n(x_{\mathbf{b}}(x,v))\cdot v_1>0 }|h^{m-1}(s-t_{\mathbf{b}}(x,v),x_{\mathbf{b}}(x,v),v_1 )| n(x_{\mathbf{b}}(x,v)\cdot v_1)dv_1\\
&& \hspace{5mm}+ \mathbf{1}_{\{s-t_{\mathbf{b}}(x,v)>\t \}}\mathbf{1}_{\{1/k<|n(x_{\mathbf{b}}(x,v))\cdot v|\}} \Big(1+\frac{1}{|n(x_{\mathbf{b}}(x,v))\cdot v|} \Big)|R^{m-1}(s-t_{\mathbf{b}}(x,v),x_{\mathbf{b}}(x,v),v)|.
\end{eqnarray*}

We plug this estimate into the left hand side of (\ref{33Claim for Almost Grazing set}) to have
\begin{eqnarray}
&& \int_{\t}^t\iint_{x\in\partial\Omega,n(x)\cdot v>0} \mathbf{1}_{\{(x,v) \in \gamma_+^{\delta}\}}\mathbf{1}_{\{1/k<|n(x_{\mathbf{b}}(x,v))\cdot v|\}} |h^m(s,x,v)|\{n(x)\cdot v \}dvdS_x ds \notag\\
&&\hspace{3mm} \leq \int_{\t}^t\iint_{\gamma_+^{\delta}} \mathbf{1}_{\{s- t_{\mathbf{b}}(x,v)<\t \}} \mathbf{1}_{\{1/k<|n(x_{\mathbf{b}}(x,v))\cdot v|\}}|h^m(\t,x-(s-\t)v,v)|\{n(x)\cdot v\} dvdS_x ds \label{33Estimate of Initial Data}\\
&&\hspace{5mm} +\int_{\t}^t\iint_{\gamma_+^{\delta}} \mathbf{1}_{\{1/k<|n(x_{\mathbf{b}}(x,v))\cdot v|\}}\int_{\max\{\t,s-t_{\mathbf{b}(x,v)}\}}^{s}|H^{m-1}(\tau',x-(s-\tau')v,v)|\{ n(x)\cdot v\}d\tau' dvdS_x ds \label{33Estimate of Nonlinear term}\\
&&\hspace{5mm}+ \int_{\t}^t\iint_{\gamma_+^{\delta}} \mathbf{1}_{\{1/k<|n(x_{\mathbf{b}}(x,v))\cdot v|\}}  \frac{|n(x)\cdot v|}{|n(x_{\mathbf{b}}(x,v))\cdot v|}
 \int_{n(x_{\mathbf{b}}(x,v))\cdot v_1>0 } \mathbf{1}_{\{s-t_{\mathbf{b}}(x,v) >\t \} }\notag\\
 & & \hspace{1cm} \times|h^{m-1}(s-t_{\mathbf{b}}(x,v),x_{\mathbf{b}}(x,v),v_1 )|\, |n(x_{\mathbf{b}}(x,v)\cdot v_1)|dv_1dS_xdvds \label{33Estimate of subsequence}\\
&& \hspace{5mm}+ \int_{\t}^t\iint_{\gamma_+^{\delta}} \mathbf{1}_{\{s-t_{\mathbf{b}}(x,v)>\t \}}\mathbf{1}_{\{1/k<|n(x_{\mathbf{b}}(x,v))\cdot v|\}} \frac{ |n(x)\cdot v|}{|n(x_{\mathbf{b}}(x,v))\cdot v|} \notag \\
&&\hspace{1cm} \times |R^{m-1}(s-t_{\mathbf{b}}(x,v),x_{\mathbf{b}}(x,v),v)|dvdS_xds.\label{33Estimate of Remainder term}
\end{eqnarray}

\underline{\it Estimate of (\ref{33Estimate of Initial Data})}: Note that $x \in \partial\Omega$ in (\ref{33Estimate of Initial Data}). Without loss of generality we may assume that there exists $\eta: \mathbb{R}^2 \rightarrow \mathbb{R}$ such that $x^3 = \eta(x^1, x^2)$. We apply the following change of variables:
for fixed $v \in \mathbb{R}^3$,
\begin{eqnarray*}
  (x^1, x^2; s) \mapsto  y = (x^1 - (s-\t)v^1, x^2 - (s-\t) v^2, \eta(x^1, x^2) - (s-\t)v^3).
\end{eqnarray*}
 It maps $ \mathbb{R}^2 \times \{0 \leq  s-\t \leq t_{\mathbf{b}}(x, v)\} $ into $ \overline{\Omega}$.  We compute the Jacobian:
\begin{eqnarray*}
\det\left(\frac{\partial(y^1,y^2,y^3)}{\partial(x^1,x^2,s)}\right) &=&
\det\left(\begin{array}{ccc} 1&0&-v^1\\0&1&-v^2\\ \partial_{x^1}\eta(x^1,x^2)& \partial_{x^2}\eta(x^1,x^2)&-v^3 \end{array}\right)\\
&=& v\cdot\left(\begin{array}{cc} \partial_{x^1}\eta\\ \partial_{x^2}\eta \\-1\end{array}\right) =
v\cdot n\sqrt{1+|\partial_{x^1}\eta|^2+|\partial_{x^2}\eta|^2}.
\end{eqnarray*}
Therefore, such mapping $(x^1,x^2,s-\t) \rightarrow y$ is one-to-one when $ (x^1,x^2) \in \gamma_+$ and
\begin{eqnarray*}
\{n(x)\cdot v\}dS_xds =  \{n(x)\cdot v\} \sqrt{1+|\partial_{x^1}\eta|^2+|\partial_{x^2}\eta|^2}dx^1dx^2 ds =dy =dy^1dy^2dy^3,
\end{eqnarray*}
and
\begin{eqnarray}\label{33Estimate of Inital Data F}
(\ref{33Estimate of Initial Data}) &\leq& \int_{V}\int_{\t}^t \int_{\partial\Omega} \mathbf{1}_{\{(x,v)\in \gamma_+\}} \mathbf{1}_{\{1/k <|n(x_{\mathbf{b}(x,v)})\cdot v |\}} |h^m(\t,x-(s-\t)v,v)|\, |n(x)\cdot v|dS_xdsdv \notag\\
& \leq & \int_{V} \int_{\Omega} |h^m(\t, y,v)| dydv
\leq ||h(\t)||_1.
\end{eqnarray}

\underline{\it Estimate of (\ref{33Estimate of Nonlinear term})}: Considering the region of $\{(\tau', s) \in [\t, t] \times [\t, t] : \max\{\t, s - t_{\mathbf{b}}(x, v)\} \leq \tau' \leq s\}$, it is bounded by
\begin{eqnarray}\label{33Estimate of Nonlinear term Stp1}
\int_{V} dv \int_{\t}^t d\tau' \int_{\tau'}^{\min\{t,\tau'+ t_{\mathbf{b}}(x, v)\}} \int_{\partial\Omega} |H^{m-1}(\tau',x-(s-\tau')v,v)|\,|n(x)\cdot v|dS_x ds .
\end{eqnarray}
Note that $x \in \partial \Omega$, without loss of generality, we may assume that $x^3 = \eta(x^1, x^2)$ for $¦Ç :
\mathbb{R}^2 \mapsto \mathbb{R}$. We apply the change of variables: for fixed $v \in V$ and $\tau' \in[0, t]$,
\begin{eqnarray*}
(x^1, x^2; s) \mapsto y \equiv (x^1-(s-\tau')v^1, x^2-(s-\tau')v^2, \eta(x^1, x^2)-(s-\tau')v^3).
\end{eqnarray*}
Clearly, it maps $ \mathbb{R}^2 \times [\tau', \min\{t, \tau'+ t_{\mathbf{b}}(x, v)\}]$ into $ \overline{\Omega}$ since $0\leq (s-\tau')\leq t_{\mathbf{b}}(x,v)$.  The Jacobian of this change variable is $\{v\cdot n(x)\}\sqrt{1+|\partial_{x^1}\eta|^2+|\partial_{x^2}\eta|^2}$ and $\{v\cdot n(x)\}dS_x ds \leq dy$. Applying the change
of variables to (\ref{33Estimate of Nonlinear term Stp1}) to have
\begin{eqnarray}\label{33Estimate of Nonlinear term F}
 (\ref{33Estimate of Nonlinear term}) \leq \int_{\t}^t\int_{V}\int_{\Omega} |H^{m-1}(\tau',y,v)| dy dvd\tau' = \int_{\t}^t ||H^{m-1}(\tau')||_1 d\tau'.
\end{eqnarray}

\underline{\it Estimate of (\ref{33Estimate of subsequence})}: This part is the most delicate. We rewrite (\ref{33Estimate of subsequence}) as
\begin{eqnarray}\label{33Estimate of subsequence Step1}
 && \int_{\t}^t ds\int_{\partial \Omega}dS_x\int_{V}dv\int_{V}dv_1 \mathbf{1}_{\{(x,v)\in \gamma_+^{\delta}\}} \mathbf{1}_{\{n(x_{\mathbf{b}}(x,v))\cdot v_1 >0\}} \mathbf{1}_{\{s- t_{\mathbf{b}}(x,v)>\t\}} \mathbf{1}_{\{1/k<|n(x_{\mathbf{b}}(x,v))\cdot v|\}} \notag  \\
  & & \hspace{1cm}\times \frac{|n(x)\cdot v|}{|n(x_{\mathbf{b}}(x,v))\cdot v|}|n(x_{\mathbf{b}}(x,v))\cdot v_1|
 |h^{m-1}(s-t_{\mathbf{b}}(x,v),x_{\mathbf{b}}(x,v),v_1 )|.
\end{eqnarray}
We apply the following change of variables
\begin{eqnarray*}
s \in [0, t] \mapsto \tilde{s} = s - t_{\mathbf{b}}(x, v) \in [\t, t - t_{\mathbf{b}}(x, v)] \subset [\t,t], \label{Change of variable1}
\end{eqnarray*}
where we have used the fact that $s\in [t_{\mathbf{b}}(x, v)+\t, t]$. Clearly the Jacobian is $1$ so
that $d\tilde{s} = ds$ and hence
\begin{eqnarray}\label{33Estimate of subsequence Step2}
&& (\ref{33Estimate of subsequence Step1}) \leq \int_{\t}^{t} d\tilde{s}\int_{\partial \Omega}dS_x\int_{V}dv\int_{V}dv_1 \mathbf{1}_{\{(x,v)\in \gamma_+^{\delta}\}} \mathbf{1}_{\{n(x_{\mathbf{b}}(x,v))\cdot v_1 >0\}} \mathbf{1}_{\{1/k<|n(x_{\mathbf{b}}(x,v))\cdot v|\}} \notag  \\
  & & \hspace{2cm}\times  \frac{|n(x)\cdot v|}{|n(x_{\mathbf{b}}(x,v))\cdot v|}|n(x_{\mathbf{b}}(x,v))\cdot v_1|
 |h^{m-1}(\tilde{s},x_{\mathbf{b}}(x,v),v_1 )|.
\end{eqnarray}
Let us denote
%\begin{eqnarray}
$\tilde{x}:= x_{\mathbf{b}}(x, v).$
%\end{eqnarray}
In the case $n_3(x_{\mathbf{b}}(x,v)) \neq 0$, there exists some function $\phi: \mathbb{R}^2 \rightarrow \mathbb{R}$ such that
\begin{eqnarray*}
\tilde{x}= x_{\mathbf{b}}(x,v)=(\tilde{x}_1,\tilde{x}_2,\tilde{x}_3)=(\tilde{x}_1,\tilde{x}_2,\phi(\tilde{x}_1,\tilde{x}_2)) \in \partial \Omega.
\end{eqnarray*}
Note that since $(x, v) \in \gamma_+$ and $|n(x_{\mathbf{b}}(x, v))\cdot v| > 1/k$, from Lemma \ref{33Estimate of coordinate transform}, the mapping
$(x, v) \mapsto (\tilde{x}, v)$ is one-to-one.
We apply the change of variables of Lemma \ref{33Estimate of coordinate transform}: for $(x, v) \in \gamma^+$ and $|n(x_{\mathbf{b}}(x, v)) \cdot v| =
|n(\tilde{x}) \cdot v| > 1/k$, we apply the change of variables
\begin{eqnarray}\label{33Change of variable3}
(x, v) \mapsto (\tilde{x}, v) := (x_{\mathbf{b}}(x, v), v).
\end{eqnarray}
The Jacobian is
\begin{eqnarray*}
 \det\left( \frac{\partial(\tilde{x},v)}{\partial(x,v)} \right) =\det\left( \frac{\partial \tilde{x}}{\partial x} \right)=\left|\frac{n(x)\cdot v}{n(\tilde{x})\cdot v}\right| \frac{\sqrt{1+|\nabla \eta|^2}}{\sqrt{1+|\nabla \phi|^2}}
,~~\text{and}~~dS_{\tilde{x}}:= \left|\frac{n(x)\cdot v}{n(\tilde{x})\cdot v}\right| dS_x.
\end{eqnarray*}
At the same time, we have
\begin{eqnarray*}t_{\mathbf{b}}(x, v) &=& t_{\mathbf{b}}(x_{\mathbf{b}}, -v),\\
x &=& x_{\mathbf{b}}(x, v) + t_{\mathbf{b}}(x, v)v = x_{\mathbf{b}}(x, v) + t_{\mathbf{b}}(x_{\mathbf{b}}(x, v), -v)v
\\ &=& x_{\mathbf{b}}(x, v) - t_{\mathbf{b}}(x_{\mathbf{b}}(x, v), -v)(-v)
= \tilde{x} - x_{\mathbf{b}}(\tilde{x}, -v)(-v).
\end{eqnarray*}
So we rewrite $\{(x,v) \in \gamma_+^{\delta}\} $ as
\begin{eqnarray*}
\mathbf{1}_{\{(x,v) \in \gamma_+^{\delta}\}} = \mathbf{1}_{\{0<n(\tilde{x}-t_{\mathbf{b}}(\tilde{x},-v)(-v))\cdot v\leq \delta\}}.
\end{eqnarray*}

Then, from (\ref{33Estimate of subsequence Step2}),
\begin{eqnarray}
(\ref{33Estimate of subsequence Step2})&\leq &\int_{\t}^{t}d\tilde{s}\int_{V}dv_1\int_{V}dv\int_{\partial \Omega} dS_{\tilde{x}} \,\mathbf{1}_{\{0<n(\tilde{x}-t_{\mathbf{b}}(\tilde{x},-v)(-v))\cdot v\leq \delta \}}\notag \\
&& \hspace{1cm} \times \mathbf{1}_{\{n(\tilde{x}\cdot v_1 >0\}} \mathbf{1}_{\{1/k<|n(\tilde{x})\cdot v|\}} \,|n(\tilde{x})\cdot v_1|
 |h^{m-1}(\tilde{s},\tilde{x},v_1 )| \notag\\
 & \leq &\int_{\t}^t \int_{\gamma_+} |h^{m-1}(\tilde{s},\tilde{x},v_1 )| \, |n(\tilde{x})\cdot v_1|dS_{\tilde{x}} dv_1 d\tilde{s}\sup_{\tilde{x}\in \partial \Omega}\int_{V}\mathbf{1}_{\{0<n(\tilde{x}-t_{\mathbf{b}}(\tilde{x},-v)(-v))\cdot v\leq \delta \}} dv  \label{33Estimate of subsequence Step3}
\end{eqnarray}

Due to Lemma \ref{Guo Covering}, for given $\delta > 0$, there exists $\delta_{\delta,V}$ and $l_{\delta,\Omega,V}$ balls $B_{x_1;r_1},B(x_2;r_2),\cdots, B(x_l;r_l)$ covering $\overline{\Omega}$, as well as $l$ open sets $\mathcal{O}_{x_1},\mathcal{O}_{x_2},\cdots,\mathcal{O}_{x_l} \subset V$, with $m_3(\mathcal{O}_{x_i}) \leq \delta$ for all $1 \leq i\leq l$  such that
\begin{eqnarray*}
\sup_{\tilde{x}\in \partial \Omega}\int_{V}\mathbf{1}_{\{0<n(\tilde{x}-t_{\mathbf{b}}(\tilde{x},-v)(-v))\cdot v\leq \delta \}} dv &\leq & \max_{i}\sup_{\tilde{x} \in B_{\mathbb{R}^3}(x_1;r_i)}m_3 \{v \in V: |n_{\mathbf{b}}(\tilde{x},-v)\cdot (-v)|
\leq \delta \} \\
&\leq& \max_i m_3(\mathcal{O}_i) \leq \delta.
\end{eqnarray*}

Therefore, for $0 < \delta \ll 1,$  such that
\begin{eqnarray}\label{33Estimate of subsequence F}
(\ref{33Estimate of subsequence}) &\lesssim & O(\delta)\int_{\t}^t \int_{\partial\Omega}\int_V |h^{m-1}(\tilde{s},\tilde{x},v_1)|\, |n(\tilde{x})\cdot v_1| dS_{\tilde{x}}dv_1 d\tilde{s} \notag\\
&= & O(\delta) \int_{\t}^t |h^{m-1}(s)|_{\gamma_+,1}ds.
\end{eqnarray}
$\vspace{3pt}$

\underline{\it Estimate of (\ref{33Estimate of Remainder term})}: We apply the change of variables $\tilde{s}=s-t_{\mathbf{b}}(x,v) $ and (\ref{33Change of variable3}), then by using Lemma \ref{Trace theorem} to bound as
\begin{eqnarray}\label{33Estimate of Remainder term F}
(\ref{33Estimate of Remainder term}) \lesssim \int_{\t}^t\int_{\partial \Omega}\int_{V}  |R^{m-1}(\tilde{s},\tilde{x},v)|dS_{\tilde{s}}dvd\tilde{s} = \int_{\t}^t |R^{m-1}(s)|_1 ds.
\end{eqnarray}
Finally from (\ref{33Estimate of Inital Data F}), (\ref{33Estimate of Nonlinear term F}), (\ref{33Estimate of subsequence F}) and (\ref{33Estimate of Remainder term F}), we prove our claim (\ref{33Claim for Almost Grazing set}).

The last step is to pass a limit $k \rightarrow \infty$. Clearly the sequence is non-decreasing in $k$:
\begin{eqnarray*}
0 \leq \mathbf{1}_
{ 1/k <|n(x_{\mathbf{b}}(x,v))\cdot v|}|h^m(s, x, v)| \leq  \mathbf{1}_{ 1/(k+1) <|n(x_{\mathbf{b}}(x,v))\cdot v|}|h^m(s, x, v)|.
\end{eqnarray*}
For $ \varepsilon> 0$, we choose $k \gg 1$ such that $1/ k < \varepsilon $. Then
\begin{eqnarray*}
&&\int_{\gamma_+}
\Big[1-\mathbf{1}_{\{ |n(x_{\mathbf{b}}(x,v'))\cdot v'|>1/k \}}(x,v') \Big] d\gamma \\
&&\hspace{3mm}\leq  \int_{\partial \Omega} \int_{n(x_{\mathbf{b}}(x,v'))\cdot v'>0}  \mathbf{1}_{\{ |n(x_{\mathbf{b}}(x,v'))\cdot v'|<1/k \}} |n(x_{\mathbf{b}}(x,v'))\cdot v'|dv' dS_x \\
&&\hspace{3mm}\leq \frac{1}{k}\int_{\partial \Omega}  \int_{n(x_{\mathbf{b}}(x,v')) \cdot v'>0}  dv' dS_x \\
& & \hspace{3mm}\lesssim \varepsilon.
\end{eqnarray*}
It concludes that
\begin{eqnarray*}
\mathbf{1}_{\{\frac{1}{k}<|n(x_{\mathbf{b}})\cdot v|\}} h^{m}(s,x,v) \rightarrow |h^m(s,x,v)|,~~\,\,a.e.\,(x,v)\in \gamma_+ \, \text{with} \,d\gamma.
\end{eqnarray*}
Now, we use the monotone convergence theorem to conclude
\begin{eqnarray}
\int_0^t\int_{\gamma_+^{\delta}} \mathbf{1}_{\{1/k<|n(x_{\mathbf{b}}(x,v))\cdot v|\}}|h^m(s, x, v)| d\gamma ds \rightarrow \int_0^t \int_{\gamma_+^{\delta}} |h^m(s, x, v)| d\gamma ds,
\end{eqnarray}
as $k \rightarrow \infty $ and therefore $ \int_0^t \int_{\gamma_+^{\delta}} |h^m(s, x, v)| d\gamma ds$ has the same upper bound of (\ref{33Claim for Almost Grazing set}). Together with (\ref{33Estimate on non-grazing set}) we conclude (\ref{33Estimate on In-Flow Boundary}). $\hfill\square$

\subsection{Estimates of the total variation}
The purpose of this subsection is to prove Theorem 2. To give the estimate of solution in total variation, we use following approximation scheme.
For $u_0 \in BV (\Omega \times V)$ and $||u_0||_{\infty}<\infty $ we choose $u_0^{\varepsilon}\in BV (\Omega \times V) \cap  C^{\infty}(\Omega \times V)$
satisfying $||[\hat{u}_0 - u_0]||_{\infty} \rightarrow 0$ and $ ||\nabla_{x,v} \hat{u}_0||_1 \rightarrow ||u_0||_{\widetilde{BV}}$. At the same time, for fixed $0\leq t\leq T$, we choose $\hat{q} \in BV (\Omega \times V) \cap  C^{\infty}(\Omega \times V)$ satisfying $||[\hat{q}- q]||_{\infty} \rightarrow 0$ and $ ||\nabla_{x,v} \hat{q}||_1 \rightarrow ||q||_{\widetilde{BV}}$.

Consider the sequence $u^{\varepsilon, m}$ defined by $u^{\varepsilon,0} = \chi_{\varepsilon} u_0 $ and for all $m \geq 0$,
\begin{eqnarray}\label{34Approximation sequences}
\begin{array}{rlllll}
\partial_t u^{\varepsilon,m+1}+ v\cdot \nabla_x u^{\varepsilon,m+1} +\Sigma u^{\varepsilon,m+1} &=&\chi_{\varepsilon} [ K u^{\varepsilon,m}+\hat{q}],~~& \text{in} &~~\Omega\times V,\vspace{3pt}\\
u^{\varepsilon,m+1}(0,x,v)&=&\chi_{\varepsilon} \hat{u}_0(x,v), ~~&\text{in}&~~\Omega\times V,\vspace{3pt}\\
u^{\varepsilon,m+1} &= &\chi_{\varepsilon} \mathcal{P}_{\gamma}u^{\varepsilon,m} + \chi_{\varepsilon} r ,~~&\text{on}&~~\gamma_-,
\end{array}\end{eqnarray}
where $\chi_{\varepsilon}$ is defined in (\ref{3Definition of Mollifier3}).

In order to study the derives of $u^{\varepsilon, m+1}(t,x,v)$ with respect to $x,v$, we need to consider the derivatives on the boundary.
For the purpose of it, we assume that $u$ satisfies the following neutron transport equation with the diffusive-inflow boundary condition
\begin{eqnarray*}
(\partial_t + v\cdot \nabla_x +(\Sigma -K))u = q,~~\quad u|_{\gamma_-} = \mathcal{P}_{\gamma} u +r.
\end{eqnarray*}
Let ${\tau_1(x), \tau_2(x)} $ be a basis of the tangent space at $x \in\partial \Omega$ (therefore
${\tau_1(x), \tau_2(x), n(x)}$ is
an orthonormal basis of $\mathbb{R}^3$), i.e. $\tau_1(x) \cdot n(x) = \tau_2 \cdot n(x)=0$ and
$\tau_1 \times \tau_2 =n(x)$. Define the orthonormal transformation from $n(x), \tau_1, \tau_2$ to the standard bases $(e_1,e_2,e_3)$, i.e. $\mathcal{T} n(x) =e_1$, $\mathcal{T} \tau_1(x) =e_2$, $\mathcal{T} \tau_2 =e_2$  and $\mathcal{T}^{-1}=\mathcal{T}^t$. Upon a change of variable:
$\xi =\mathcal{T} v'$, we have
\begin{eqnarray}\label{34Orthonormal transformation}
 n(x) \cdot v' =n(x) \cdot \mathcal{T}^t \xi =n(x)\cdot n(x)^t \mathcal{T}^t \xi =[\mathcal{T}n(x)]^t \xi =e_1 \cdot \xi =\xi_1,
\end{eqnarray}
then denote $\partial_{\tau_i}$ to be the (tangential) $\tau_i$-directional derivative and $\partial_n$ to be the normal derivative.
For all $(x,v)\in \gamma_-$, both $t$ and
$v$ derivatives behave nicely for the diffusive boundary condition,
\begin{eqnarray}
(\partial_t u)|_{\gamma_-} &=& c  \int_{n(x)\cdot v'>0} \partial_t u(t,x,v')\{n(x)\cdot v'\}dv'+ \partial_t r,\label{34Derivative at the boundary t}\\
 (\nabla_v  u)|_{\gamma_-} &=& \nabla_v r.\label{34Derivative at the boundary v}
\end{eqnarray}
From the choose of $\mathcal{T}$ in (\ref{34Orthonormal transformation}),
\begin{eqnarray*}
\int_{n(x)\cdot v'>0} u(t,x,v')\{n(x)\cdot v'\}dv' = \int_{\xi_1>0} u(t,x,\mathcal{T}^t(x) \xi)\,\xi_1 d\xi,
\end{eqnarray*}
So, we can further take the tangential derivatives $\partial_{\tau_i}~(i=1,2)$ as, for $(x,v) \in \gamma _- $,
\begin{eqnarray} \label{34Derivative at the tangential}
&& \partial_{\tau_i} u(t,x,v)
 = c\int_{\xi_1>0} \Big[\partial_{\tau_i} u(t,x,\mathcal{T}^t\xi) + \nabla_v u(t,x,\mathcal{T}^t \xi)\frac{\partial \mathcal{T}^t (x)}{\partial \tau_i} \Big] \,\xi_1 d\xi + \partial_{\tau_i} r \notag\\
&&\hspace{18mm} = c \int_{n(x)\cdot v'>0}
\Big[ \partial_{\tau_i} u(t,x,v')+  \nabla_v u(t,x,v')\frac{\partial \mathcal{T}^t (x)}{\partial \tau_i} \mathcal{T}v'\Big] \notag\\
&& \hspace{2.5cm} \times \{n(x)\cdot v'\}dv' + \nabla r \frac{\partial \mathcal{T}^t (x)}{\partial \tau_i}.
\end{eqnarray}
The difficulty is always the control of the normal spatial derivative $\partial_n$. Near the boundary $\partial \Omega$, it is natural to use the original equation to solve $\partial_n u$ inside the region, in terms of $\partial_t u,~\nabla_v u$ and
$\partial_{\tau_i} u$ as
\begin{eqnarray*}
\partial_n u(t,x,v) =-\frac{1}{n(x)\cdot v} \Big\{ \partial_t u +\sum_{i=1}^2 (v\cdot \tau_i)\partial_{\tau_i} u -(\Sigma- K)u+q\Big\}.
\end{eqnarray*}
 From (\ref{34Derivative at the boundary t}), (\ref{34Derivative at the boundary v}) and (\ref{34Derivative at the tangential}), we can express $\partial _n u$ at $(x,v) \in \gamma_- $ as
\begin{eqnarray}\label{34Normal derivative at the boundary}
&& \partial_n u(t,x,v) =  -\frac{c}{n(x)\cdot v} \bigg\{  \int_{n(x)\cdot v'>0} dv'\{n(x)\cdot v'\} \notag\\
 &&\hspace{2cm}\times \bigg( \partial_t u(t,x,v') + \sum_{i=1}^2(v\cdot \tau_i) \Big[\partial_{\tau_i} u(t,x,v')+  \nabla_v u(t,x,v')\frac{\partial \mathcal{T}^t (x)}{\partial \tau_i}\mathcal{T}v' \Big]\bigg) \\
 &&\hspace{2cm}  + \partial_t r + \sum_{i=1}^2(v\cdot \tau_i)  \nabla r \frac{\partial \mathcal{T}^t (x)}{\partial \tau_i}  - [(\Sigma- K)u- q]|_{\gamma_-} \bigg\}.\notag
\end{eqnarray}
Moreover,  the equation gives
\begin{eqnarray*}
 \partial_t u(t,x,v') =q -\Big[\sum_{i=1}^2 (v'\cdot \tau_i) \partial_{\tau_i} + (v'\cdot n) \partial_n  + (\Sigma-K)\Big] u(t,x,v').
\end{eqnarray*}
Submitting this equality into (\ref{34Normal derivative at the boundary}), it derives to
\begin{eqnarray}\label{34Normal derivative at the boundaryF}
&& \partial_n u(t,x,v) =  -\frac{ c }{n(x)\cdot v} \bigg\{ \int_{n(x)\cdot v'>0} \{n(x)\cdot v'\}dv' \notag\\
 &&\hspace{3cm}\times \bigg( q -\Big[\sum_{i=1}^2 (v'\cdot \tau_i) \partial_{\tau_i} + (v'\cdot n) \partial_n  + (\Sigma-K)\Big]u(t,x,v') \notag\\
 &&\hspace{3.5cm}+ \sum_{i=1}^2(v\cdot \tau_i) \Big[\partial_{\tau_i} u(t,x,v')+  \nabla_v u(t,x,v')\frac{\partial \mathcal{T}^t (x)}{\partial \tau_i}\mathcal{T}v' \Big]\bigg) \\
 &&\hspace{3cm}  + \partial_t r + \sum_{i=1}^2(v\cdot \tau_i)  \nabla r \frac{\partial \mathcal{T}^t (x)}{\partial \tau_i}  - [(\Sigma- K)u- q]|_{\gamma_-} \bigg\}.\notag
\end{eqnarray}

We firstly study the estimates of the derivatives of the solution for the following simpler neutron transport equation with in-flow boundary condition
\begin{eqnarray}\label{34Simper approximate inf-low Boundary}
u_t + v\cdot \nabla u +\Sigma u= \mathbf{Q},\quad u(t, x, v)|_{\gamma_-} = \mathbf{R}(t, x, v),\quad u(0,x,v)=u_0(x,v).
\end{eqnarray}

\begin{lem}\label{34Estimate of approximate sequaence} Assume $\mathcal{U}$ is an open subset of $\mathbb{R}^3 \times \mathbb{R}^3$ such that $\mathfrak{S}_B\subset \mathcal{U}$. For $ (t, x, v) \in [0, T] \times \{\mathcal{U} \cap (\overline{\Omega} \times V\}$, we assume
\begin{eqnarray}
u_0(x, v) \equiv 0,\, \mathbf{R}(t, x, v) \equiv 0, \,\mathbf{Q}(t, x, v) \equiv 0.
\end{eqnarray}
Assume further that
\begin{eqnarray*}
u_0\in L^{\infty}(\Omega\times V), \quad \mathbf{R} \in L^{\infty}([0, T] \times \gamma_-), \quad  \mathbf{Q} \in L^{\infty}([0, T] \times \Omega \times V),
\end{eqnarray*}
and
\begin{eqnarray*}
&&\nabla_{x,v} u_0 \in L^1(\Omega\times V),\\
&& \partial_{\tau_i}\mathbf{R},\, \frac{1}{n(x)\cdot v}\Big\{ -[\partial_t + \sum_{i=1}^2(v\cdot\tau_i)\partial_{\tau_i}+\Sigma]\mathbf{R}+\mathbf{Q}\Big\}, \,\nabla_v \mathbf{R}, \,\nabla_{x,v}\Sigma \in L^1([0,T]\times \gamma_-),  \\
&&  \nabla_{x,v}\Sigma,\,\, \nabla_{x,v}\mathbf{Q}~\in L^1([0,T]\times \Omega\times V).
\end{eqnarray*}
Then there exists a unique solution $u$ to the transport equation (\ref{34Simper approximate inf-low Boundary}) such that $u\in C^0([0, T]\times \overline{\Omega}\times V)$ and $ \nabla_{x,v}u \in C^0([0, T],L^1(\Omega\times V))$ and
the traces satisfy
\begin{eqnarray*}
 \nabla_{x,v} u= \nabla_{x,v}\mathbf{R}, ~~\text{on}~~\gamma_-,
~~~~ \nabla_{x,v}u(0,x,v)=\nabla_{x,v}u_0(x,v),~~\text{in}~~\Omega\times V,
\end{eqnarray*}
where $\nabla_x \mathbf{R}$ is defined by
\begin{eqnarray*}
\nabla_x \mathbf{R} = \sum_{i=1}^2\tau_i\partial_{\tau_i} \mathbf{R} +\frac{n}{n\cdot v}\bigg\{-\Big[\partial_t+\sum_{i=1}^2(v\cdot \tau_i)\partial_{\tau_i}+ \Sigma \Big]\mathbf{R}+\mathbf{Q}\bigg\}.
\end{eqnarray*}
Moreover,
\begin{eqnarray}
%&&||\partial_t u(t)||_1 +\int_0^t |\partial_t u|_{\gamma_+,1}+\int_0^t||\Sigma \partial_t u||_1 \notag\\
%&& \hspace{1cm} =|| \partial_t u_0||_1 + \int_0^t |\nabla_x \mathbf{R}|_{\gamma_-,1}+\int_0^t\iint_{\Omega\times V} \text{sgn}(\partial_t u)\{ %\nabla_x \mathbf{Q}- \nabla_x \Sigma\, u\}, \label{34Estimates of t Partial Differential}\\
&&||\nabla_x u(t)||_1 +\int_0^t |\nabla_x u|_{\gamma_+,1}+\int_0^t||\Sigma \nabla_x u||_1 \notag\\
&& \hspace{1cm} =||\nabla_x u_0||_1+\int_0^t |\nabla_x \mathbf{R}|_{\gamma_-,1}+\int_0^t\iint_{\Omega\times V} \text{sgn}(\nabla_x u)\{ \nabla_x \mathbf{Q}- \nabla_x \Sigma\, u\}, \label{34Estimates of x Partial Differential}\\
&&||\nabla_v u(t)||_1 +\int_0^t |\nabla_v u|_{\gamma_+,1}+\int_0^t||\Sigma \nabla_v u||_1 \notag\\
&& \hspace{1cm} =||\nabla_v u_0||_1+\int_0^t |\nabla_v \mathbf{R}|_{\gamma_-,1}+\int_0^t\iint_{\Omega\times V} \text{sgn}(\nabla_v u)\{ \nabla_v \mathbf{Q}- \nabla_x u-\nabla_v \Sigma \,u\}. \label{34Estimates of v Partial Differential}
\end{eqnarray}
\end{lem}

\prof We use the Duhamel formula of $u$:
\begin{eqnarray}
u(t,x,v) &=& \mathbf{1}_{\{t< t_{\mathbf{b}}(x,v)\}} e^{-\int_0^t \Sigma(t-\tau,x-\tau v,v)d\tau} u_0(x-tv,v) \notag\\
 & & + \mathbf{1}_{\{t >t_{\mathbf{b}}(x,v)\}} e^{-\int_0^{t_{\mathbf{b}}(x,v)} \Sigma(t-\tau,x-\tau v,v)d\tau} \mathbf{R}(t-
 t_{\mathbf{b}}(x,v),x_{\mathbf{b}}(x,v),v) \label{34Expression of Solution}\\
& & + \int_0^{\min\{t,t_{\mathbf{b}}(x,v)\}}e^{-\int_0^s\Sigma(t-\tau,x-\tau v,v)d\tau}\mathbf{Q}(t-s,x-sv,v)ds. \notag
\end{eqnarray}
Recall the derivatives of $x_{\mathbf{b}}$ and $t_{\mathbf{b}}$ in Lemma \ref{Derivatives},
following Proposition 1 in \cite{[GKTT1]}, %the derivative of $u$ with respect to $t$ is
%\begin{eqnarray*}
%&&\partial_t u(t,x,v)\mathbf{1}_{\{t\neq t_{\mathbf{b}}\}}\notag \\
%&&=\mathbf{1}_{\{t < t_{\mathbf{b}}\}}  e^{-\int_0^t \Sigma(t-\tau,x-\tau v,v)d\tau}\bigg\{ \Sigma(0, x-tv,v) u_0 (x-tv,v)+v\cdot \nabla_x u_0
%\notag\\&&\hspace{2cm} -v \Big(\int_0^t\partial_t \Sigma(t-\tau,x-\tau v,v)d\tau \Big) u_0(x-tv,v) -\mathbf{Q}|_{t=0}\bigg\}\notag\\
%&&\hspace{3mm} + \mathbf{1}_{\{t >t_{\mathbf{b}}\}} e^{-\int_0^{t_{\mathbf{b}}} \Sigma(t-\tau,x-\tau v,v)d\tau}\bigg\{ %\int_0^{t_{\mathbf{b}}}\partial_t \Sigma(t-\tau,x-\tau v,v) + \partial_t \bigg\} \mathbf{R}(t-t_{\mathbf{b}},x_{\mathbf{b}},v)\\
%&&\hspace{3mm}- \int_0^{\min\{t,t_{\mathbf{b}}(x,v)\}}e^{-\int_0^s\Sigma(t-\tau,x-\tau v,v)d\tau}\Big(\int_0^s \partial_t \Sigma(t-\tau,x-\tau %v,v)d\tau + \partial_t\Big) \mathbf{Q}(t-s,x-sv,v)ds.\notag \label{34Expression of T Derives}
%\end{eqnarray*}
the derivative of $u$ with respect to $x,v$ are
{\small \begin{eqnarray*}
&&\nabla_x u(t,x,v)\mathbf{1}_{\{t\neq t_{\mathbf{b}}\}}\notag \\
&&=\mathbf{1}_{\{t < t_{\mathbf{b}}\}}  e^{-\int_0^t \Sigma(t-\tau,x-\tau v,v)d\tau}\bigg\{\nabla_x u_0(x-tv,v)-\int_0^t\nabla_x \Sigma(t-\tau,x-\tau v,v)d\tau  u_0(x-tv,v)\bigg\}\notag\\
&&\hspace{3mm} + \mathbf{1}_{\{t >t_{\mathbf{b}}\}} e^{-\int_0^{t_{\mathbf{b}}} \Sigma(t-\tau,x-\tau v,v)d\tau}\bigg\{\sum_{i=1}^2\tau_i \partial_{\tau_i}\mathbf{R} \notag\\
&& \hspace{2cm}-\frac{n(x_{\mathbf{b}})}{n(x_{\mathbf{b}})\cdot v}\Big\{\Big[\partial_t+\sum_{i=1}^2(v\cdot \tau_i)\partial_{\tau_i}+ \Sigma \Big]\mathbf{R}-\mathbf{Q}\Big\} \bigg\}(t-t_{\mathbf{b}},x_{\mathbf{b}},v)\notag\\
&&\hspace{3mm}- \mathbf{1}_{\{t >t_{\mathbf{b}}\}} e^{-\int_0^{t_{\mathbf{b}}}\Sigma(t-\tau,x-\tau v,v)d\tau}
\Big(\int_0^{t_{\mathbf{b}}}\nabla_x \Sigma(t-\tau,x-\tau v,v)d\tau \Big) \mathbf{R}(t-t_{\mathbf{b}},x_{\mathbf{b}},v)\\
&&\hspace{3mm}+ \int_0^{\min\{t,t_{\mathbf{b}}\}}e^{-\int_0^s\Sigma(t-\tau,x-\tau v,v)d\tau}\nabla_x \mathbf{Q}(t-s,x-sv,v)ds \notag\\
&&\hspace{3mm}- \int_0^{\min\{t,t_{\mathbf{b}}(x,v)\}}e^{-\int_0^s\Sigma(t-\tau,x-\tau v,v)d\tau}\Big(\int_0^s \nabla_x \Sigma(s-\tau,x-\tau v,v)d\tau \Big) \mathbf{Q}(t-s,x-sv,v)ds,\notag\label{34Expression of X Derives}
\end{eqnarray*}}
and
{\small \begin{eqnarray*}
&&\nabla_v u(t,x,v)\mathbf{1}_{\{t\neq t_{\mathbf{b}}\}}\notag \\
&&=\mathbf{1}_{\{t < t_{\mathbf{b}}\}}  e^{-\int_0^t \Sigma(t-\tau,x-\tau v,v)d\tau}[-t\nabla_x u_0+\nabla_v u_0](x-tv,v)\notag\\
&&\hspace{3mm}  - \mathbf{1}_{\{t < t_{\mathbf{b}}\}}  e^{-\int_0^t \Sigma(t-\tau,x-\tau v,v)d\tau}\left\{\int_0^t(-\tau\nabla_x \Sigma +\nabla_v\Sigma)(t-\tau,x-\tau v,v)d\tau\right\} u_0(x-tv,v)\notag\\
&& \hspace{3mm} - \mathbf{1}_{\{t >t_{\mathbf{b}}\}} t_{\mathbf{b}} e^{-\int_0^{t_{\mathbf{b}}} \Sigma(t-\tau,x-\tau v,v)d\tau}\bigg\{\sum_{i=1}^2\tau_i \partial_{\tau_i}\mathbf{R}\notag\\
& &\hspace{2cm}-\frac{n(x_{\mathbf{b}})}{n(x_{\mathbf{b}})\cdot v}\Big\{\Big[\partial_t+\sum_{i=1}^2(v\cdot \tau_i)\partial_{\tau_i}+ \Sigma \Big]\mathbf{R}-\mathbf{Q}\Big\} \bigg\}(t-t_{\mathbf{b}},x_{\mathbf{b}},v)\notag\\
& &\hspace{3mm}  + \mathbf{1}_{\{t >t_{\mathbf{b}}\}} t_{\mathbf{b}} e^{-\int_0^{t_{\mathbf{b}}} \Sigma(t-\tau,x-\tau v,v)d\tau}
 \nabla_v \mathbf{R} (t-t_{\mathbf{b}},x_{\mathbf{b}},v)\\
&&\hspace{3mm} - \mathbf{1}_{\{t >t_{\mathbf{b}}\}} e^{-\int_0^{t_{\mathbf{b}}}\Sigma(t-\tau,x-\tau v,v)d\tau}
\Big(\int_0^{t_{\mathbf{b}}}\{-\tau \nabla_x \Sigma+\nabla_v\Sigma\}(t-\tau,x-\tau v,v)d\tau \Big) \mathbf{R}(t-t_{\mathbf{b}},x_{\mathbf{b}},v)\notag\\
&&\hspace{3mm}  + \int_0^{\min\{t,t_{\mathbf{b}}\}}e^{-\int_0^s\Sigma(t-\tau,x-\tau v,v)d\tau}\{\nabla_v \mathbf{Q}-s\nabla_x \mathbf{Q}\}(t-s,x-sv,v)ds \notag\\
&&\hspace{3mm} - \int_0^{\min\{t,t_{\mathbf{b}}(x,v)\}}e^{-\int_0^s\Sigma(t-\tau,x-\tau v,v)d\tau}\Big(\int_0^s \{-\tau \nabla_x+\nabla_v\Sigma\}(s-\tau,x-\tau v,v)d\tau \Big) \notag\\
& &\hspace{2cm}\times \mathbf{Q}(t-s,x-sv,v)ds.\notag
\label{34Expression of V Derives}
\end{eqnarray*}}
Therefore, we have, for all $0\leq t\leq T$
\begin{eqnarray}
% || \partial_t u(t)\mathbf{1}_{\{t\neq t_{\mathbf{b}}\}}||_1 &\lesssim& ||u_0||_{\infty}+||\nabla_x u_0||_1 + ||\mathbf{Q}(s)(0)||_{\infty} %+\sup_{0\leq s\leq T} ||\mathbf{R}(s)||_{\infty} + \int_0^t ||\partial_t \mathbf{Q} (s)||_1,\label{34Estimate of T Derives}
%\\
||\nabla_x u(t)\mathbf{1}_{\{t\neq t_{\mathbf{b}}\}}||_1 &\lesssim& ||\nabla_x u_0||_1 +(||u_0||_{\infty}+||\mathbf{R}||_{\infty}) \notag\\
&& + \int_0^t \Big|\sum_{i=1}^2\tau_i \partial_{\tau_i}\mathbf{R}-\frac{n(x_{\mathbf{b}})}{n(x_{\mathbf{b}})\cdot v}\{[\partial_t+\sum_{i=1}^2(v\cdot \tau_i)\partial_{\tau_i}+ \Sigma ]\mathbf{R}-\mathbf{Q}\}\Big|_{\gamma_-,1}\\
&& + \int_0^t||\nabla_x \mathbf{Q}(s)||_1+\int_0^ts||\mathbf{Q}(s)||_{\infty},\notag\label{34Estimate od X detives}
\end{eqnarray}
and
\begin{eqnarray}
||\nabla_v u(t)\mathbf{1}_{\{t\neq t_{\mathbf{b}}\}}||_1 &\lesssim& t||\nabla_x u_0||_1+||\nabla_v u_0||_1+t||u_0||_{\infty} \notag\\
&& + t\int_0^t \Big|\sum_{i=1}^2\tau_i \partial_{\tau_i}\mathbf{R}-\frac{n(x_{\mathbf{b}})}{n(x_{\mathbf{b}})\cdot v}\Big\{[\partial_t+\sum_{i=1}^2(v\cdot \tau_i)\partial_{\tau_i}+ \Sigma ]\mathbf{R}-\mathbf{Q}\Big\}\Big|_{\gamma_-,1}\notag\\
&& + \int_0^t|\nabla_v \mathbf{R}|_{\gamma_-,1}+t^2 \sup_{0\leq s\leq t}|\mathbf{R}(s)|_{\gamma_-,\infty}\\
&& + t\int_0^t||\nabla_x \mathbf{Q}(s)||_1+\int_0^t||\nabla_v\mathbf{Q}||_1+ \int_0^ts||\mathbf{Q}(s)||_{\infty}.\notag\label{34Estimate od V detives}
\end{eqnarray}
Since $u_0, ~\mathbf{R}$, and $\mathbf{Q}$ have compact supports and the RHS of (\ref{34Estimate od X detives}) and ({\ref{34Estimate od V detives}}) are bounded. Therefore
\begin{eqnarray*}
\partial u\mathbf{1}_{\{t\neq t_{\mathbf{b}}\}} = [
%\partial_t u \mathbf{1}_{\{t\neq t_{\mathbf{b}}\}},
\nabla_x
 u \mathbf{1}_{\{t\neq t_{\mathbf{b}}\}},\nabla_v u \mathbf{1}_{\{t\neq t_{\mathbf{b}}\}}] \in L^{\infty}([0, T]; L^1(\Omega\times \mathbb{R}^3)).
\end{eqnarray*}
Since $\partial u \equiv 0$ around $t= t_{\mathbf{b}}$, clearly $\partial u\mathbf{1}_{\{t\neq t_{\mathbf{b}}\}}$ is the distributional derivative of $u$. Therefore
$\nabla_x u $ and $\nabla_v u$ lie in $L^{\infty}([0, T]; L^1(\Omega\times \mathbb{R}^3))$, this allows us to apply Lemma \ref{Trace theorem} to compute the traces on the incoming boundary in $L^1([0, T]; L^1(\gamma_-,d\gamma))$ (by taking limits of the flow along the characteristics: see the proof of Proposition 1 in \cite{[GKTT1]} for details). Then, by Green's identity
\ref{Green Identity lemma} we know that $\nabla_x u$ and $\nabla_v u$ lie in $C^0([0, T]; L^1(\Omega\times \mathbb{R}^3))$. Then we get (\ref{34Estimates of x Partial Differential}) and (\ref{34Estimates of v Partial Differential}). $\hfill \square\\$

\noindent{\it \bf Proof of Theorem \ref{1BV regularity fro neutron transport}}. We firstly consider the bound in $||\cdot||_{\infty}$ of the solution for the approximation scheme (\ref{34Approximation sequences}). For a fixed $0 < \varepsilon\ll 1$, it is clear that $(u^{\varepsilon,m})_m$ is a Cauchy series for the norm $\sup_{0\leq t \leq T}||\cdot||_{\infty} $ for fixed $0 < T $ from Theorem \ref{1Existence of the solution}.
More precisely, the sequence $(u^{\varepsilon,m})_m$ satisfy
\begin{eqnarray}\label{34Uniform bound for sequaence}
|| u^{\varepsilon,m}||_{\infty} \lesssim_{T,\Omega,V} ||u_0||_{\infty} +\sup_{0\leq s \leq T}|r(s)|_{\infty}+ \sup_{0\leq s\leq T}||q(s) ||_{\infty}.
\end{eqnarray}
Therefore $u^{\varepsilon,m} \rightarrow u^{\varepsilon}$ up to subsequence for the norm $\sup_{0¡Üt¡ÜT}|| \cdot||_{\infty}$ and $u^{\varepsilon}$ satisfies (\ref{34Approximation sequences}) with both $u^{\varepsilon,m+1}$ and $u^{\varepsilon,m}$ replaced by $u^{\varepsilon}$ by the Green theorem. Since $|\chi_{\varepsilon}| \leq 1$ for $0 < \varepsilon \ll 1$, $\sup_{0\leq t\leq T}|| u^{\varepsilon}||_{\infty}$ is uniformly bounded in $\varepsilon$ for fixed $T$. Therefore $u^{\varepsilon} \rightarrow u$ weak $-*$ up to a subsequence and the limiting function $u$ solves the original neutron transport equation in the sense of distributions.

Secondly, we consider the derivatives of the solution $u^{\varepsilon,m}$ of (\ref{34Approximation sequences}). Recall that
$BV (\Omega\times V)$ has

i) a compactness property: Suppose $g^k \in BV$ and $sup_k ||g^k||_{BV}<\infty $,
then there exists $g \in BV$ with $g^k \rightarrow g$ in $L^1$ up to subsequence,

ii) a lower semicontinuity property: Suppose $g^k \in BV$ and $g^k \rightarrow g$ in $L^1_{ \text{loc}}$
then $||g||_{\tilde{BV}} \leq \liminf_{k \rightarrow \infty}||g^k||_{\tilde{BV}}$.\\

Due to the smooth approximation $\hat{u}_0$ of the initial datum $u_0$ and the cut-off
$\chi_{\varepsilon}$, $u^{\varepsilon,m}$ is smooth by Lemma \ref{34Estimate of approximate sequaence}.
On one hand, by Lemma \ref{34Estimate of approximate sequaence} with $\Sigma \geq 0$,
\begin{eqnarray}
&&||\partial u^{\varepsilon,m+1}(t)||_1 +\int_0^t|\partial u^{\varepsilon,m+1}|_{\gamma_+,1} \notag\\
&& \lesssim ||\hat{u}_0 ||_{\infty}+||\partial \hat{u}_0||_1+ \int_0^t \Big( ||u^{\varepsilon,m+1} ||_{\infty}+||u^{\varepsilon,m}||_{\infty} +|| q||_{BV} \Big) \notag\\
&& \hspace{3mm} +
\int_0^t|\partial u^{\varepsilon,m+1}(s)|_{\gamma_-,1}+ \int_0^t \Big(||\partial u^{\varepsilon,m+1}||_1+ ||\partial u^{\varepsilon,m}||_1\Big), \label{34Estimate of Derive of Sequence}
\end{eqnarray}
where we have used the assumptions that
\begin{eqnarray*}
M_a' = ||\partial \Sigma||_{\infty} <\infty, \quad M'_b=\sup_{x,v}\int_{v'}|\partial f(x,v,v')| dv' <\infty.
\end{eqnarray*}
 From the uniform estimate in (\ref{34Uniform bound for sequaence}), we obtain
\begin{eqnarray}
&&||\partial u^{\varepsilon,m+1}(t)||_1 +\int_0^t|\partial u^{\varepsilon,m+1}|_{\gamma_+,1} \notag\\
&&\lesssim ||u_0||_{BV} + \sup_{0\leq s\leq t}(|r(s)|_{\infty}+ |\partial_t r(s)|_1 + |\partial r(s)|_1 +||q(s)||_{BV})
\notag\\
&& \hspace{3mm} + \int_0^t|\partial u^{\varepsilon,m+1}(s)|_{\gamma_-,1} + t\Big[ \sup_{0\leq s\leq t}||\partial u^{\varepsilon,m+1}(s)||_1 +\sup_{0\leq s\leq t}||\partial u^{\varepsilon,m}(s)||_1 \Big]. \label{34Estimate of Derive of Sequence1}
\end{eqnarray}

On the other hand, by taking derivatives $\partial \in \{\nabla_x, \nabla_v\}$ to (\ref{34Approximation sequences}), we have
\begin{eqnarray*}
&&[\partial_t + v\cdot \nabla_x +\Sigma](\nabla_{x} u^{\varepsilon,m+1}) =- \nabla_{x} \Sigma \, u^{\varepsilon,m+1} + \nabla_{x} (\chi_{\varepsilon}K u^{\varepsilon,m} + \chi_{\varepsilon}\hat{q})\\
&&[\partial_t + v\cdot \nabla_x +\Sigma](\nabla_v u^{\varepsilon,m+1}) = -\Big(\nabla_x + \nabla_v \Sigma\Big) u^{\varepsilon,m+1} + \nabla_v (\chi_{\varepsilon} K u^{\varepsilon,m}+ \chi_{\varepsilon}\hat{q} )\\
&& \partial u^{\varepsilon,m+1}(0,x,v) =\partial \chi_{\varepsilon} \hat{u}_0(x,v)+\chi_{\varepsilon}\partial \hat{u}_0(x,v),
\end{eqnarray*}
and, from (\ref{34Approximation sequences})-(\ref{34Normal derivative at the boundaryF}) as well as (\ref{34Estimate od V detives}), we have, for all $(x,v) \in \gamma_-$,
\begin{eqnarray*}
&&|\partial u^{\varepsilon,m+1}(t, x, v)|
 \lesssim \Big(1+\frac{1}{|n(x)\cdot v|}\Big) \int_{n(x)\cdot v'>0}|\partial u^{\varepsilon,m}(t,x,v')| \{n(x)\cdot v'\} dv'  \notag\\
&& \hspace{26mm}+(1+ \frac{1}{|n(x)\cdot v|})\Big( || u^{\varepsilon,m}||_{\infty} + ||\hat{q}||_{\infty} + |\partial_t r|+|\partial r| \Big).
\end{eqnarray*}
We apply Proposition \ref{33Estimate of sequence} to bound
\begin{eqnarray}
 &&\int_0^t|\partial u^{\varepsilon,m+1}|_{\gamma_-,1}\notag\\
&&\hspace{3mm} \lesssim  O(\delta)\int_0^t |\partial u^{\varepsilon,m-1}|_{\gamma_+,1}
 +C_{\delta}t \sup_{0\leq s\leq t} || \partial u^{\varepsilon,m+1}||_1 +C_{\delta}t  \max_{i=m,m-1}\sup_{0\leq s\leq t} ||\partial u^{\varepsilon,i}(s)||_1 \notag\\
& &\hspace{1cm} + C_{\delta}\Big\{||u_0||_{BV} + \sup_{0\leq s\leq t}(|r(s)|_{\infty}+ |\partial_t r(s)|_1 + |\partial r(s)|_1 +||q(s)||_{BV})  \Big\}
.\label{34Estimate of Derive of Sequence2}
\end{eqnarray}
Finally from (\ref{34Estimate of Derive of Sequence1}) and (\ref{34Estimate of Derive of Sequence2}), choosing
$\delta \ll 1$ and $T_0:= T(u_0)$ is small enough, we have for all
$0 \leq t \leq T_0$
\begin{eqnarray*}
&&||\partial u^{\varepsilon,m+1}(t)||_1 +\int_0^t|\partial u^{\varepsilon,m+1}|_{\gamma,1}\\
&&\hspace{3mm}\leq   \frac{1}{8}\max_{i=m,m-1}\Big\{\sup_{0\leq s\leq t}||\partial u^{\varepsilon,i}(s)||_1 +\int_0^t |\partial u^{\varepsilon,i}|_{\gamma,1}\Big\}\\
&& \hspace{1cm} + C\Big\{||u_0||_{BV} + \sup_{0\leq s\leq t}(|r(s)|_{\infty}+ |\partial_t r(s)|_1 + |\partial r(s)|_1 +||q(s)||_{BV})  \Big\}. \label{34Estimate of Derive of Sequence3}
\end{eqnarray*}
Now, using (\ref{2Main result for sequaences}) with $k=2$, we conclude, for all $m\in \mathds{N}$
\begin{eqnarray}
&&||\partial u^{\varepsilon,m+1}(t)||_1 +\int_0^t|\partial u^{\varepsilon,m+1}|_{\gamma,1} \notag \\
&&\hspace{3mm} \lesssim ||u_0||_{BV} + \sup_{0\leq s\leq t}(|r(s)|_{\infty}+ |\partial_t r(s)|_1 + |\partial r(s)|_1 +||q(s)||_{BV})  . \label{34Main Estimate of BV}
\end{eqnarray}
Now, we pass the to limit in $m \rightarrow \infty$ and then in $\varepsilon \rightarrow 0$ to conclude the main theorem when $0\leq t\leq T_0$. Repeat the same procedure for $[T_0,2T_0],~[2T_0,3T_0]\cdots, $ to conclude the main theorem for all $0\leq t\leq T$. From the compactness and a lower semicontinuity we conclude
\begin{eqnarray*}
\sup_{0\leq s\leq t}||u(s)||_{BV} \lesssim ||u_0||_{BV} + \sup_{0\leq s\leq T}(|r(s)|_{\infty}+ |\partial_t r(s)|_1 + |\partial r(s)|_1 +||q(s)||_{BV}) ,\quad 0\leq t\leq T.
\end{eqnarray*}

For the boundary term we use the weak compactness of measures: If $\sigma^k$ is a signed Radon
measures on $\partial\Omega\times V$ satisfying $sup_k \sigma^k(\partial\Omega\times V ) < \infty$ then there exists a Radon measure $\sigma$
such that $\sigma^k \rightarrow \sigma$ in $\mathfrak{M}$.
More precisely we define, for almost-every $s$, and for any Lebesgue-measurable set $A \subset
\partial\Omega\times V$,
\begin{eqnarray*}
\sigma_s^{\varepsilon,m}(A) &=& \Big(\sigma_{s,x^1}^{\varepsilon,m}(A),\sigma_{s,x^2}^{\varepsilon,m}(A),\sigma_{s,x^3}^{\varepsilon,m}(A),
\sigma_{s,v^1}^{\varepsilon,m}(A),\sigma_{s,v^2}^{\varepsilon,m}(A),\sigma_{s,v^3}^{\varepsilon,m}(A)  \Big)^T\\
&:=& \int_A \nabla_{x,v}u^{\varepsilon,m}(s)d\gamma \in V\times V.
\end{eqnarray*}
Then there exists a Radon measure $\sigma_s$  such that $ \sigma_s^{\varepsilon,m} \rightharpoonup \sigma_s $ in
$\mathcal{M}$, i.e.
\begin{eqnarray}
\int_{\partial \Omega \times V}g \partial u^{\varepsilon,m}d\gamma \rightarrow
\int_{\partial \Omega \times V}g d\sigma_s ~~~~{\text{
for~ all}}~~~ g \in C_c^0(\partial \Omega \times V).  \label{34Convergence of Boundary term}
\end{eqnarray}
It is standard (Hahn¡¯s decomposition theorem) to decompose $\sigma_s = \sigma_{s,+} - \sigma_{s,-}$ with
$\sigma_{s,\pm}\geq  0$.
Denote $|\sigma_s|_{\mathcal{M}(\gamma)} = \sigma_{s,+}(\partial \Omega \times V)
 + \sigma_{s,-}(\partial \Omega \times V)$. Then by the lower semicontinuity property of measures we have \begin{eqnarray*}
 |\sigma_s|_{\mathcal{M}(\gamma)}  \leq  \liminf |\sigma_s^{\varepsilon,m}|_{\mathcal{M}(\gamma)} = \liminf |\partial u_s^{\varepsilon,m}|_{L^1(\gamma)},
 \end{eqnarray*}
 So that by (\ref{34Main Estimate of BV})
\begin{eqnarray*}
\int_0^t |\sigma_s|_{\mathcal{M}(\gamma)}ds \lesssim ||u_0||_{BV}+\sup_{0\leq s\leq t}|\partial r(s)|_1 +\sup_{0\leq s\leq t}||g(s)||_{\infty}).
\end{eqnarray*}
 Due to (\ref{34Convergence of Boundary term}), the (distributional) derivatives $\nabla_{x,v}u(s)|_{\gamma}$
equal the Radon measure $\sigma_s$ on $\partial \Omega \times \mathbb{R}^3$ in the sense of distributions. $\hfill \square$

\section{Appendix}
\setcounter{section}{4} \setcounter{equation}{0}
\begin{lem}\label{Derivatives} (\cite{[Guo2],[GKTT1]}).
If $v\cdot n(x_b(x,v))<0$, then $(t_b(x,v),x_b(x,v))$ are smooth functions of $x,v$ such that
\begin{eqnarray*}
&&\nabla_x t_b = \frac{n(x_b)}{v\cdot n(x_b)}, ~~~~~~~~~~~~\nabla_v t_b=-\frac{t_b n(x_b)}{v\cdot n(x_b)},\\
&& \nabla_x x_b=I-\frac{n(x_b)}{ v\cdot n(x_b)} \otimes v,~~~\nabla_v x_b= - t_b I +\frac{t_b n(x_b)}{v\cdot n(x_b(x,v))}\otimes v.
\end{eqnarray*}
Let $x_i \in \partial \Omega$ for $i=1,2$ and $(t_1,x_1,v)$ and $(t_2,x_2,v)$ be connected with the trajectory $\frac{dX(s)}{ds}=V(s)$, $\frac{dV(s)}{ds}=0$ which lies inside $\overline{\Omega}$. Then there exists a constant $C_{\xi}>0$ such that
\begin{eqnarray*}
|t_1-t_2| \geq \frac{n(x_1)\cdot v}{C_{\xi}|v|^2}.
\end{eqnarray*}
\end{lem}
For the estimate of the outing trace on $\gamma_+ \setminus \gamma_+^{\delta}$, we need the following trace theorem.
\begin{lem}\label{Trace theorem} (Outgoing trace theorem,\cite{[EGKM]}). Assume $\psi \geq 0$. For any small parameter $\delta>0$, there exists a constant $C_{\delta,T,\Omega}$ such that for any $h \in L^1([0,T]\times \Omega\times V)$ with $ (\partial_t  +v\cdot \nabla_x + \psi) h$ lying in $L^1([0,T] \times \Omega\times V)$, we have for all $0\leq t\leq T$,
\begin{eqnarray}
\int_0^T \int_{\gamma_+ \setminus \gamma_+^{\delta}}|h|d\gamma da \leq C_{\delta,T,\Omega}\bigg[||h_0||_1 +\int_0^t \bigg( ||h(s)||_1+||[\partial_t +v\cdot \nabla_x + \psi ]h(s)||_1\bigg) ds\bigg].
\end{eqnarray}
Furthermore, for any $(s,x,v)\in [0,T]\times \Omega\times V$, the function $h(s+s',x+s'v,v)$ is absolutely continuous in $s'$ in the interval $[-\min\{t_b(x,v),s\},\min\{t_b(x,-v),T-s\}] $.
\end{lem}
\begin{lem}\label{Green Identity lemma} (Green Identity, \cite{[Guo2]}, \cite{[GKTT1]}).
Let $p\in [1,\infty)$. Assume that $f, (\partial_t  +v\cdot \nabla_x + \psi) f \in L^p([0,T]\times \Omega\times V)$ with $\psi \geq 0$ and $f|_{\gamma_-} \in L^p([0,T]\times \partial\Omega \times V, dt d\gamma)$. Then $f\in C^0([0,T],L^p(\Omega\times V))$ and $f|_{\gamma_+} \in L^p([0,T]\times \Omega\times V,dtd\gamma)$ and for almost every $t\in [0,T]$,
\begin{eqnarray*}\label{Green Identity}
||f||_p^p + \int_0^t |f|^p_{\gamma_+,p} = ||f(0)||_p^p +\int_0^t|f|^p_{\gamma_-,p} + \int_0^t\iint_{\Omega \times V}\{\partial_t f  +v\cdot \nabla_x f+ \psi f\}|f|^{p-2}f.
\end{eqnarray*}
\end{lem}

The covering lemma has proved in \cite{[Guo2],[GKTT2]}, here we have the similar result by replacing $B_N = \{v\in \mathbb{R}^3: |v|\leq N\}$ with the compact set $V \in \mathbb{R}^3$.
\begin{lem}\label{Guo Covering} (Covering Lemma, \cite{[Guo2],[GKTT2]}). Let $\Omega \subset R^3$ be an open bounded set with
a smooth boundary $\partial\Omega$ and $V$ is a compact set in $\mathbb{R}^3$. Then, for all $x\in \overline{\Omega}$, we have
\begin{eqnarray}\label{Measure of singular set}
m_3 \{v\in V:~ n(x_b(x,v))\cdot v=0 \}=0.
\end{eqnarray}
Moreover, for any $\varepsilon$, there exist $\delta_{\varepsilon,V}>0$ and $l=l_{\varepsilon,\Omega,V}$ balls
$B(x_1,r_1)$, $B(x_2,r_2)$, $\cdots ,B(x_l,r_l)$ with $x_i \in \overline{\Omega}$ and covering $\overline{\Omega}$
(i.e. $\overline{\Omega} \subset \bigcup B(x_i,r_i)$), as well as $l$ open sets $\mathcal{O}_{x_1},\mathcal{O}_{x_2},\cdots
\mathcal{O}_{x_l} \subset V $, with $m_3(\mathcal{O}_{x_i}) <\varepsilon$ for all $1\leq i\leq l_{\varepsilon,\Omega,V}$ such that  for any $x\in \overline{\Omega}$, there exists $i=1,2,\cdots,l_{\varepsilon,\Omega,V} $ such that
$x\in B(x_i,r_i)$ and
\begin{eqnarray}
|v\cdot n(x_b(x,v))| > \delta_{\varepsilon,V}, ~~~\text{for all}~~ v ~\notin \mathcal{O}_{x_i}.
\end{eqnarray}
In particular,
\begin{eqnarray}\label{Covering}
\bigcup_{x\in B(x_i,r_i)}\bigg\{ v\in V: |v\cdot n(x_b(x,v))| \leq \delta_{\varepsilon,V}\bigg\} \subset \mathcal{O}_{x_i}.
\end{eqnarray}
\end{lem}

Acknowledgments. The second author would like to express deep thank to S.Q. Liu and F.J. Zhou for the fruitful discussion during the visit to Brown University.

\end{document}